\documentclass[11pt]{article}

\usepackage{amsmath,amssymb,amsfonts,xcolor}
\usepackage{epsfig}

\numberwithin{equation}{section}
\newtheorem{thm}{Theorem}[section]

\newtheorem{alem}[thm]{Lemma}
\newtheorem{aprop}[thm]{Proposition}
\newtheorem{acor}[thm]{Corollary}
\newtheorem{arem}[thm]{Remark}

\newtheorem{algo}{Algorithm}

\newenvironment{adem}[1][]%
   {\ \\ {\bf Proof#1: }}%
   {\hfill\mbox{\rule{2 true mm}{3 true mm}}\vskip 2 ex\noindent}
\newcommand{\rat}{\rho}
\newcommand{\eqdef}{:=}
\def\eqsp{\;}
\newcommand{\un}{\mathbf{1}}
\newcommand{\hatH}{\widehat{H}}
\def\tv{\mathrm{TV}} 
\def\calK{\mathcal{K}}
\def\wl{\mathrm{WL}}
\def\wt{\mathrm{wt}}

\newcommand{\E}{{\mathbb E}}

\newcommand{\N}{{\mathbb N}}
\newcommand{\R}{{\mathbb R}}

\renewcommand{\P}{{\mathbb P}}
\renewcommand{\t}{\theta}
\def\Xset{\mathsf{X}}
\def\tn{\theta}   
 
\def\tu{\tilde \theta}  
\def\tstar{\theta_\star} 
\newcounter{hypoconbis}
\newcounter{saveconbis}
\newcommand\debutA{\begin{list} {\textbf{A\arabic{hypoconbis}}}{\usecounter{hypoconbis}}\setcounter{hypoconbis}{\value{saveconbis}}}
\newcommand\finA{\end{list}\setcounter{saveconbis}{\value{hypoconbis}}}
\def\rmd{\mathrm{d}}
\newcounter{hypoconbisf}
\newcounter{saveconbisf}
\newcommand\debutF{\begin{list} {\textbf{R\arabic{hypoconbisf}}}{\usecounter{hypoconbisf}}\setcounter{hypoconbisf}{\value{saveconbisf}}}
\newcommand\finF{\end{list}\setcounter{saveconbisf}{\value{hypoconbisf}}}
\def\rmd{\mathrm{d}}

\newcommand\debutS{\begin{list} {\textbf{S\arabic{hypoconbisf}}}{\usecounter{hypoconbisf}}\setcounter{hypoconbisf}{\value{saveconbisf}}}
\newcommand\finS{\end{list}\setcounter{saveconbisf}{\value{hypoconbisf}}}

\renewcommand{\leq}{\leqslant}
\renewcommand{\le}{\leqslant}
\renewcommand{\geq}{\geqslant}
\renewcommand{\ge}{\geqslant}
\newcommand{\dps}{\displaystyle}

\title{Convergence and efficiency of adaptive importance sampling techniques with partial biasing}

\author{G. Fort\thanks{LTCI, CNRS, T\'el\'ecom ParisTech,
    Universit\'e Paris-Saclay, 75013, Paris, France; IMT, Universit\'e de Toulouse et CNRS, F-31062 Toulouse Cedex 9, France. email:
     gersende.fort@math.univ-toulouse.fr},
  B.Jourdain\thanks{Universit\'e Paris-Est, CERMICS (ENPC), INRIA
  F-77455 Marne-la-Vall\'ee, 
  France. emails: benjamin.jourdain@enpc.fr, tony.lelievre@enpc.fr, gabriel.stoltz@enpc.fr}, T. Leli\`evre\footnotemark[2] and
  G. Stoltz\footnotemark[2]}

\begin{document}

\maketitle

\abstract{We consider a generalization of the discrete-time Self Healing
  Umbrella Sampling method, which is an adaptive importance technique useful to
  sample multimodal target distributions. The importance function is based on
  the weights (namely the relative probabilities) of disjoint sets which form a partition of the space. These
  weights are unknown but are learnt on the fly yielding an adaptive algorithm.
  In the context of computational statistical physics, the logarithm of these
  weights is, up to a multiplicative constant, the free energy, and the
  discrete valued function defining the partition is called the
  collective variable. The algorithm falls into the general class of Wang-Landau type methods, and is a generalization of
  the original Self Healing Umbrella Sampling method in two ways: (i) the
  updating strategy leads to a larger penalization strength of already visited
  sets in order to escape more quickly from metastable states, 
  and (ii) the target distribution is biased using only a fraction of the
  free energy, in order to increase the effective sample size and reduce the
  variance of importance sampling estimators. The algorithm can also be seen as
  a generalization of well-tempered metadynamics. We prove the convergence of
  the algorithm and analyze numerically its efficiency on a toy example.}

\section{Introduction}

In many situations, sampling methods are considered in order to compute
expectations of given observables with respect to a distribution $\pi
\, \rmd \lambda$ with support a subset $\Xset$ of $\R^D$.  We denote by $\pi:\R^D \to \R$
the density of the target distribution, with respect to a reference non
negative measure $\lambda$ on $\R^D$. We are interested here in the case
when the target distribution is highly multimodal.

Typically, the expectations under consideration are approximated by empirical averages
of the observables computed along a path of a Markov chain or process,
ergodic with respect to $\pi \, \rmd \lambda$. This is the principle
of Markov chain Monte Carlo (MCMC) methods (see
e.g.~\cite{brooks:gelman:jones:meng:2011}); the famous
Metropolis-Hastings algorithm is one instance of a general approach to
build a Markov chain having a given probability measure as invariant
distribution~\cite{MRRTT53,Hastings70}.  In many
situations of interest however, the target probability measure $\pi \,
\rmd\lambda$ is multimodal: the most likely regions are separated by
low probability regions, which makes the design of efficient numerical
sampling methods difficult. Indeed, without {\em a priori} knowledge
on $\pi \, \rmd \lambda$, MCMC methods typically rely on local moves
and the algorithms are stuck in high probability regions: the
dynamics of the Markov process is metastable in the sense that it
remains trapped for a very long time in some region of the space,
called a metastable state, before hopping to another metastable state.
The aim of this paper is to
discuss free energy-based adaptive importance sampling techniques
which have been developed in the physics and chemistry literature 
to efficiently sample such a multimodal probability measure in high dimension.

In free energy-based importance sampling techniques, the auxiliary
distribution from which the samples are drawn is obtained by a local
reweighting of the target distribution. More precisely, assume we are
given a partition of the state space $\Xset$ into $d$ disjoint subsets
(called strata hereafter):
$$\Xset=\bigcup_{i=1}^d \Xset_i.$$
For future reference, let us introduce the so-called collective variable $I:\Xset \to \{1,\ldots,d\}$ associated with this partition:
\begin{equation}
\label{eq:def_I}
\forall x \in \Xset, \, I(x)\eqdef i \ \ \text{if and only if}  \ \ x \in \Xset_i.
\end{equation}
Denote by $\tn_\star = ( \tn_\star(1), \ldots, \tn_\star(d))$ the
vector collecting the weights of each stratum under $\pi \, \rmd
\lambda$:
\begin{equation}\label{eq:theta_star}
  \forall i \in \{1, \ldots, d\}, \, \tn_\star(i) \eqdef \int_{\Xset_i} \pi
  \, \rmd\lambda.
\end{equation}
In the context of computational statistical physics, minus the
logarithm of these weights is called the free energy. More generally,
the free energy is, up to a multiplicative constant, the log marginal of the target distribution along some chosen degrees
of freedom, see Section~\ref{sec:WT} for a more precise definition. Finally, for all $\tn \in \Theta$, where
\[
\Theta \eqdef \left\{\tn = \big(\tn(1), \ldots, \tn(d)\big) \in (0,1)^d, \ \sum_{i=1}^d \tn(i)=1 \right\},
\]
define a probability
measure $\pi_\tn \, \rmd \lambda$ on $\Xset$ by:
\begin{equation}
  \label{eq:def:pitheta}
  \forall x \in \Xset, \, \pi_\tn(x) \eqdef \frac{1}{Z_\tn} \sum_{i=1}^d
  \frac{ \pi(x)}{\theta(i)}\un_{\Xset_i}(x), \qquad \text{with}  \ Z_\tn \eqdef \sum_{i=1}^d \frac{\tn_\star(i)}{\tn(i)}.
\end{equation}
By definition, all the strata have the same weight $1/d$ under
$\pi_{\tn_\star} \, \rmd \lambda$. As a consequence, if the strata are well chosen,
$\pi_{\tn_\star} \, \rmd \lambda$ is less multimodal than the original
target $\pi \, \rmd \lambda$, and the sampling of $\pi_{\tn_\star} \,
\rmd \lambda$ is thus easier: a Markov chain sampling
$\pi_{\tn_\star} \, \rmd \lambda$ easily visits the whole space $\Xset$. Since we are interested in efficient Monte Carlo approximations
of expectations under the distribution $\pi \, \rmd \lambda$, a standard
reweighting (or importance sampling) strategy provides an estimator of such an expectation from samples
approximating $\pi_{\tn_\star} \, \rmd \lambda$, upon noting that we have
for any measurable and bounded function $f:\Xset \to \R$,
\begin{equation}\label{eq:IS}
  \int_{\Xset} f \, \pi \, \rmd \lambda = d \sum_{i=1}^d  \tn_\star(i)  \int_{\Xset_i} 
  f   \, \pi_{\theta_\star} \, \rmd \lambda.
\end{equation}
Notice that for $x \in \Xset_i$, the importance ratio $\pi(x) /
\pi_{\tn_\star}(x)$ is equal to $d \tn_\star(i)$ which justifies formula~\eqref{eq:IS}.

Before discussing how to apply this method for a fixed partition, let
us explain how the partition can be built. Free-energy based techniques have
originally been designed in the field of Monte Carlo
simulation of materials and molecular dynamics. In this context a continuous collective variable $\xi: \Xset \to \R$
is chosen, and the partitions are designed as level sets of $\xi$ (typically,
$\Xset_i = \xi^{-1}([a_i,a_{i+1}))$ for $a_1 < a_2 < \ldots <a_{d+1}$). The
choice of a good function $\xi$ is the subject of many papers: it
typically relies on some {\em a priori} knowledge of some ``slow'' degrees of
freedom, which index transitions between metastable states. We refer for
example to the
monographs~\cite{lelievre-rousset-stoltz-book-10,chipot-pohorille-07} for
discussions on this subject. For applications to Bayesian inference, the choice
of a good partition is discussed in~\cite{chopin-lelievre-stoltz-12}. From now
on, we assume that the partition is given.

In practice, besides the choice of the partition, there are two major
difficulties when applying the free-energy biased sampling method described above: first, the vector $\tn_\star$ is unknown;
second, the discrepancy between the weights $(\tn_\star(i))_{i=1,\ldots,d}$ may
yield a large variance in the importance sampling estimator of the quantity
$\int_{\Xset} f \pi \rmd \lambda$ deduced from \eqref{eq:IS}.  This discrepancy between the reweighting factors can be
quantified through the so-called effective sample size (see~\cite{KongLiuWong}
and formula~\eqref{eq:ESS} below for a precise definition): the larger the
discrepancy between the weights, the smaller the effective sample size.


To overcome the first difficulty, namely the fact that the vector $\tn_\star$
is unknown, the idea is to learn it on-the-fly. This yields a so-called
adaptive importance sampling algorithm.  The sampler is an iterative procedure
and each iteration combines a sampling step and an update step: the sampling
step samples a configuration $X_{n+1}$ under a distribution approximating
$\pi_{\tn_n} \, \rmd \lambda$; the update step builds a new approximation
$\tn_{n+1}$ of $\tn_\star$, by using the past of the algorithm $\{\tn_0, X_0,
\ldots, X_n \}$ and the new draw $X_{n+1}$. These two steps are designed in
such a way that (in some sense to be made precise) $\tn_n$ converges to
$\tn_\star$ in the longtime limit $n \to +\infty$ and thus, the distribution of
$X_n$ converges to $\pi_{\tn_\star} \, \rmd \lambda$.  Many free energy-based
adaptive importance sampling techniques have been proposed in the statistical
physics literature, first with the scope of computing the vector $\tn_\star$. 
As a byproduct, they also provide a
sampler targeting the distribution $\pi_{\tn_\star} \, \rmd \lambda$. These
algorithms essentially differ in the way the updating strategy is implemented,
see e.g.  the Wang Landau algorithm~\cite{wang-landau-01,wang-landau-01-PRL},
the Self-Healing Umbrella Sampling (SHUS) algorithm~\cite{SHUS}, the Adaptive
Biasing Force algorithm~\cite{DP01,HC04,JLR10,minoukadeh-chipot-lelievre-10},
the metadynamics algorithm~\cite{laio-parrinello-02,bussi-laio-parrinello-06},
the well-tempered metadynamics algorithm~\cite{barducci-bussi-parrinello-08},
etc...

To overcome the second difficulty, we apply the adaptive strategy described
above, but for another target than the free-energy biased target density
$\pi_{\theta_\star}$. More precisely, we consider a density
$\pi^\rho_{\theta_\star}$ where a measurable non-decreasing function
$\rat:(0,1) \to (0,+\infty)$ is introduced in order to make
$\pi^\rho_{\theta_\star}$ and the original target $\pi$ closer, and thus to
lower the discrepancy of the weights in the importance sampling estimator. For
all $\tn \in \Theta$, the probability measure $\pi^\rat_{\theta}$ is defined by
\begin{equation}\label{defpithetaf}
  \forall x \in \Xset, \ \    \pi^\rat_{\theta}(x) \eqdef  \left(Z^\rat_\theta \right)^{-1} \ \sum_{i=1}^d \frac{\pi(x)}{ \rat(\theta(i))}\un_{\Xset_i}(x), \qquad \text{with} \  Z^\rat_\theta \eqdef \sum_{i=1}^d\frac{\theta_\star(i)}{
    \rat(\theta(i))}.
\end{equation}
The idea
of using the biased measure~\eqref{defpithetaf} instead of~\eqref{eq:def:pitheta} 
dates back to the well-tempered metadynamics
algorithm~\cite{barducci-bussi-parrinello-08}, where the function
$\rat$ is chosen as $\rat: t \mapsto t^a$ for some $a \in (0,1)$; in some sense, we are
studying here a discrete variant (discrete in time and discrete in
terms of the collective variable) of the
well-tempered metadynamics sampler, see Section~\ref{sec:WT} below for a more
detailed discussion. 

Since the biasing measure changes, importance sampling estimators such as~\eqref{eq:IS} should be modified accordingly: for any measurable and bounded function
$f:\Xset \to \R$,
\begin{equation}\label{eq:IS_rho}
  \int_{\Xset} f \, \pi \, \rmd\lambda = \left(\sum_{j=1}^d \frac{\tn_\star(j)}{\rat(\tn_\star(j))}  \right) \sum_{i=1}^d \rat(\tn_\star(i)) \,  \int_{\Xset_i}  f \,   \pi^\rat_{\tn_\star} \, \rmd \lambda,
\end{equation}
from which an importance sampling estimator of $\int_{\Xset} f \, \pi \, \rmd
\lambda$ can be deduced from samples approximating $\pi_{\tn_\star}^\rat \,
\rmd \lambda$. The aim of the function $\rat$ is to make the
importance ratio $\pi/\pi^\rat_{\tn_\star}$ (and thus the weights $\rho(\tn_\star(i))$) closer to a constant. More precisely, $\rho$ should be such that
\[
\frac{\dps \max_{i=1,\dots,d} \rho(\tn_\star(i))}{\dps \min_{i=1,\dots,d} \rho(\tn_\star(i))} \leq \frac{\dps \max_{i=1,\dots,d} \tn_\star(i)}{\dps \min_{i=1,\dots,d} \tn_\star(i)}.
\]
This allows to spend more time in strata~$\Xset_i$ with larger weights
$\tn_\star(i)$ in the estimation of averages with respect to~$\pi$. To better
understand the interest of the function $\rat$, consider the example $\rat: t
\mapsto t^a$ for some $a \in (0,1]$. In that case,
$\pi_{\theta_\star}^{t^a}(\Xset_i)$ is proportional to $\left(\tn_{\star}(i)
\right)^{1-a}$, and for $x \in \Xset_i$, the importance ratio
$\pi(x)/\pi^\rat_{\tn_\star}(x)$ is equal to
$\left( \sum_{i=1}^d \left(\tn_\star(i)\right)^{1-a}\right) \left(\tn_\star(i)\right)^a $. On the one hand, the closer $a$ is to $1$, the more uniform the
weights of the strata are and the less metastable the sampler targeting
$\pi_{\tn_\star}^{t^a}$ is (at least if the strata are well chosen, see the
discussion above). But when $a$ is close to $1$, the importance ratio
is far from a constant, and thus the effective sample size associated
with the estimator~\eqref{eq:IS_rho} is small.
On the other hand, when $a$ gets close to $0$, the function
$\rat$ is close to a constant, and thus the effective sample size associated
with the estimator~\eqref{eq:IS_rho} is large.  But when $a$ gets close to $0$,
$\pi_\theta^{t^a}$ gets close to the original target density $\pi$, and thus
the sampling dynamics becomes as metastable as the original non adaptive one.  There
is thus a compromise to find between two contradictory objectives: biasing the
dynamics in order to leave the metastable states more quickly and thus converge
faster to equilibrium; not modifying the original target probability $\pi$ too much,
since this will give too large weights to originally unlikely regions, which
will lead to a very small effective sample size.

\bigskip

The first contribution of this paper is to propose a free-energy based adaptive
importance sampling algorithm, denoted SHUS$_\rat^\alpha$, for the sampling of
a metastable distribution $\pi \, \rmd \lambda$, which combines the two
ingredients presented above. The lowerscript $\rho$ refers to the function
$\rho$ discussed above, while the parameter $\alpha \in (1/2,1]$ enters the
updating formula of the sequence $(\theta_n)_{n \ge 0}$ in order to control its
rate of convergence to $\tn_\star$. This method is designed to \textit{(i)} 
learn on-the-fly the weights $(\tn_\star(i))_{i=1, \ldots,d}$ of the strata,
and \textit{(ii)} provide draws sampling $\pi_{\tn_\star}^\rat$. This
algorithm is described in Section~\ref{sec:methodo}, where its relationships
with the SHUS algorithm and the Wang-Landau algorithm are also discussed. We
show that, like all these algorithms, SHUS$_\rat^\alpha$ updates the weight
vector $\tn_n$ based on the frequency of visit of each stratum in such a way
that it penalizes already visited strata when sampling the next configuration
$X_{n+1}$. As in SHUS and well-tempered metadynamics, SHUS$_\rat^\alpha$
automatically computes, based on its past behavior, the strength of the
penalization.  As in well-tempered metadynamics, SHUS$_\rat^\alpha$ also allows
the use of a function $\rho$ in order to improve the
quality of the importance sampling estimators based on the algorithm, as discussed above. In
particular, we explain that, for $\alpha=1$, SHUS$_\rat^\alpha$ with
$\rho(t)=t$ is the standard SHUS algorithm (see Section~\ref{sec:relation_FE})
while, with $\rho(t)=t^a$ (where $a \in (0,1)$), it can be seen as a version of
the well-tempered metadynamics algorithm with a discrete collective variable
and a discrete-in-time stochastic dynamics (see
Section~\ref{sec:WT}). Nevertheless, it differs from these samplers by introducing the additional
degree of freedom~$\alpha$, which can be tuned so that the algorithm escapes
far more quickly from metastable states. This design parameter therefore
improves the transient phase of the algorithm. One motivation of this work is actually the study
of the convergence and efficiency of the well-tempered metadynamics, in a
slightly different setting than the original one, and to propose and study accelerated
versions of the SHUS algorithm and well-tempered metadynamics, thanks to the
introduction of the parameter $\alpha$.


The second contribution of this paper is to mathematically analyze the
asymptotic behavior of the SHUS$_\rat^\alpha$ algorithm.  Our work belongs to a
series of contributions where free energy-based adaptive importance algorithms
are mathematically analyzed in order to prove their convergence and to measure
their efficiency; see
e.g.~\cite{atchade-liu-10,JR11,fort:jourdain:kuhn:lelievre:stoltz:2014,amrx}
for Wang-Landau algorithms, \cite{LRS08,lelievre-minoukadeh-09} for ABF,
\cite{fort:jourdain:lelievre:stoltz:2015} for the SHUS algorithm and
\cite{dama-parrinello-voth-14} for the well-tempered metadynamics algorithm.

We provide in Section~\ref{sec:convergence:results} sufficient conditions on
the function $\rat$, on the parameter $\alpha$ and on the sampling
step, in order to obtain the convergence of $(\tn_n)_{n \ge 0}$ to
$\tn_\star$ and the
consistency of an importance sampling estimator of $\int_{\Xset} f \, \pi \,
\rmd \lambda$ computed "online" (i.e. from the points $(X_n)_{n \ge 0}$ and the
sequence $(\tn_n)_{n \ge0}$ produced by the iterative algorithm). For that
purpose, a crucial step is to provide a recurrence result for the random sequence $(\tn_n)_{n \ge
  0}$ showing that, with probability one, it enters infinitely often a compact
subset of $\Theta$.  Finally, we show that the update rule for the
vectors of weights
$\tn_n$ can be seen as a stochastic approximation scheme with a (random)
step-size sequence, self-tuned by the algorithm; we prove that this sequence
converges to zero at the rate $O(n^{-\alpha})$.

After discussing in details the link between SHUS$^\alpha_\rho$ and the
well-tempered metadynamics in Section~\ref{sec:WT}, we then numerically
illustrate the efficiency of this algorithm on a toy model in
Section~\ref{sec:application}. The roles of the function $\rat$ and of the
parameter $\alpha$ are highlighted. 
Finally, Section~\ref{sec:proofs} is devoted to the proofs of the asymptotic
results stated in Section~\ref{sec:convergence:results}.

\section{The SHUS$_\rat^\alpha$ algorithm}
\label{sec:methodo}
We introduce the SHUS$_\rat^\alpha$ algorithm in Section~\ref{sec:algo}. We
then discuss its connections with the well-known Wang-Landau algorithm in
Section~\ref{sec:relation_FE}, where we also compare the SHUS$_\rat^\alpha$
algorithm with other free energy adaptive methods. Finally, in
Section~\ref{sec:discuss_algo}, we present how this algorithm
is derived and explain its expected properties using heuristic arguments. We would like to stress that
the arguments used in this section are not intended to be fully rigorous, but
hopefully give some intuition on the SHUS$_\rat^\alpha$  algorithm. Rigorous statements about
the convergence of the algorithm are provided in
Section~\ref{sec:convergence:results}.

\subsection{The algorithm}\label{sec:algo}
Let $\rat:(0,1) \to \R_+$ be a measurable and non-decreasing function.
Let $\alpha
\in (1/2,1]$, $\gamma >0$ and $\mu>0$ be three constants. Define the function $g_\alpha:(0,+\infty)\to
(0,+\infty)$ by:
\begin{equation}\label{hyp:g}
\forall s >0, \, g_\alpha(s) 
= \left\{
\begin{array}{ll}
\big(\ln(1+s)\big)^{\frac{\alpha}{1-\alpha}} & \text{if} \ \alpha \in \left(\frac{1}{2},1\right), \\[5pt]
s^\mu& \text{if} \ \alpha =1.
 \end{array}
\right.
\end{equation}
For any measurable function $\rat:(0,1) \to \R_+$ and for any $\tn \in \Theta$,
let us denote by $P^\rat_\tn$ a Markov transition kernel ergodic with respect
to the probability measure $\pi_\tn^\rat \, \rmd \lambda$, where $\pi_\tn^\rat$
is given by \eqref{defpithetaf}. For example, $P^\rat_\tn$ stands for a
Metropolis Hastings
kernel~\cite{brooks:gelman:jones:meng:2011,Hastings70,MRRTT53}. This will
actually be our setting in the sequel.

\begin{algo}[SHUS$_\rat^\alpha$]\label{algo:GenSHUS}
  For a (possibly random) initial condition $(\tu_0,X_0)$ in
  $(\R_+^*)^d\times\Xset$, the SHUS$_\rat^\alpha$ algorithm consists in
  iterating the following three steps over $n \ge 0$: given $(\tu_n,X_n)\in
  (\R_+^*)^d\times \Xset$,
\begin{itemize}
\item Compute the normalizing constant $S_n$ and the probability measure
  $\tn_n$ on $\{1,\ldots,d\}$, obtained by normalizing $\tu_n$:
     \begin{equation}
       \label{eq:def:theta}
      S_n \eqdef \sum_{i=1}^d \tu_n(i), \qquad  \tn_n \eqdef \frac{\tu_n}{S_n} \in \Theta,
     \end{equation}
   \item Sample $X_{n+1}$ according to the distribution $P^\rat_{\tn_n}(X_{n},\cdot)$,
   \item Compute, for all $i \in \{1, \ldots d\}$,
  \begin{equation}
    \label{eq:def:thetatilde}
  \tu_{n+1}(i) = \tu_n(i) +  \frac{\gamma}{g_\alpha(S_n)} \, S_n  \, \rat(\t_n(i))\,  \un_{\Xset_i}(X_{n+1}).
 \end{equation}
\end{itemize}
\end{algo}
Roughly speaking, $\tu_n(i)$ stands for an "occupation measure" of the stratum
$i$ at the end of iteration~$n$. The weight vector
$\tn_n$ is the normalized vector associated with this occupation
measure. If the sample $X_{n+1}$ is in $\Xset_{i_0}$, then
$\tu_{n+1}(i_0)>\tu_{n}(i_0)$ while for $i \neq i_0$, $\tu_{n+1}(i) =
\tu_{n}(i)$. Therefore, at the next iteration, the probability to be in
the $i_0$-th stratum is lower under the probability
$\pi^{\rho}_{\theta_{n+1}} \, \rmd \lambda$ (which is the invariant probability of the
kernel $P^{\rho}_{\theta_{n+1}}$) than under~$\pi^{\rho}_{\theta_{n}} \, \rmd \lambda$.

As will become clear below, the main parameters of the
SHUS$_\rat^\alpha$ algorithm are the parameter $\alpha \in (1/2,1]$
and the function $\rho$.
In particular, we do not explicitly mention the dependence on $\gamma$ in the notation
SHUS$_\rat^\alpha$ since this parameter does not play an important role in the
mathematical analysis. However, this parameter will play a role in the
numerical tests in Section~\ref{sec:application}, when studying the efficiency
of the algorithm as $\alpha$ varies (see the choice
in~\eqref{eq:choice_alpha} which allows to obtain a continuous behavior of the
algorithm in the limit $\alpha \to 1$, in spite of the discontinuity of $\alpha
\mapsto g_\alpha$ at $\alpha=1$).
Likewise, we do not make explicit in the notation
SHUS$_\rat^\alpha$ and $g_\alpha$ the dependence
on the parameter $\mu$ when $\alpha=1$ because, compared to $\alpha$, this
parameter has a weak influence on the behavior of the algorithm (see for
instance Corollary \ref{cor:asymptotic} below).  Nevertheless, this parameter
is needed to enclose in our analysis the discrete version of the well-tempered
metadynamics algorithm, as explained in Section \ref{sec:WT}.

\subsection{Relationship with other free energy adaptive
  techniques}\label{sec:relation_FE}

In this section, we observe that the SHUS$^\alpha_\rho$ algorithm can be
seen as one example of a generalized Wang-Landau 
algorithm WL$_\rho$. This is
useful to understand the basic principles underlying the algorithm, and to
discuss the differences and similarities of SHUS$^\alpha_\rho$ with other free
energy adaptive techniques. In addition, we state convergence results for generalized Wang-Landau 
algorithms (see Propositions~\ref{prop:cvggal} and~\ref{prop:recurrence} below), 
which allows us to prove the convergence of SHUS$_\rat^\alpha$. 

\paragraph{A generalized Wang-Landau algorithm.} Let us first
introduce a generalization of the original Wang-Landau algorithm~\cite{wang-landau-01}.
Let $\rat:(0,1) \to \R_+$ be a measurable and non-decreasing function.
\begin{algo}[WL$_\rat$]\label{algo:WL}
For a (possibly random) sequence $(\gamma_n)_{n \ge 1}$ and
  initial conditions $(\tu_0,X_0)$ in
$(\R_+^*)^d\times\Xset$, the WL$_\rat$ algorithm consists in iterating the following
three steps over $n \ge 0$: given $(\tu_n,X_n)\in (\R_+^*)^d\times \Xset$,
\begin{itemize}
\item Compute the normalizing constant $S_n$ and the probability measure
  $\tn_n$ on $\{1,\ldots,d\}$:
     \begin{equation}
       \label{eq:def:theta_WL}
      S_n = \sum_{i=1}^d \tu_n(i), \qquad  \tn_n = \frac{\tu_n}{S_n} \in \Theta.
     \end{equation}
   \item Sample $X_{n+1}$ according to the distribution
     $P^\rat_{\tn_n}(X_{n},\cdot)$.
   \item Compute, for all $i \in \{1, \ldots d\}$,
  \begin{equation}
    \label{eq:evol_algo_gen}
\tu_{n+1}(i)=\tu_n(i) \left(1 + \gamma_{n+1}  \, \frac{\rho(\tn_n(i))}{\tn_n(i)}
  \un_{\Xset_i}(X_{n+1}) \right).
 \end{equation}
\end{itemize}
\end{algo}
The principle of generalized Wang-Landau algorithms is to penalize
already visited strata in order to favor transitions to new regions of the state space.
The weight $\frac{\rat(\t_n(i))}{\t_n(i)}$ in~\eqref{eq:evol_algo_gen} is introduced in order to
compensate for the biasing term $\frac{1}{\rho(\theta(i))}$ in
$\pi^\rho_\theta$ (see~\eqref{eq:IS_rho}), so
that $(\t_n)_{n \ge 0}$ is expected to converge to $\t_\star$
(see~\eqref{eq:theta_star}), and the stationary state of the algorithm
is expected to be $\pi^\rho_{\theta_\star}$. This will be proven below (see
Propositions~\ref{prop:cvggal} and~\ref{prop:recurrence}) under
appropriate assumptions on the stepsize sequence $(\gamma_n)_{n \ge 1}$.

Notice that to adapt the stepsize sequence to the
already visited states, it is natural to choose a stepsize $\gamma_{n+1}$ which
depends on the past (namely $\gamma_{n+1}$ is a function of $(\tilde{\theta}_0,X_1, \ldots
X_n)$).  Our convergence analysis of the WL$_\rho$ algorithm allows for such
random stepsize sequences satisfying some summability assumptions.

\paragraph{Two examples of generalized Wang-Landau algorithms.}

For the original Wang-Landau
algorithm~\cite{wang-landau-01}, the function $\rho$ is $\rho(t)=t$ so that the target measure
at convergence is $\pi_{\theta_\star}$ which gives equal weight to all
the strata. The updating rule considered in the mathematical analysis provided in~\cite{fort:jourdain:kuhn:lelievre:stoltz:2014} relies on 
the following recurrence relation
(compare to~\eqref{eq:evol_algo_gen}):
\begin{equation}\label{eq:WL}
\tu^{\wl}_{n+1}(i) = \tu^{\wl}_{n}(i) \Big( 1 + \gamma^{\wl}_{n+1} \, \un_{\Xset_i}(X_{n+1}) \Big), \qquad \tn^{\wl}_{n+1}(i) \eqdef  \frac{\tu^{\wl}_{n+1}(i)}{\sum_{j=1}^d \tu^{\wl}_{n+1}(j)} 
\end{equation}
where $(\gamma^{\wl}_n)_{n \ge 1}$ is a deterministic positive stepsize sequence chosen
by the user. 
The stepsize sequence $(\gamma^{\wl}_n)_{n \ge 1}$ gives the penalization
strength. This sequence should decrease to zero in order for the normalized
sequence $(\t^{\wl}_{n})_{n \ge 0}$ to have a limit, but not too fast since one
wants $(\t^{\wl}_{n})_{n \ge 0}$ to converge to $\theta_\star$. The convergence
of the sequence $(\tn^{\wl}_n)_{n \ge 0}$ implies the convergence of the
distribution of $X_n$ to $\pi_{\tn_\star}$ (see e.g.~\cite[Theorems 3.3 and
3.4]{fort:jourdain:kuhn:lelievre:stoltz:2014}). Actually, in the original Wang Landau
algorithm~\cite{wang-landau-01}, the stepsize is divided by
2 each time the occupation measure of the strata is close to uniform, up to
an error related to the current value of the stepsize (see~\cite{JR11} for a
mathematical analysis). The updating rule is thus not completely
deterministic, since it involves a random time. 
The original Wang-Landau algorithm~\cite{wang-landau-01} is one instance of a generalized
Wang-Landau algorithm. 
Adapting our analysis to handle such an updating rule (by checking
that the hypotheses of Proposition~\ref{prop:cvggal} below are
satisfied for this updating rule) would be an
interesting contribution.

The SHUS$_\rat^\alpha$ algorithm is
the WL$_\rat$ algorithm for the specific stepsize
sequence:
\begin{equation}
\label{eq:DefinitionGamma}
\gamma_{n+1}=\frac{\gamma}{g_\alpha(S_n)},
\end{equation} 
since~\eqref{eq:evol_algo_gen}--\eqref{eq:DefinitionGamma} is equivalent
to~\eqref{eq:def:thetatilde}.
 We will see below
that the parameter $\alpha \in (\frac12,1]$ gives the limiting behaviour of the
stepsize sequence~\eqref{eq:DefinitionGamma}: in the large $n$ limit, $\gamma_n
\simeq n^{-\alpha}$ (see Section~\ref{sec:discuss_algo} for a formal argument,
and Corollary~\ref{cor:asymptotic} below for a rigorous derivation).

The SHUS$_\rat^\alpha$ algorithm thus differs from the original Wang Landau algorithm since the target density is not $\pi_{\theta_\star}$ but $\pi^\rho_{\theta_\star}$, so
that at convergence the $i$-th stratum has probability
$(Z^\rho_{\theta_\star})^{-1}\frac{\theta_\star(i)}{\rho(\theta_\star(i))}$.
Compared to Wang-Landau, the two main parameters $\rat$ and $\alpha$ of the SHUS$_\rat^\alpha$
algorithm thus introduce flexibility in the algorithm.
The function $\rho$ allows to balance the two objectives of the
 importance sampling strategy: reducing the metastable features of a Markov
 chain targeting the original probability measure $\pi \, \rmd \lambda$ without
 reducing too much the effective sample size of the weighted samples in~\eqref{eq:IS_rho}.
 The parameter $\alpha$ provides a control on the step-size sequence
 of the stochastic approximation algorithm (see Remark~\ref{rem:stepsize} below
 for a discussion on the interest of controlling the step-size sequence).


\paragraph{Comparison with other free energy adaptive techniques.}
The SHUS$_\rat^\alpha$ algorithm is a generalization of algorithms which are
widely used in practice as efficient sampling tools in molecular dynamics. In
particular, the original SHUS algorithm~\cite{SHUS} corresponds to the choices
$\alpha=1$ and $\rat(t)=t$, see~\cite{fort:jourdain:lelievre:stoltz:2015} for a
mathematical analysis. Besides, the well-tempered
metadynamics~\cite{barducci-bussi-parrinello-08} corresponds to the choices
$\alpha=1$, $\mu=1-a$ and $\rat(t)=t^a$ with $a \in (0,1)$, as shown in
Section~\ref{sec:WT} below.

\subsection{Discussion of the construction of the algorithm}\label{sec:discuss_algo}

As explained in the introduction, adaptive free energy biasing algorithms are
designed to update the parameter $\tn_n$ in such a way that $(\tn_n)_{n \ge 0}$
converges to $\tn_\star$ (defined by~\eqref{eq:theta_star}). We motivate in
this section the choices of the updating rule~\eqref{eq:def:thetatilde} (in the
case $\alpha = \mu = 1$) and of the function $g_\alpha$ defined in~\eqref{hyp:g}.

\paragraph{Motivation of the updating rule~\eqref{eq:def:thetatilde} when $\alpha =\mu= 1$.}
Let us explain heuristically, in the case $\alpha=\mu=1$, the reason why
$\lim_n \tn_n = \theta_\star$, assuming that $(\tn_n)_{n \ge 0}$ converges to 
some $\theta_\infty \in \Theta$. Then, asymptotically, 
everything happens as if the states $X_k$ were sampled under
$\pi_{\tn_\infty}^\rat \, \rmd \lambda$ and were satisfying a strong
law of large number.
Hence, for any $i \in \{1, \ldots,d\}$, almost-surely,
\[
\lim_{n \to +\infty} \frac1n \sum_{k=1}^n \un_{\Xset_i}(X_k) = \int_{\Xset_i}
\pi^\rat_{\theta_\infty} \, \rmd\lambda=(Z^\rat_{\theta_\infty})^{-1} \ 
\frac{\tn_\star(i)}{\rat(\theta_\infty(i))}.
\]
Since $g_\alpha(s) =s$, \eqref{eq:def:thetatilde} implies
$$
\tu_{n}(i) = \tu_0(i) + \gamma \sum_{k=1}^n \rat(\t_{k-1}(i))\,
\un_{\Xset_i}(X_{k}),$$
so that, almost-surely, for any $i \in \{1,
\ldots,d\}$,
\[
\lim_{n \to +\infty} \frac{ \tu_n(i)}{n} = \gamma \, \rat(\theta_\infty(i)) \,
(Z^\rat_{\theta_\infty})^{-1}
\frac{\theta_\star(i)}{\rat(\theta_\infty(i))}=\gamma \,
\frac{\theta_\star(i)}{Z^\rat_{\theta_\infty}}.
\]
By summing over $i=1, \ldots, d$, one thus gets that
\begin{equation}\label{eq:heuristic:Sn}
  \lim_{n \to +\infty} \frac{ S_n}{n} = \frac{\gamma}{
  Z^\rat_{\theta_\infty}}.
\end{equation}
Therefore, since $\tn_n = \tu_n /S_n$, we have $\lim_n \tn_n = \tn_\star$;
hence, $\tn_\infty = \tn_\star$. This is not a proof of convergence, but an
indication that the only reasonable limit for $(\tn_n)_{n \ge 0}$, when it
exists, is $\tn_\star$.


\paragraph{Choice of the function $g_\alpha$.}
Let us now explain the role of the function~$g_\alpha$
in~\eqref{hyp:g}, by looking at the asymptotic behavior of
the stepsize sequence $(\gamma_{n})_{n \ge 1} = (\gamma/g_\alpha(S_{n-1}))_{n \ge 1}$ when $n \to \infty$. As explained above, when $\alpha=\mu=1$, $S_n$ scales as $\gamma
(Z^\rho_{\theta_\star})^{-1} n$ (see~\eqref{eq:heuristic:Sn}). Therefore the stepsize
$\gamma_{n+1}=\gamma/S_n$ scales like
$Z^\rho_{\theta_\star}/n$. As discussed
in~\cite{amrx,fort:jourdain:kuhn:lelievre:stoltz:2014,fort:jourdain:lelievre:stoltz:2015}, it may be
interesting in practice to use larger stepsizes, of order $n^{-\alpha}$ with $\alpha\in(\frac{1}{2},1)$, in order to leave more
quickly metastable states (this will be discussed in more details in
Remark~\ref{rem:stepsize} below). Let us check that this can be performed by
the choice~\eqref{hyp:g} of the function~$g_\alpha$. 

In order to understand the possible choices for the function $g_\alpha$, we
consider a generalized updating rule
\begin{equation}
   S_{n+1} = S_n +  \frac{\gamma}{g(S_n)} \, S_n  \, \rat\Big(\t_n(I(X_{n+1}))\Big),\label{eq:evolsn}
\end{equation}
where $g_\alpha$ in~\eqref{eq:def:thetatilde} has been replaced by any function
$g:\R^*_+ \to \R^*_+$, and
where, we recall, $I$ is defined by (see~\eqref{eq:def_I}):
\begin{equation*}
I(x) = i \qquad \text{if and only if}  \ \ x \in \Xset_i.
\end{equation*}
We define accordingly the generalized stepsize sequence by
$\gamma_{n+1}=\frac{\gamma}{g(S_n)}$. We now follow a formal reasoning,
comparing the asymptotic behaviors of sequences with the asymptotic behaviors
of the associated ordinary differential equations. All these computations will
be rigorously justified in Section~\ref{sec:convergence:results}.  Since
$\rho(\theta_n(I(X_{n+1})))$ is expected to average out, in the longtime limit,
at
$$\int_{\Xset} \rho(\theta_\star(I(x))) \pi^\rho_{\theta_\star}(x) \rmd
\lambda(x) = (Z_{\theta_\star}^\rho)^{-1},$$
one thus expects $S_n$ and
$\gamma_n$ to behave when $n \to \infty$ like $s(n)$ and $\gamma(n)$ where $t
\mapsto s(t)$ and $t \mapsto \gamma(t)$ satisfy
\begin{equation}\label{eq:s}
\dot{s}(t) =\frac{\gamma}{g(s(t))} s(t) (Z_{\theta_\star}^\rho)^{-1}
\text{ and } \gamma(t) = \frac{\gamma}{g(s(t))}.
\end{equation}
As explained in~\cite{fort:jourdain:kuhn:lelievre:stoltz:2014}, classical
results on stochastic approximation algorithms require that the stepsize
sequence satisfies $\sum_{n \ge 1} \gamma_n = \infty$ and $\sum_{n \ge 1}
\gamma_n ^2< \infty$ in order to ensure the almost-sure convergence. In the
continous-time setting introduced above, the question is thus: which functions
$g: \R_+ \to \R_+$ are such that
\begin{equation}\label{eq:gamma_CV}
\int_{\R_+} \gamma(t) \, \rmd t = \infty \text{ and } \int_{\R_+}
\gamma^2(t) \, \rmd t < \infty
\end{equation}
where the functions $\gamma$ and $s$ are defined by~\eqref{eq:s}? In the
literature on stochastic approximation algorithms, functions which
satisfy~\eqref{eq:gamma_CV} are classically chosen as
\begin{equation}\label{eq:gamma_ex}
\gamma(t)=\left(\frac{\nu}{  t}\right)^\alpha
\end{equation}
where $\nu >0$ and $\alpha \in (\frac12,1]$ are two parameters. By
using the fact that, from the first equation in~\eqref{eq:s}, it holds 
$\gamma(t)=Z_{\theta_\star}^\rho \frac{\rmd \ln (s(t))}{\rmd t}$, one obtains that $s(t)=C\exp\left(\frac{\nu^\alpha t^{1-\alpha}}{
    Z_{\theta_\star}^\rho (1-\alpha)}\right)$ when $\alpha \in (\frac12,1)$, and
$s(t)=C t^{\frac{\nu}{Z_{\theta_\star}^\rho}}$ when $\alpha=1$ (for some
positive constant $C$). Since $g(s(t))=\frac{\gamma}{\gamma(t)}=\frac{\gamma
  t^\alpha}{\nu^\alpha}$, one thus gets: $\forall t > 0$,
\begin{equation}\label{eq:possible_g}
\left\{
\begin{aligned}
g\left(C\exp\left(\frac{\nu^\alpha t^{1-\alpha}}{Z_{\theta_\star}^\rho (1-\alpha)}\right)\right)
= \frac{\gamma t^\alpha}{\nu^\alpha} &\text{ when } \alpha \in
\left(\frac12,1\right) \, , \\
g\left(C t^{\frac{\nu}{Z_{\theta_\star}^\rho}} \right)
= \frac{\gamma t}{\nu} &\text{ when } \alpha =1 \, .\\
\end{aligned}
\right.
\end{equation}
Taking into account the fact that $\lim_{t \to \infty} s(t) = \infty$ to get
rid of irrelevant constants, one can check that the choices~\eqref{hyp:g} of
$g$ are consistent with the equalities~\eqref{eq:possible_g}, for the following
choices of $\nu$:
$$\nu=\left\{
\begin{aligned}
 \gamma^{\frac{1}{\alpha}-1}
  Z_{\theta_\star}^\rho (1-\alpha) &\text{ for } \alpha \in
  \left(\frac12,1\right) \, ,\\
\frac{Z_{\theta_\star}^\rho}{\mu} &\text{ for } \alpha =1 \, .\\
\end{aligned}
\right.
$$
We will see that these heuristics are compatible with rigorous mathematical
analysis (see Corollary~\ref{cor:asymptotic}) and numerical analysis (see
Section~\ref{sec:application}).

\begin{arem}\label{rem:stepsize}
The objective of this remark is to discuss the interest of considering
larger stepsizes ($\gamma_n \simeq n^{-\alpha}$ with $\alpha$ close to
$1/2$) than in the original SHUS or well tempered metadynamics
(for which $\alpha=1$). It is known that in a Stochastic Approximation
  updating rule, the choice of the stepsize sequence $(\gamma_n)_{n \geq 1}$
  plays a role on the efficiency of the algorithm.

On the one hand, before reaching
  equilibrium (namely in the transient phase), it is better to choose
  $(\gamma_n)_{n \geq 1}$ slowly decreasing in order to overcome a possible poor initialization (see
  e.g.~\cite{amrx,fort:jourdain:kuhn:lelievre:stoltz:2014,fort:jourdain:lelievre:stoltz:2015}). This
  is illustrated in Section~\ref{sec:exit_num}
  where it is shown that choosing $\alpha$ close to $1/2$ leads to
  exit times from metastable states which are much smaller than for $\alpha=1$.

On the other hand,  when the system gets closer to equilibrium,
$(\gamma_n)_{n \geq 1}$  should decrease rapidly to reduce the
asymptotic fluctuations of $\tn_n$ around $\tn_\star$ (see
  e.g.~\cite{benveniste:metivier:priouret:1987}). More precisely, it
  can typically be shown that $(\gamma_n^{-1/2} (\tn_n -
  \tn_\star))_{n \ge 1}$
  converges in distribution to a centered Gaussian distribution in the
  limit $n \to \infty$, see
  for example~\cite[Theorem
  3.6]{fort:jourdain:kuhn:lelievre:stoltz:2014} for the case of the
  Wang Landau algorithm. Thus, in this regime, the smaller $\gamma_n$, the better.

In practice, the above considerations indicate that one should use large timesteps in a first stage of the algorithm, and then
smaller ones (hence the interest of being able to control the
parameter~$\alpha$). Another idea is to combine an updating rule with large
timesteps with an averaging technique to recover an asymptotic
variance of order $1/n$ (see for example~\cite{R88,PJ92}
or~\cite[Theorem 3.2]{fort-2015}). It is not our objective to further explore such issues in this article.
\end{arem}

\section{Long-time behavior of the SHUS$_\rat^\alpha$ algorithm} 
\label{sec:convergence:results}

\subsection{General assumptions}
We study the convergence of the algorithm under the following assumption on the
target probability measure $\pi \, \rmd \lambda$ and on the strata
$(\Xset_i)_{i \in \{1, \ldots ,d\}}$: 
\debutA
\item \label{hyp:targetpi} The density $\pi$ of the target distribution is such
  that $\sup_\Xset \pi < \infty$ and the strata
  $(\Xset_i)_{i \in \{1, \ldots ,d\}}$ satisfy $\min_{1 \leq i \leq d }
  \pi(\Xset_i)> 0$.  
\finA
It is assumed that the Markov transition kernels $\{P^\rat_\tn, \tn \in \Theta \}$
satisfy:
\debutA
\item \label{hyp:kernel} For any $\tn \in \Theta$, $P^\rat_{\tn}$ is a
  Metropolis-Hastings transition kernel with proposal kernel $q(x,y)\,\rmd
  \lambda(y)$ where $q(x,y)$ is symmetric and satisfies $\inf_{\Xset^2} q > 0$,
  and with invariant distribution $\pi^\rat_{\tn} \, \rmd \lambda$, where
  $\pi^\rat_{\tn}$ is given by \eqref{defpithetaf}. 
\finA 
The assumption $\inf_{\Xset^2} q > 0$ is particularly useful for the proof of
recurrence\footnote{We recall that the algorithm is recurrent if
  with probability one, the sequence $(\tn_n)_{n \ge 0}$ enters infinitely often a compact
  subset of $\Theta$, see~\eqref{condstab} below.}.
It is unclear how to
adapt our argument to a setting where the support of $q$ is not $\Xset^2$.
However, it is expected, when the recurrence property holds, that the
convergence proof can be adapted without this assumption (see
e.g.~\cite{andrieu:moulines:priouret:2005,fort:moulines:priouret:2011,fort:moulines:schreck:vihola:2015}).

  The convergence results for SHUS$_\rat^\alpha$ are established for any
  measurable function $\rat: (0,1) \to (0, +\infty)$ satisfying:  
 \debutF
\item \label{hyp:LypFctRat} For any compact subset $\calK$ of $(0,1)$, there
  exists a constant $C$ such that
\[
\forall t\in\calK,\;\forall t'\in (0,1),\;\left| 1 - \frac{\rat(t')}{\rat(t)} \right| \leq C \left| t'-t\right|.
\]
\item \label{hyp:BoundFctRat} $\sup_{t\in (0,1)}\rat(t) < \infty.$
\item \label{hyp:fdec} $\rat$ is non-decreasing on $(0,1)$ and  there exists $R >1$ such that
\[
\inf_{t\in(0,1/R)}\frac{\rat(t)}{\rat(R t)}>0.
\]
\item
  \label{hyp:ffininf} $t \mapsto \rat(t) /t $  is  non-increasing on $(0,1)$ and $\lim_{t\to 0^+}\rat(t)/t=+\infty$.
\item\label{hyp:finf} $\inf_{t\in(0,1)}\rat(t)/t>0.$ \finF

The function $\rho: t \mapsto (\max(t,t_0))^a$ for some $a \in [0,1)$
and $t_0\in [0,1]$ (as well as for $a=1$ and $t_0 \in (0,1]$) satisfies the
above assumptions. It clearly satisfies Assumptions
  R\ref{hyp:BoundFctRat} to R\ref{hyp:finf}. It also satisfies Assumption
  R\ref{hyp:LypFctRat}: this is obvious for $a=0$, and for $a\in
  (0,1]$, it is checked as follows. For $t_0=0$ and $t,t'\in (0,1)$ either $t'\le \frac{t}{2}$ and
  $|1-(\frac{t'}{t})^a|\le 1\le \frac{2}{t}|t'-t|$ or $t'\ge \frac{t}{2}$ and
  $|1-(\frac{t'}{t})^a|=\frac{|(t')^a-t^a|}{t^a}\le
  \frac{a}{t^a}\left(\frac{t}{2}\right)^{a-1}|t'-t|$. The latter
  calculation also shows that, for $t_0>0$, $\left| 1 - \frac{\rat(t')}{\rat(t)}
  \right| \leq C \left| \max(t',t_0)-\max(t,t_0)\right| \le C \left| t'-t\right|$.
  As already mentioned above, the case $t_0=0$ and $a \in (0,1]$ is typically the
  case of interest in practice. Notice that
 recent papers~\cite{DHSV15,MVTP15} consider the case of a function $\rho$ which,
 like $t \mapsto (\max(t,t_0))^a$ with $t_0>0$, is constant in a neighborhood of~$0$.
When $a=0$, the algorithm actually
  corresponds to the naive sampling, without adaptation.

  Note that these assumptions are not satisfied for the function $\rat(t)=t$.
  The SHUS$^\alpha_{\rho(t)=t}$ case has actually already been studied in~\cite{fort:jourdain:lelievre:stoltz:2015},
  where similar asymptotic results as those presented below are proven
  under the additional hypothesis $\inf_\Xset\pi>0$ used to check the
  recurrence property. The proofs in our case follow the same lines as in~\cite{fort:jourdain:lelievre:stoltz:2015}. Compared to~\cite{fort:jourdain:lelievre:stoltz:2015}, the fact that $\lim_{t \to 0} \rat(t)/t=+\infty$ implies that strata which have not been visited  are more penalized (see~\eqref{eq:evol_algo_gen}). This makes the recurrence of the algorithm easier to establish, so that we could get rid of the assumption $\inf_\Xset \pi  > 0$ needed in~\cite{fort:jourdain:lelievre:stoltz:2015}.

  \subsection{Convergence results}

  Our main result is the following convergence result.
  \begin{thm}\label{theo:cvggal}
    Let $\gamma >0$. Assume A\ref{hyp:targetpi},
    A\ref{hyp:kernel} and R\ref{hyp:LypFctRat} to R\ref{hyp:finf}.
    Assume moreover that either (i) $\alpha \in (\frac12,1)$ or (ii) $\alpha=1$
  and $\mu \ge 1$  or (iii) $\alpha=1$, $\mu \in (0,1)$ and
  $\rho(t)=t^a$ for some $a \in [0,1)$. 
 Then the
    SHUS$_\rat^\alpha$ algorithm starting from any
    $(\R_+^*)^d\times\Xset$-valued random initial condition $(\tu_0,X_0)$
    converges in the following sense:
\begin{enumerate}
  \item[(i)] $ \lim_{n \to +\infty} \tn_n
      = \tstar$ $\P$-a.s.
\item[(ii)] For any bounded measurable function $f: \Xset \to \R$,
\begin{align*}
  \lim_{n\to\infty} \E\left[f(X_n)\right] = \int_\Xset f \,
  \pi^\rat_{\tn_\star} \, \rmd \lambda \qquad \mbox{ and } \qquad \lim_{n\to\infty}
  \frac{1}{n} \sum_{k=1}^n f(X_k) = \int_\Xset f \, \pi^\rat_{\tn_\star} \,
  \rmd\lambda \quad \mathrm{\P-a.s.}\eqsp \nonumber
  \end{align*}
\item[(iii)] For any bounded measurable function $f: \Xset \to \R$,
\begin{align*}
  & \lim_{n\to\infty} \E\left[\left(\sum_{i=1}^d
    \frac{\theta_{n-1}(i)}{\rat(\theta_{n-1}(i))}\right)\left(\sum_{j=1}^d
    \rat(\theta_{n-1}(j))\,f(X_n) \un_{\Xset_j}(X_n)\right)\right] = \int_\Xset f \,
  \pi \, \rmd \lambda \eqsp, \nonumber
  \\
  & \lim_{n\to\infty}\frac{1}{n} \sum_{k=1}^n \left( \sum_{i=1}^d
  \frac{\theta_{k-1}(i)}{\rat(\theta_{k-1}(i))} \right)\left( \sum_{j=1}^d
  \rat(\theta_{k-1}(j)) \, f(X_k) \un_{\Xset_j}(X_k) \right) = \int_\Xset f \, \pi \,
  \rmd\lambda\quad \mathrm{\P-a.s.}
  \label{theo:ergoLLN}
  \end{align*}
\end{enumerate}
\end{thm}
Let us recall that $\tstar$ and $\pi_{\tn}^\rat$ are respectively defined in \eqref{eq:theta_star} and
\eqref{defpithetaf}.
 Here and in the following, $\P$ and $\E$ respectively denote the
probability measure and the expectation on any probability space which
supports all the random variables needed to define the algorithm. 

Notice that we have not been able to prove convergence for $\alpha=1$,
$\mu \in (0,1)$ and a general function $\rat$ satisfying
R\ref{hyp:LypFctRat} to R\ref{hyp:finf}. In this case, convergence
however holds for the specific choice $\rat(t)=t^a$ which is the case
of interest in practice. We refer to Section~\ref{sec:gamma_SHUS} for more comments
on that.

The key property for the proof of convergence of the sequence $(\tn_n)_{n \ge 0}$ to $\tn_\star$ is to
rewrite the updating rule of the weight sequence $(\tn_n)_{n \ge 0}$ as in
Lemma~\ref{lem:SHUSalpharho:AS} below. This allows us to consider it as a stochastic approximation
algorithm (see
e.g.~\cite{benveniste:metivier:priouret:1987,chen:2002,kushner:yin:2003,borkar:2008,kushner:2010})
with \textit{(i)} a random stepsize sequence $(\gamma_n)_{n \ge 1}$ and \textit{(ii)}
Markovian inputs $X_{n+1}$ with a distribution conditional to the past,
depending on the current estimate $\tn_n$ of $\t_\star$. Let us recall
the definition~\eqref{eq:DefinitionGamma} of the stepsize
\begin{equation*}
\gamma_{n+1} \eqdef  \frac{\gamma}{g_\alpha(S_n)}.
\end{equation*}
Let us introduce the function $H: \Xset \times \Theta \to \R^d$ with
$i$-th component 
\begin{align}
H_i(x,\theta)  &\eqdef  \rat(\theta(i)) \, \un_{\Xset_i}(x)-  \t(i) \,
\rat (\theta(I(x)))  \label{eq:DefinitionH} \\
&= \rat(\theta(I(x))) \, (\un_{\Xset_i}(x)-  \t(i)), \nonumber
\end{align}
where $I(x)$ is defined in~\eqref{eq:def_I}. We do not indicate explicitly the dependence
of $H$ on $\rho$ for the ease of notation.

The following lemma is useful to rewrite the evolution of the
sequence $(\tn_n)_{n \ge 0}$. We state it in the  general
setting of the WL$_\rho$ algorithm (see Section~\ref{sec:relation_FE}).
\begin{alem}
  \label{lem:SHUSalpharho:AS} 
Let us consider a sequence $(\tu_n)_{n \ge 0}$ satisfying the
recurrence relation~\eqref{eq:evol_algo_gen} for some stepsize
sequence $(\gamma_n)_{n \ge 1}$. Then the associated normalized sequence
$(\tn_n)_{n \ge 0}$ defined by~\eqref{eq:def:theta_WL}
  satisfies the equation
  \begin{equation}\label{eq:AlgoStoch}
  \tn_{n+1} = \t_n + \gamma_{n+1} H(X_{n+1}, \tn_n) + \gamma_{n+1}
  \Lambda_{n+1},
\end{equation}
where $(\Lambda_n)_{n\ge 1}$ is a random sequence with values in
$\R^d$ such that 
\begin{equation}\label{eq:hyp_Lambda}
\sup_{n \ge 0} \frac{ |\Lambda_{n+1}| }{\gamma_{n+1}}
\leq \sqrt{2}\, \left(\sup_{(0,1)} \rat\right)^2 \quad \mathrm{\P-a.s.}
\end{equation}
\end{alem}
The proof is a consequence of simple computations made in Section~\ref{sec:proof:lem:SHUSalpharho:AS}.

From the property stated in Lemma~\ref{lem:SHUSalpharho:AS}, the general
strategy to prove Theorem~\ref{theo:cvggal} is to combine convergence results
for stochastic approximation algorithms (see e.g.
\cite{benveniste:metivier:priouret:1987,kushner:yin:2003,andrieu:moulines:priouret:2005})
and for so-called adaptive MCMC algorithms (see in
particular~\cite{fort:moulines:priouret:2011}). Indeed, the WL$_\rat$
algorithm can essentially be seen as an iterative procedure on $(\tn_n)_{n \ge
  0}$ (see Lemma~\ref{lem:SHUSalpharho:AS}) with $X_{n+1}$ generated according
to the distribution $P^\rat_{\tn_n}(X_n, \cdot)$.

Let us introduce the mean field function $h: \Theta \to \R^d$
(omitting again to explicitly indicate the dependence on $\rat$)
\begin{equation}\label{eq:Definitionh}
h(\theta) \eqdef  \int_{\Xset}
H(x,\tn)\, \pi^\rat_{\theta}(x) \, \rmd \lambda(x)=\frac{\tstar-\tn}{\sum_{i=1}^d\frac{\theta_\star(i)}{\rat(\theta(i))}}.
\end{equation}
In the physics literature, $h$ is sometimes called the quasi-equilibrium
average of $H$.
The recursion~\eqref{eq:AlgoStoch} is a noisy version of the dynamics driven by
the mean field function (we use here the notation
of~\cite{andrieu:moulines:priouret:2005}):
$$\tn_{n+1}(i)= \tn_n(i) + \gamma_{n+1}
h_i(\tn_n)+\gamma_{n+1} \xi_{n+1}(i),
$$
where
$$\xi_{n+1}(i)\eqdef H_i(X_{n+1},\tn_n)-h_i(\tn_n)+\Lambda_{n+1}(i).
$$
For the proof of convergence of $(\tn_n)_{n \ge 0}$, the main steps are to show that
\textit{(i)} the noise $\gamma_{n+1} \xi_{n+1}$ is sufficiently small so that
the sequence inherits the behavior of the sequence $(\tau_n)_{n \ge
  0}$ satisfying the recurrence relation $\tau_{n+1} = \tau_n + \gamma_{n+1} h(\tau_n)$, and
\textit{(ii)} that the sequence $(\tn_n)_{n \ge 0}$ converges to the zero of $h$,
namely $\theta_\star$ (see~\eqref{eq:Definitionh}), by identifying the
recurrence relation as a time discretization of the ordinary
differential equation $\dot{\tau}=h(\tau)$.  More precisely,
Theorem~\ref{theo:cvggal} is a consequence of Proposition~\ref{prop:cvggal} and
Proposition~\ref{prop:pasetstab} stated below.
Proposition~\ref{prop:cvggal} shows that convergence holds as soon as the algorithm is recurrent
and the stepsize sequence $(\gamma_n)_{n\ge 1}$ a.s.  satisfies the usual
conditions $\sum_{n\geq 1}\gamma_n^2<\infty$ and $\sum_{n\ge 1}\gamma_n=\infty$
for the convergence of stochastic approximation algorithms, respectively to control
the noise and to mimic the asymptotic behavior of the ordinary differential
equation. Proposition~\ref{prop:pasetstab}
ensures that those two conditions are actually satisfied in our
setting. The proof of  Proposition~\ref{prop:pasetstab} is based in
particular on some sufficient conditions on the
sequence $(\gamma_n)_{n \ge 1}$ for the recurrence of the algorithm,
stated in Proposition~\ref{prop:recurrence}.

We state Propositions~\ref{prop:cvggal} and~\ref{prop:recurrence} in the  general
setting of the WL$_\rho$ algorithm which encompasses the  SHUS$_\rat^\alpha$
algorithm  (see Section~\ref{sec:relation_FE}).
\begin{aprop}\label{prop:cvggal} 
Assume we are given a density $\pi$ and a family of kernels
$P^\rho_\theta$ satisfying  A\ref{hyp:targetpi} and A\ref{hyp:kernel},
for a function $\rho:(0,1) \to (0,+\infty)$ satisfying R\ref{hyp:LypFctRat} and
  R\ref{hyp:BoundFctRat}. Assume that we are given sequences $(\tu_n,X_n)_{n \ge 0}$ and
$(\gamma_n)_{n\ge 1}$ generated by the WL$_\rho$ algorithm (see Algorithm~\ref{algo:WL}). 
In particular, for all $n \ge 0$, conditionally on
\begin{equation}\label{eq:Fn}
{\cal F}_n \eqdef\sigma\left(\tu_0, X_0, X_1, \cdots, X_n \right),
\end{equation}
$X_{n+1}$ is generated
according to the distribution $P^\rat_{\tn_n}(X_n, \cdot)$.

  If, moreover, \begin{itemize}
\item the stepsize sequence $(\gamma_n)_{n\ge 1}$ is predictable with
  respect to the filtration $({\cal F}_{n})_{n \ge 0}$ (i.e. for all $n \ge 1$,
  $\gamma_{n}$ is ${\cal F}_{n-1}$-measurable) and such that \begin{equation}
   \label{hyp:Gamma}\text{$\P$-a.s., } (\gamma_n)_{n\ge 1}\mbox{ is
     non-increasing},\;\sum_{n\ge 1} \gamma_n =\infty \text{ and }\sum_{n \geq
     1} \gamma_n^{2} < \infty ,
\end{equation}
\item the algorithm is recurrent, in the following
  sense: 
\begin{equation}\label{condstab}
\text{$\P$-a.s., the sequence $(\tn_n)_{n \ge 0}$ returns
  infinitely often to a compact set of $\Theta$,}
\end{equation}\end{itemize} then the conclusions (i)-(ii)-(iii) of Theorem \ref{theo:cvggal} hold. 
\end{aprop}

Some sufficient conditions on $(\gamma_n)_{n \ge 1}$ to ensure the recurrence of the WL$_\rho$
algorithm are given in the following proposition.
\begin{aprop}\label{prop:recurrence} 
Assume we are given a density $\pi$ and a family of kernels
$P^\rho_\theta$ satisfying  A\ref{hyp:targetpi} and A\ref{hyp:kernel},
for a function $\rho:(0,1) \to (0,+\infty)$ satisfying R\ref{hyp:BoundFctRat} to
  R\ref{hyp:ffininf}. Assume that we are given sequences $(\tu_n,X_n)_{n \ge 0}$ and
$(\gamma_n)_{n\ge 1}$ satisfying the recurrence relation~\eqref{eq:evol_algo_gen} of the
WL$_\rho$ algorithm. 

If the sequence $(\gamma_n)_{n\geq 1}$ is
  non-increasing, bounded from above by a deterministic sequence converging to
  $0$ as $n\to\infty$ and such that
  \begin{equation}\label{eq:hypr}
\bar{r}_{d,\gamma} \eqdef \sup_{n\geq
    1}\frac{\gamma_{n}}{\gamma_{n+d-1}}<+\infty,
\end{equation} 
then the algorithm is recurrent, in the sense of~\eqref{condstab}.
\end{aprop}

Notice that Propositions~\ref{prop:cvggal} and~\ref{prop:recurrence} give
some sufficient conditions on $(\gamma_n)_{n\geq 1}$ for the
convergence of the WL$_\rho$ algorithm.

\begin{aprop}\label{prop:pasetstab}  Let $\gamma >0$. 
  Assume A\ref{hyp:targetpi}, A\ref{hyp:kernel}, R\ref{hyp:BoundFctRat} to
  R\ref{hyp:finf}.  Then the sequences 
  $\left(\tu_{n},X_n\right)_{n \ge 0}$ and
  $\left(\gamma_{n+1}=\frac{\gamma}{g_\alpha(S_n)}\right)_{n \ge 0} $ generated by the SHUS$_\rat^\alpha$ algorithm
  starting from any $(\R_+^*)^d\times\Xset$-valued random initial condition
  $(\tu_0,X_0)$ are such that~\eqref{hyp:Gamma} and~\eqref{condstab}
  hold in the following cases: (i) $\alpha \in (\frac12,1)$ or (ii) $\alpha=1$
  and $\mu \ge 1$  or (iii) $\alpha=1$, $\mu \in (0,1)$ and
  $\rho(t)=t^a$ for some $a \in [0,1)$. 
\end{aprop}

A useful corollary of the previous results gives
the effective behavior of the stepsize sequence $(\gamma_n)_{n \ge 1}$ in the longtime limit
$n \to +\infty$. 

\begin{acor}\label{cor:asymptotic}  Let $\gamma >0$. 
  Assume A\ref{hyp:targetpi}, A\ref{hyp:kernel} and R\ref{hyp:LypFctRat} to
  R\ref{hyp:finf} and that either  (i) $\alpha \in (\frac12,1)$ or (ii) $\alpha=1$
  and $\mu \ge 1$  or (iii) $\alpha=1$, $\mu \in (0,1)$ and
  $\rho(t)=t^a$ for some $a \in [0,1)$. 

Then, the stepsize sequence $(\gamma_n)_{n \ge 1}$ generated
  by the SHUS$_\rat^\alpha$ algorithm starting from any
  $(\R_+^*)^d\times\Xset$-valued random initial condition $(\tu_0,X_0)$
  has the following asymptotic behavior:
  \[
  \P-a.s., \quad \lim_{n\to\infty} \left( n^\alpha \gamma_n \right)^{1/\alpha} = \mathcal{C}_\alpha \,
 Z_{\theta_\star}^\rho  \quad \text{with}
  \quad \mathcal{C}_\alpha \eqdef \left\{
  \begin{array}{ll}
    1/\mu & \text{if $\alpha =1$,} \\
    (1-\alpha) \, \gamma^{(1/\alpha)-1} & \text{if $\alpha \in (1/2,1)$.}
  \end{array}
\right.
\]
where, we recall $Z_{\theta_\star}^\rho=\sum_{i=1}^d
\frac{\theta_\star(i)}{\rat(\theta_\star(i))}$ (see~\eqref{defpithetaf}).
\end{acor}

The proofs of all these results are gathered in Section~\ref{sec:proofs}.

\section{Relationship with well-tempered metadynamics}
\label{sec:WT}

The objective of this section is to make explicit a
connection between the
SHUS$_\rat^\alpha$ algorithm and the well-tempered
metadynamics~\cite{barducci-bussi-parrinello-08}. This leads us to
propose an accelerated version of the well-tempered
metadynamics.

\subsection{Presentation of well-tempered metadynamics.}
The well-tempered
metadynamics~\cite{barducci-bussi-parrinello-08} is an adaptive biasing
procedure used in molecular dynamics, where a biasing potential $V_{\rm bias}:\R_+
\times \R \to \R$ is updated in time according to
(see~\cite[Equation~(3)]{barducci-bussi-parrinello-08})
\begin{equation}\label{eq:WT}
\forall t \ge 0, \, \forall z \in \R, \, \frac{\rmd V_{\rm bias}(t,z)}{\rmd t} = \omega \exp\left(-\frac{V_{\rm bias}(t,z)}{\Delta T}\right) \delta_\varepsilon (\xi(X_t) - z),
\end{equation}
where $\xi: \R^D \to \R$ is the so-called collective variable, and $\Delta T$
and $\omega$ are positive parameters. In this setting, the function
$\xi$ and the dummy variable $z \in \R$ respectively play the
roles of the function $I$ and the
dummy variable $i \in \{1,\ldots,d\}$. Moreover, $X_t \in \R^D$ denotes
the configuration of the system at time $t$, and $\delta_\varepsilon$ is an approximation
to the identity, typically $\delta_\varepsilon(z) = (2\pi \varepsilon)^{-1/2}
\exp(-|z|^2/(2\varepsilon))$ for a positive parameter $\varepsilon$. The biasing
potential $V_{\rm bias}$ is thus increased around the current value $\xi(X_t)$ of the
collective variable at time~$t$, in order to favor visits to other values of
$\xi$ than the current one.  We refer to~\cite{Dickson15} for
a discussion on using the unbiasing weight $\exp\left(-\frac{V_{\rm
      bias}(t,\xi(X_t))}{\Delta T}\right)$ instead of  $\exp\left(-\frac{V_{\rm
      bias}(t,z)}{\Delta T}\right)$ (this has no impact in the limit
$\varepsilon \to 0$ we consider afterwards).

The stochastic process $(X_t)_{t \ge 0}$ follows a dynamics which, at
time $t$ and for a fixed biasing potential $V_{\rm bias}(t,\cdot)$, is ergodic with
respect to the biased probability measure
\[
\pi_t(x) = Z_t^{-1} \exp \left(-\frac{V(x)+V_{\rm bias}(t,\xi(x))}{T} \right) \text{ where }
Z_t=\int_{\R^D} \!\!\!\exp\left(-\frac{V(x)+V_{\rm bias}(t,\xi(x))}{T}\right) \, \rmd x.
\]
Here, $T >0$ is the temperature, the Boltzmann constant is taken to
$1$ for simplicity and $V:
\R^D \to \R$ is the potential energy function. The original target
density obtained when $V_{\rm bias}=0$ is thus the Boltzmann-Gibbs
density:
\begin{equation}\label{eq:Boltz}
\pi(x)=Z^{-1} \exp\left(-\frac{V(x)}{T}\right) \text{ where }
Z=\int_{\R^D} \exp\left(-\frac{V(x)}{T}\right) \, \rmd x.
\end{equation}
One example of a dynamics followed by $(X_t)_{t \ge 0}$ is the
overdamped Langevin dynamics
\begin{equation}\label{eq:OL}
\rmd X_t = - [\nabla V(X_t) +\nabla_x V_{\rm bias}(t,\xi (X_t)) ] \, \rmd t +\sqrt{2 T} \, \rmd W_t,
\end{equation}
where $(W_t)_{t \ge 0}$ is a $D$-dimensional Brownian motion.
The well-tempered metadynamics algorithm thus consists in evolving the coupled
system~\eqref{eq:WT} and \eqref{eq:OL} (using some appropriate
time-discretization schemes).

As explained
in~\cite{barducci-bussi-parrinello-08,dama-parrinello-voth-14}, it is
expected that the biasing
potential $V_{\rm bias}$ admits a longtime limit. For small $\varepsilon$,  this longtime limit should be $-\frac{\Delta
  T}{T+\Delta T} F$ up to an irrelevant additive constant, where $F$ is the so-called free energy, defined by
\begin{equation}\label{eq:def_FE}
\forall z\in \xi(\R^D), \, \exp\left(-\frac{F(z)}{T}\right) = \int_{\{x, \,
  \xi(x)=z\}} \exp\left(-\frac{V(x)}{T}\right) \delta_{\xi(x)-z}
(\rmd x),
\end{equation}
see {\em e.g.}~\cite{lelievre-rousset-stoltz-book-10} for the precise meaning of the surface measure  $\delta_{\xi(x)-z}
(\rmd x)$.
In the stationary state, the sampled density is thus $Z^{-1}
\exp\left(-\frac{V(x) - \frac{\Delta T}{T+\Delta T}
    F(\xi(x))}{T}\right) \, \rmd x$, whose marginal in $\xi$ is
$Z^{-1}\exp\left(-\frac{F(z)}{T+\Delta T}\right)$. Let us explain the
heuristic argument  which gives the  longtime limit of $V_{\rm bias}$. Let us
assume that $V_{\rm bias}(t,z)$ converges to a limiting potential
$V_{\rm bias}(\infty,z)$ up to an additive constant: in the limit $t \to \infty$,
$$V_{\rm bias}(t,z) \simeq V_{\rm bias}(\infty,z) + C(t),$$
where $C: \R_+ \to \R$. By using~\eqref{eq:WT}, one gets in the limit
$\varepsilon \to 0$, using the
fact that $X_t$ is distributed according to $Z_\infty^{-1} \exp
\left(-\frac{V(x)+V_{\rm bias}(\infty,\xi(x))}{T} \right)$,
$$\begin{aligned}
\forall z \in \xi(\R^D), \frac{\rmd C}{\rmd t}&= \omega
\exp\left(-\frac{V_{\rm bias}(\infty,z)}{\Delta T}\right) \int_{\{x,
  \xi(x)=z\}} \!\!\!\!\!\!\!\!\!
Z_\infty^{-1} \exp \left(-\frac{V(x)+V_{\rm bias}(\infty,\xi(x))}{T} \right)\delta_{\xi(x)-z}
(\rmd x)\\
&= \omega
\exp\left(-\frac{T+\Delta T}{T\Delta T}V_{\rm bias}(\infty,z)\right)
Z_\infty^{-1}\exp\left(-\frac{F(z)}{T}\right).
\end{aligned}
$$
Since the left-hand side does not depend on $z$, this yields, up to an
irrelevant additive constant
$$\forall z \in \R, \, V_{\rm bias}(\infty,z)=-\frac{\Delta T}{T+\Delta T} F(z).$$
Notice that this reasoning is very similar to the one we used at the
beginning of Section~\ref{sec:discuss_algo} to
identify the limit of $(\tu_n)_{n \ge 0}$.

\subsection{Reformulating well-tempered metadynamics in a discrete
  setting.}  Let us now make explicit the connection between the well-tempered
metadynamics and the SHUS$^\alpha_\rho$ algorithm. We already explained the
link between the target $\pi$ and the potential energy function~$V$,
see~\eqref{eq:Boltz}. In the setting we consider for the SHUS$_\rat^\alpha$
algorithm, we use a discrete collective variable $\xi(x)=I(x)$
(see~\eqref{eq:def_I} for the definition of $I$). Therefore, the term
$\delta_\varepsilon(\xi(X_t)-z)$ in~\eqref{eq:WT} is simply replaced by an
indicator function $\un_{\xi(X_t)=i}$. Moreover, from~\eqref{eq:def_FE} the
free energy is (up to an additive constant) $F(i) = -T \ln \theta_\star(i)$ for
all $i\in \{1, \ldots ,d\}$. Finally, we also need to consider evolutions which are discrete in time.
We introduce to this end a timestep size $h>0$ and consider all quantities in the well-tempered
metadynamics at times $nh$.

In order to guess the relationship between
the biasing potential $V_{\rm bias}$ and the biasing vector $\tilde{\theta}$,
let us consider these two quantities in the longtime regime. On the one
hand, as explained above, the longtime limit of the biasing potential $V_{\rm
  bias}$ is $-\frac{\Delta T}{T+\Delta T} F$ up to an irrelevant additive
constant.  On the other hand, in the SHUS$_\rat^\alpha$ algorithm, the sequence
$(\tilde{\theta}_n)_{n \ge 0}$ converges, up to a multiplicative constant, to
$\theta_\star=\exp(-F/T)$. Therefore, the natural definition of
$\tilde{\theta}$ in terms of $V_{\rm bias}$ is: for all $n \ge 0$, for all
$i\in \{1, \ldots, d\}$,
\begin{equation}\label{eq:tu_WT}
\tu ^\wt_n(i) \eqdef \exp\left(\frac{V_{\rm bias}(n h,i)}{T} \frac{T+\Delta T}{\Delta T}\right)
\end{equation}
Indeed, with this relation, if
$V_{\rm bias}$ converges to $-\frac{\Delta
  T}{T+\Delta T} F$ (up to an additive constant), $\tu ^\wt_n$ converges to
$\theta_\star=\exp(-F/T)$ (up to a multiplicative constant), as expected.

Let us now rewrite the well-tempered metadynamics in terms of $\tu
^\wt_n$, using the definition~\eqref{eq:tu_WT} of~$\tu
^\wt_n$. The
dynamics~\eqref{eq:WT} rewritten in terms of $\tu ^\wt_n$ is, after time
discretization: for all $n \in \N$, for all $i\in \{1, \ldots, d\}$,
\begin{align*}
\frac{T \Delta T}{T+\Delta T}\left( \ln (\tu_{n+1}^\wt(i)) -
  \ln(\tu_{n}^\wt(i)) \right)& = \omega   h \ 
[\tu_{n}^\wt(i)]^{-\frac{T}{T+\Delta T}} \,
\un_{\xi(X_{(n+1)h})=i}.
\end{align*}
Using the approximation $\ln (\tu_{n+1}^\wt(i)) -
  \ln(\tu_{n}^\wt(i))\simeq\frac{\tu_{n+1}^\wt(i)) -
  \tu_{n}^\wt(i)}{\tu_{n}^\wt(i)}$ valid for small $h$, this yields
\begin{align}
\tu_{n+1}^\wt(i) -\tu_{n}^\wt(i)& = \omega   h \frac{T+\Delta T}{T \Delta T}  \ 
[\tu_{n}^\wt(i)]^{\frac{\Delta T}{T+\Delta T}} \,
\un_{\xi(X_{(n+1)h})=i} \nonumber \\
&= \gamma \, \frac{S^\wt_n}{(S_n^\wt)^{1-a}}\,
\left(\t_{n}^\wt(i)\right)^a\, \un_{\Xset_i}(X_{n+1}) \label{eq:newWT}
\end{align}
with 
$$S^\wt_n=\sum_{i=1}^d  \tu_{n+1}^\wt(i) \text{ , }
\t_{n}^\wt(i)=\frac{\tu_{n+1}^\wt(i) }{S^\wt_n} $$
and
\begin{equation}\label{eq:gamma_a}
\gamma = \omega h \, \frac{T+\Delta T}{T \Delta T}\text{ , }
a = \frac{\Delta T}{T+\Delta T}.
\end{equation}
Comparing~\eqref{eq:newWT} with~\eqref{eq:def:thetatilde}, one can see
that this is the
SHUS$^\alpha_\rat$ algorithm with $\alpha=1$, $\mu=1-a$ and $\rat(t)=t^a$.

Our analysis therefore provides a proof of convergence of the well-tempered
metadynamics in the specific context where the collective variable takes
values in a finite space, and the evolution of the position vector
$X_t$ is made using a Metropolis-Hastings procedure with target $Z_t^{-1} \exp
[- (V(x)+V_{\rm bias}(t,\xi(x)))/T ]$.

\subsection{Accelerating well-tempered metadynamics.}
As explained above (see in particular Remark~\ref{rem:stepsize}), it is interesting to consider the
SHUS$^\alpha_\rat$ algorithm with $\alpha \in (\frac12,1)$ in order to obtain
larger stepsizes than for the standard  well-tempered metadynamics
(for which $\alpha=1$). A natural
question is therefore: how to modify the dynamics~\eqref{eq:WT} in order to
obtain the SHUS$^\alpha_{t^a}$ algorithm instead of~\eqref{eq:newWT}?  One can
check that the natural modification of~\eqref{eq:WT} is:
\begin{equation}\label{eq:WT_acc}
\forall t \ge 0, \, \forall z \in \xi(\R^D), \qquad \frac{\rmd V_{\rm bias}(t,z)}{\rmd
  t} =
\omega f_\alpha \left(Z_t^{T,\Delta T}\right)\,\exp\left(-\frac{V_{\rm bias}(t,z)}{\Delta T}\right) \delta_\varepsilon (\xi(X_t)- z),
\end{equation}
where, in view of~\eqref{eq:tu_WT},
\[
Z_t^{T,\Delta T} = \int_{\R} \exp\left(\frac{V_{\rm bias}(t,z)}{T}
  \frac{T+\Delta T}{\Delta T}\right) \, \rmd z
\]
is the equivalent\footnote{In practice, for this integral to be finite, one should either consider a function $\xi$ with values in a
  compact space, or restrict the biasing potential to a bounded
  subset of the collective variable values. We do not enter here into
  these practical details, see for
  example~\cite{crespo2010metadynamics,mcgovern2013boundary} for
  discussions of the proper implementation of boundary conditions in
  well-tempered metadynamics.} of $S_n = \sum_{i=1}^d \tu_n(i)$, and
\[
f_\alpha(s) = \frac{s^{\frac{T}{T+\Delta T}}}{g_\alpha(s)}
\]
for $\alpha\in (\frac{1}{2},1)$. Following the previous reasoning, it
is observed that the reformulation of~\eqref{eq:WT_acc} in a discrete
setting is the SHUS$^\alpha_\rat$ algorithm with $\rat(t)=t^a$,
$\gamma$ and $a$ given by~\eqref{eq:gamma_a}, and $\alpha \in (\frac12,1)$.

In view of the numerical experiments
presented below in Section~\ref{sec:application}, we expect that this variant of the well-tempered
metadynamics should exhibit much smaller exit times from metastable
states, and thus a quicker exploration of the state space. As
explained in Remark~\ref{rem:stepsize}, it may however be useful to
switch back to $\alpha=1$ after the transient phase, or to combine
this with an averaging technique, in order to reduce the
asymptotic fluctuations of $\tn_n$ around $\tn_\star$.

\section{Numerical illustration}
\label{sec:application}

The results stated in Section~\ref{sec:convergence:results} precisely describe the asymptotic behavior of the SHUS$_\rat^\alpha$ algorithm but do not give much information
about the efficiency of the adaptive algorithm compared to the
original non adaptive one. The aim of this section is to explore on
a specific numerical example already considered in previous works
(see for example~\cite{PSLS03,metzner-schuette-vanden-eijnden-06,amrx,fort:jourdain:lelievre:stoltz:2015}) the interest of using the SHUS$_\rat^\alpha$ algorithm
in terms of computational efficiency.

\begin{figure}
  \begin{center}
    \includegraphics[width=0.5\textwidth]{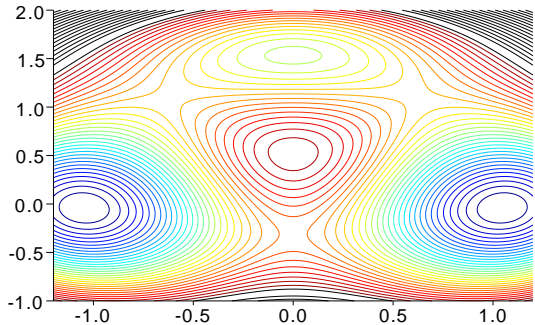}
  \end{center}
  \caption{\label{fig:contour}  Level sets of the potential~$V$ considered for the simulations. 
    The minima are located at the positions $x_\pm
  \simeq (\pm 1.05,-0.04)$, which defines two metastable states, one on the left in the vicinity of $x_-$, and one on the right
  in the vicinity of $x_+$. }
\end{figure}

We consider the system based on the two-dimensional potential suggested
in~\cite{PSLS03}. The state space is $\Xset = [-R,R] \times
\mathbb{R}$. The density of the target measure reads
\[
\forall x=(x_1,x_2)\in \Xset,\, \pi(x) =Z^{-1} \un_{[-R,R]}(x_1) \, \mathrm{e}^{-\beta V(x_1,x_2)},
\]
for some positive inverse temperature $\beta$ and $Z=\int_{\Xset}
\mathrm{e}^{-\beta V(x_1,x_2)} \, \rmd x_1 \rmd x_2$, with
\begin{align}
V(x_1,x_2)
& = 3 \exp\left(-x_1^2 - \left(x_2-\frac13\right)^2\right)
- 3 \exp\left(-x_1^2 - \left(x_2-\frac53\right)^2\right) \label{eq:pot_U} \\
& \quad - 5 \exp\left(-(x_1-1)^2 - x_2^2\right)
- 5 \exp\left(-(x_1+1)^2 - x_2^2\right)
+ 0.2 x_1^4 + 0.2 \left(x_2-\frac13\right)^4. \nonumber
\end{align}
A plot of the level sets of the potential~$V$ is presented in
Figure~\ref{fig:contour}. The global minima of the potential are located at the points
$x_-\simeq(-1.05,-0.04)$ and $x_+\simeq(+1.05,-0.04)$. This induces
two metastable states, located in the vicinities of each of the global minima.

We introduce $d$ strata $\left(\Xset_\ell = (a_\ell,a_{\ell+1}) \times
\mathbb{R}\right)_{\ell=1, \ldots,
d}$,  where $a_\ell = -R + 2(\ell-1) R/d$ for $\ell=1, \ldots,
d+1$.  In the following, we set $R=1.2$ and $d=24$. The initial weight vector
is $\tilde{\theta}_0=(1/d,\dots,1/d)$ and the
initial configuration is $X_0 = (-1,0)$. For all $\theta$, the kernel
$P^\rat_\theta$ is defined by a
Metropolis-Hastings step with target density $\pi^\rat_\theta$
and a two-dimensional Gaussian proposal distribution $q(x,y)=
\frac{1}{2\pi \sigma^2} \exp(-|x-y|^2/\sigma^2)$. We choose in the
following $\sigma^2 = 0.01$ (so that $\sigma = 2R/d$) and we use the Mersenne-Twister
random number generator as implemented in the GSL library.

 The Metropolis-Hastings dynamics without
adaptation (namely with the Gaussian proposal distribution and target $\pi$ at each iteration) is metastable:
it takes a long time (which becomes exponentially large in the limit
$\beta \to \infty$) to go from the stratum containing $x_-$ to the
stratum containing $x_+$.

The SHUS$_\rat^\alpha$ algorithm
is applied with 
$$\rat(t)=t^a$$
where $a \in [0,1]$.  
Moreover, we
choose in all numerical simulations the parameter $\gamma$ (which appears in the recurrence
relation~\eqref{eq:def:thetatilde}) as a function of $\alpha$ as follows:
\begin{equation}\label{eq:choice_alpha}
\gamma(\alpha)=
\left\{
\begin{aligned}
 1& \text{ if } \alpha = 1, \\
(1-\alpha)^{-\frac{\alpha}{1-\alpha}} & \text{ if } \alpha \in
  \left(1/2,1\right).
\end{aligned}
\right.
\end{equation}
This is to avoid the degeneracy of the constant ${\mathcal C}_\alpha(\gamma)$ which
appears in the asymptotic behavior of the stepsize sequence (see
Corollary~\ref{cor:asymptotic}) when $\alpha \to 1$. Indeed, when
$\gamma$ does not depend on $\alpha$, $\lim_{\alpha \to 1}  {\mathcal
  C}_\alpha(\gamma)=0$. With the choice~\eqref{eq:choice_alpha}, ${\mathcal
  C}_\alpha(\gamma(\alpha))$ no longer depends on $\alpha$: 
\[
\forall \alpha \in (1/2,1), \, {\mathcal C}_\alpha(\gamma(\alpha))={\mathcal C}_1(\gamma(1))=1.
\]
To avoid overflows in the values of $\widetilde{\tn}_n$ and $S_n$, we use the procedure described in~\cite[Section~5.2.1]{fort:jourdain:lelievre:stoltz:2015}.

\subsection{Asymptotic behavior of the  stepsize sequence}

\begin{figure}
\includegraphics[width=0.5\textwidth]{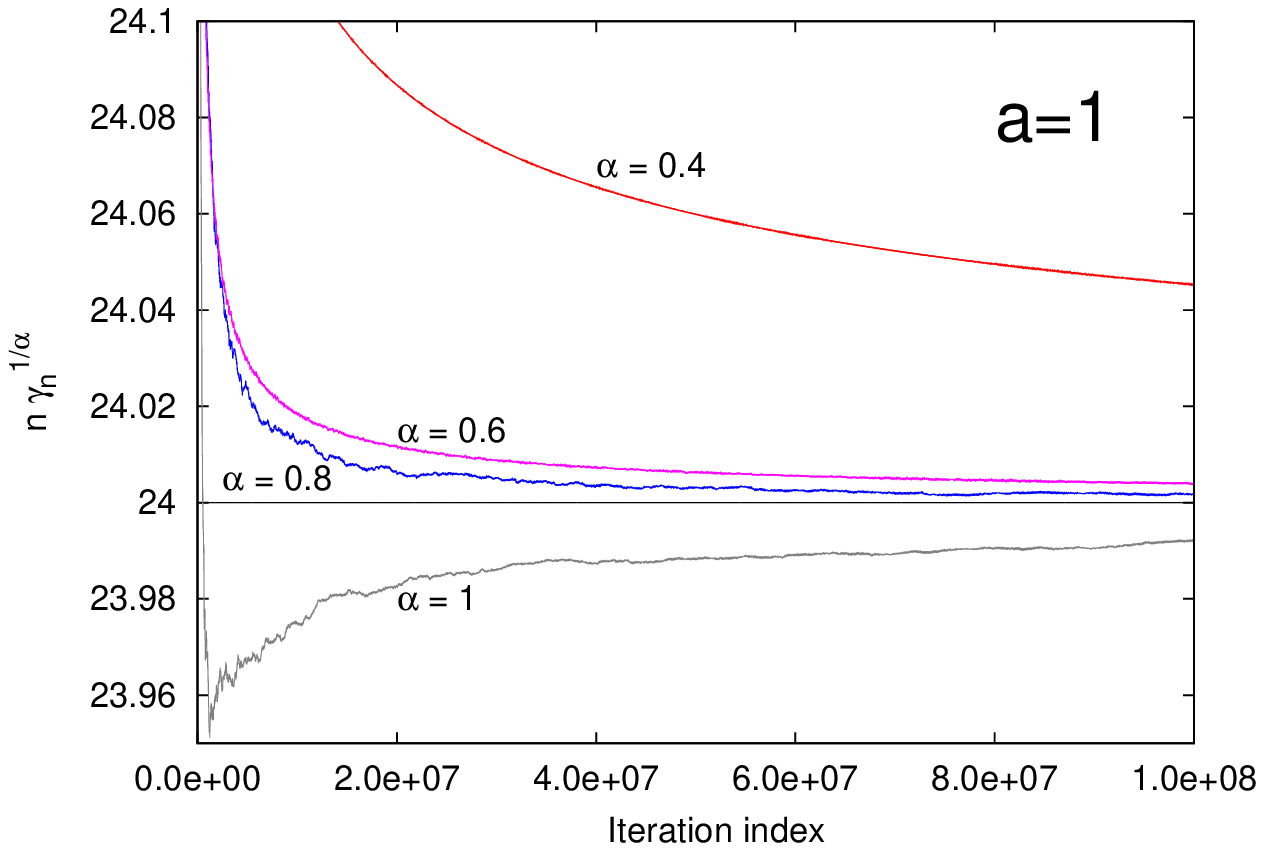}
\includegraphics[width=0.5\textwidth]{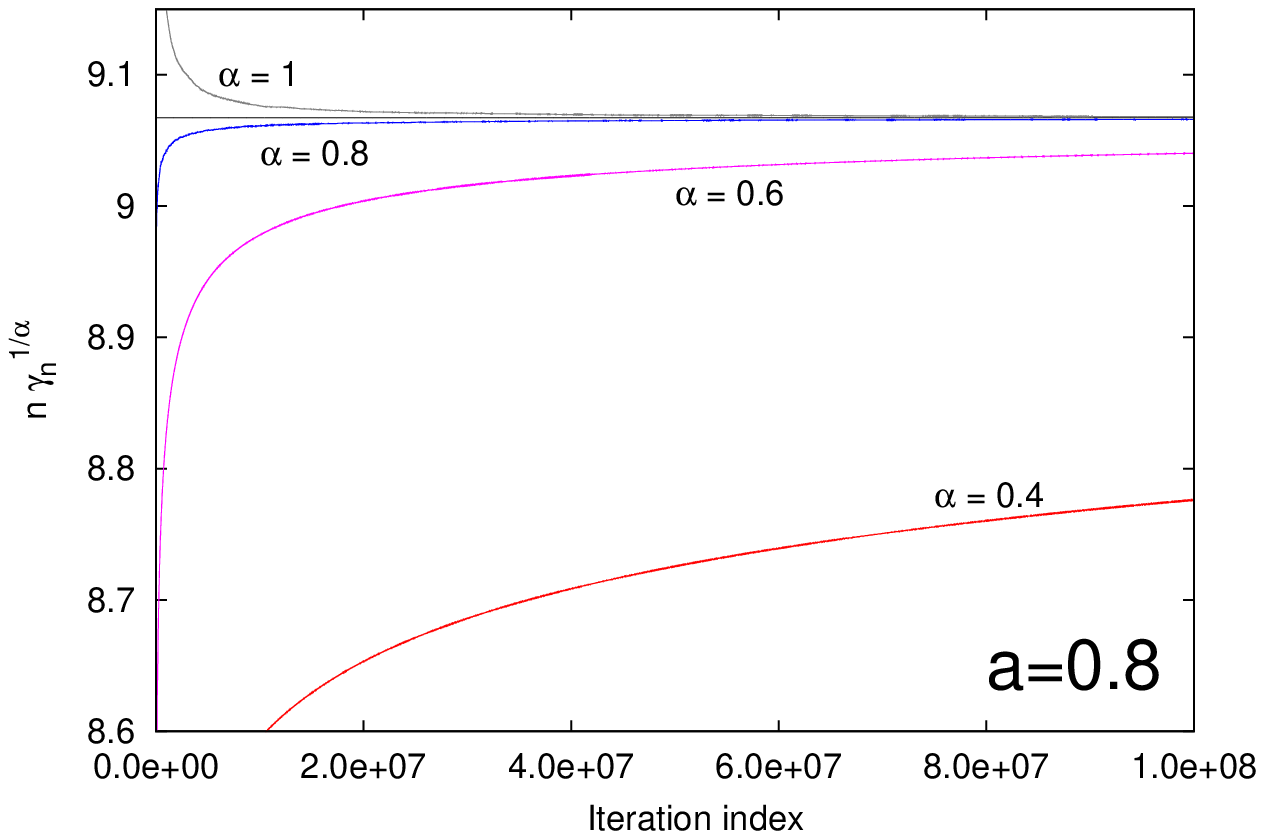}
\includegraphics[width=0.5\textwidth]{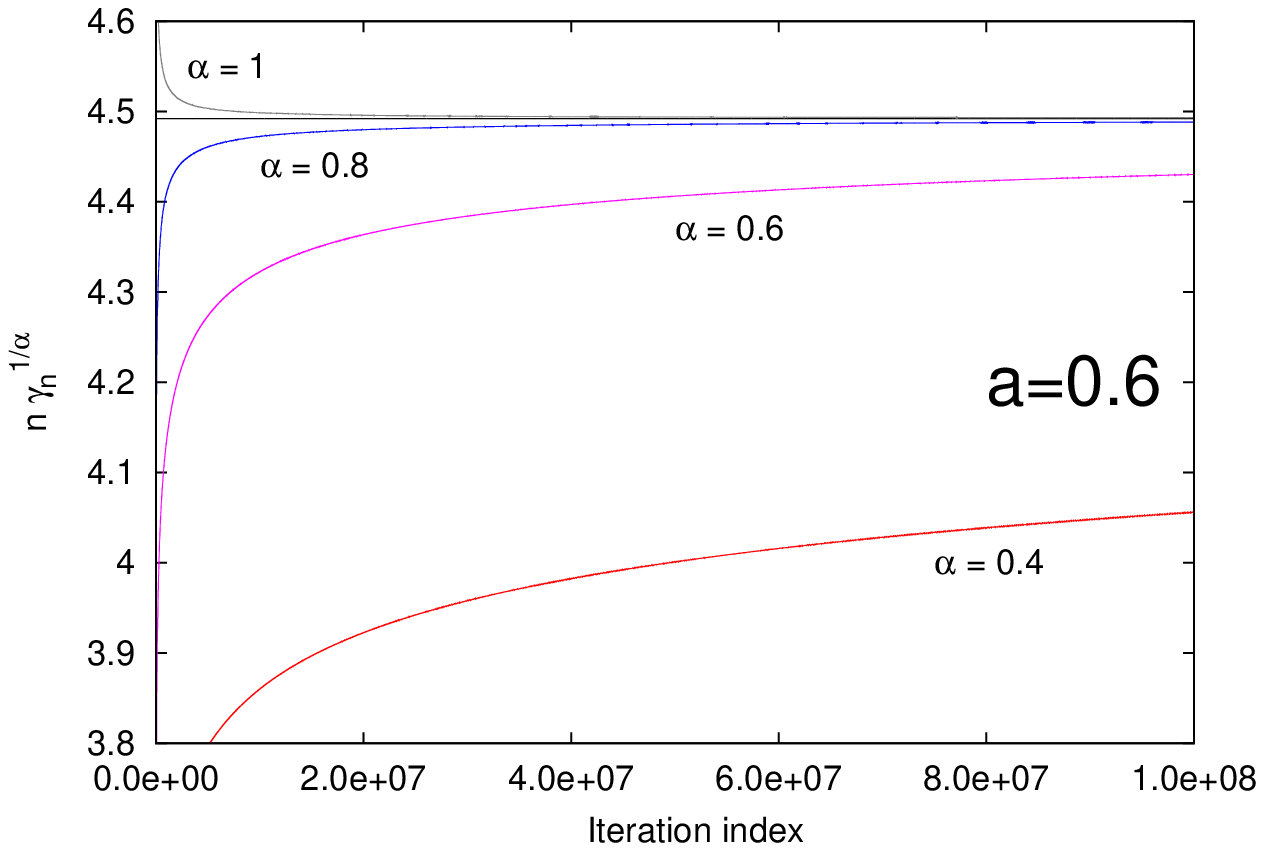}
\includegraphics[width=0.5\textwidth]{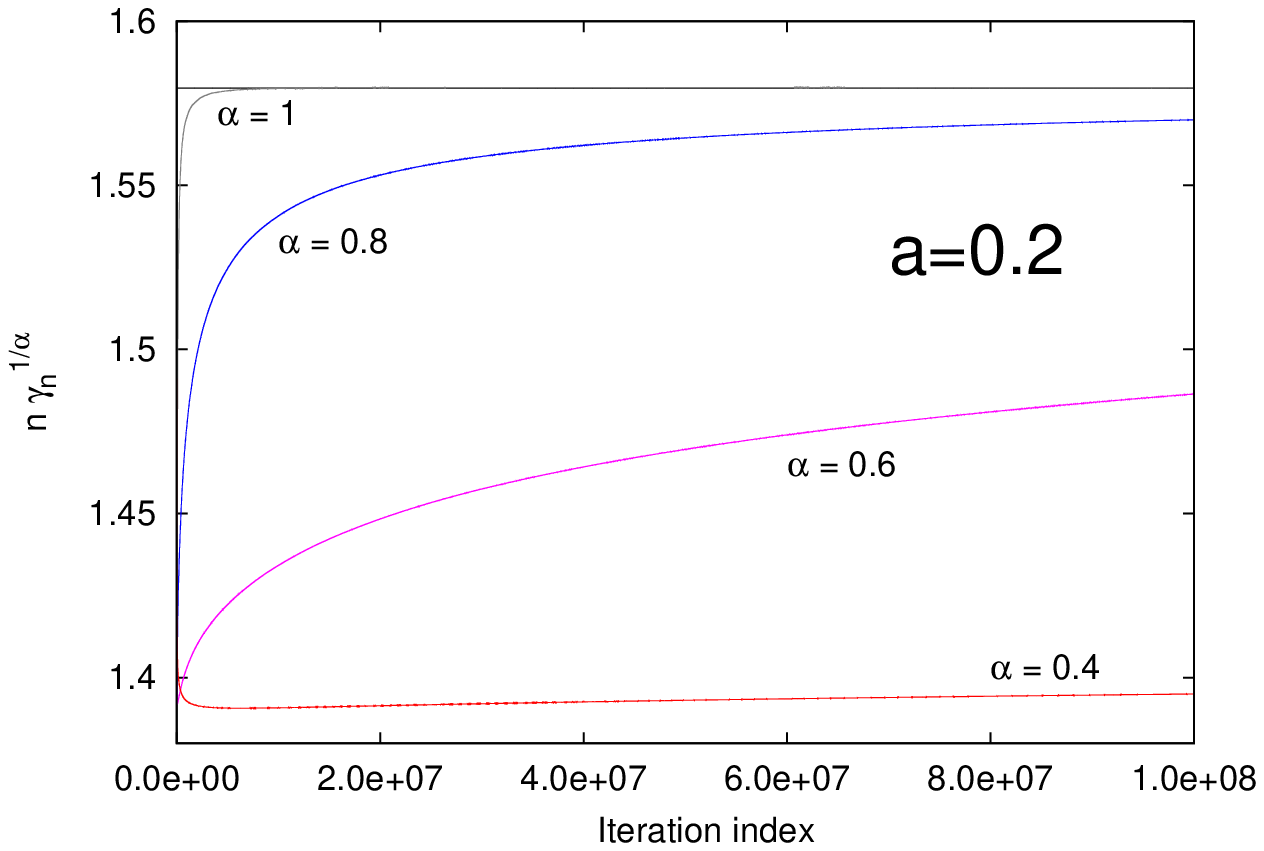}
\caption{\label{fig:scaling_steps} Behavior of the effective time step as a function of the
  iteration index, at $\beta = 4$ and averaged over 5000
  realizations. Each plot represents $n \gamma_n^{1/\alpha}$ as a function of the iteration index~$n$ for a different value of~$a$, the various
  curves representing the results obtained for $\alpha \in \{ 0.4, 0.6, 0.8, 1\}$. Note that, in all cases, the convergence is faster for larger values of~$\alpha$. The horizontal line represents the limiting value $\mathfrak{g}(a)$ (namely $\mathfrak{g}(1) = 24$, $\mathfrak{g}(0.8) = 9.07$, $\mathfrak{g}(0.6) = 4.49$ and $\mathfrak{g}(0.2) = 1.58$).
}
\end{figure}

\begin{figure}
\begin{center}
\includegraphics[width=0.5\textwidth]{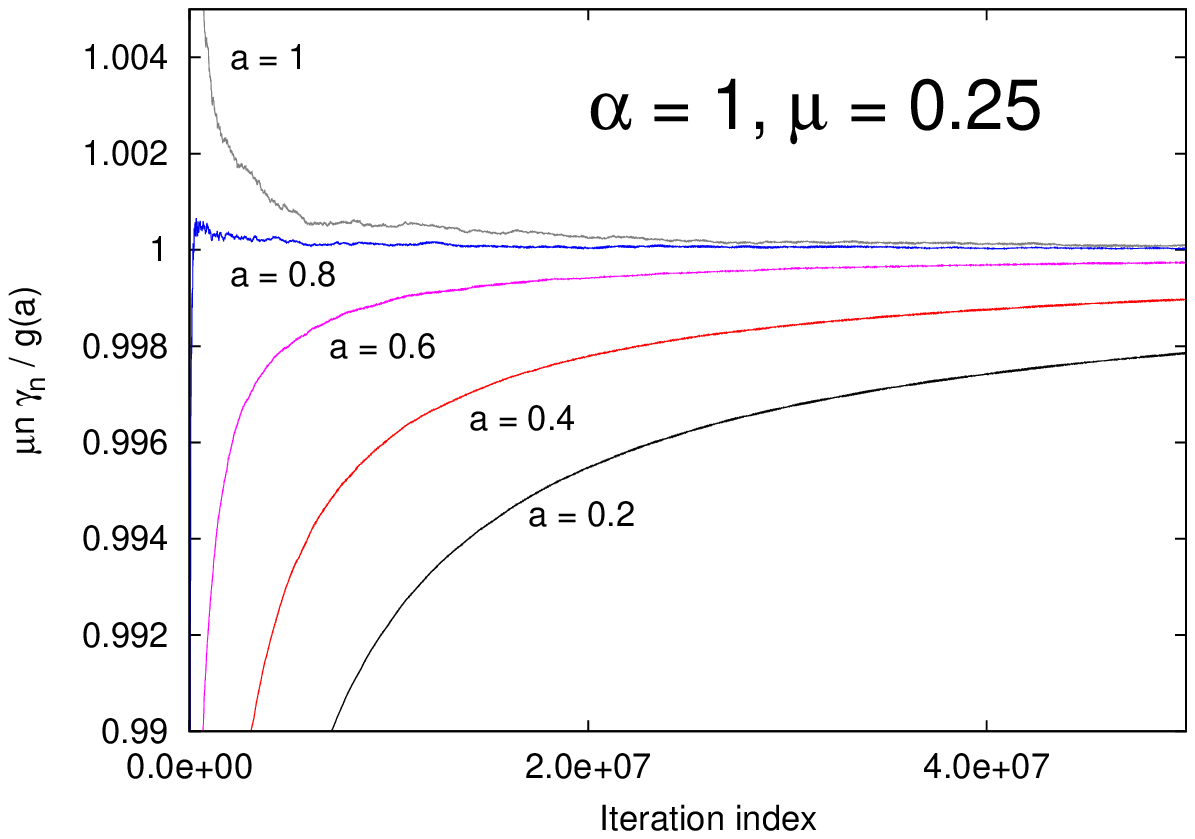}
\caption{\label{fig:scaling_steps_mu} Behavior of the effective time step as a function of the iteration index for $\alpha = 1$ and $\mu = 0.25$, at $\beta = 4$ and averaged over 5000 realizations. Each plot represents $\mu n \gamma_n / \mathfrak{g}(a)$ as a function of the iteration index~$n$ for a different value of~$a$, the various
  curves representing the results obtained for $\alpha \in \{ 0.4, 0.6, 0.8, 1\}$. The expected limit is~1 in all cases.
}
\end{center}
\end{figure}

With the choice~\eqref{eq:choice_alpha} of the parameter $\gamma$, it is expected from Corollary~\ref{cor:asymptotic} that, in the large $n$ limit and for $\alpha \in (1/2,1)$,
\[
\gamma_n \sim \left( \frac{\mathfrak{g}(a)}{n} \right)^\alpha \text{ with
}\mathfrak{g}(a) = Z_{\theta_\star}^{t^a} = \sum_{i=1}^d \theta_\star(i)^{1-a},
\]
while, for $\alpha = 1$,
\[
\mu \gamma_n \sim \frac{\mathfrak{g}(a)}{n}. 
\]
This is checked numerically on Figures~\ref{fig:scaling_steps} (for various values of $\alpha$, with $\mu = 1$ when $\alpha=1$) and~\ref{fig:scaling_steps_mu} (where the dependence on~$\mu$ for $\alpha=1$ is investigated). The reference values $\tn_\star(i)$ are obtained by a two-dimensional
quadrature. The numerical results are in excellent agreement with the
theoretical findings, and allow to go even beyond the theoretical predictions since the scaling of the steps also hold for values of $\alpha < 1/2$. Note however that the convergence is slower for smaller values of~$a$ and~$\alpha$. In particular, although not apparent on the scale of the plots, the values of $n \gamma_n^{1/\alpha}$ are very slowly increasing for $a=0.2$ and $\alpha=0.4$. 

\subsection{Exit times}
\label{sec:exit_num}

In order to show the interest of using the SHUS$_\rat^\alpha$
algorithm to get out of metastable states, we now study the exit time from the left metastable state in the small
temperature regime $\beta \to \infty$. More precisely, starting from the initial condition $X_0 = (-1,0)$ close to the
global minimum~$x_-$, we consider the time it takes to the system to go to the vicinity of the global minimum~$x_+$.

This section is organized as follows. In
Section~\ref{sec:strata}, exit times from metastable states for the SHUS$_\rat^\alpha$
algorithm are studied using a toy model with only three
states. Following previous results (see~\cite{amrx}),
this gives heuristic scalings for the exit times in the small
temperature regime on the two-dimensional model presented above. We then compare in
Section~\ref{sec:exit_numerics} these expected asymptotic behaviors with those
numerically observed on the two-dimensional potential.

\subsubsection{Heuristic on a toy model}\label{sec:strata}

Following~\cite{amrx}, one can derive a scaling of the time needed to
leave a metastable state on a toy model with only three states.

Let us recall the model we consider in~\cite{amrx}. The toy model
consists of
three states and three strata, so that $\Xset=\{1,2,3\}$ and
$\Xset_i=\{i\}$ for $i \in \{1,2,3\}$. The target probability is defined, for
a small positive parameter $\varepsilon \in (0,1)$, by
\begin{equation}
\label{eq:def_pi_toy_model}
\pi(\{1\})=\pi(\{3\})=\frac{1}{2+\varepsilon} \text{ and } \pi(\{2\})=\frac{\varepsilon}{2+\varepsilon}.
\end{equation}
Notice that in this setting, for all $i \in \{1,2,3\}$,
$\pi(\{i\})=\theta_\star(i)$. The unbiased dynamics is a
Metropolis-Hastings algorithm with target $\pi$ and proposal matrix $Q \in
\R^{3\times 3}$ defined by
$$Q=\left[\begin{array}{ccc} 2/3 & 1/3 & 0 \\ 1/3 & 1/3 & 1/3 \\ 0 & 1/3 &
  2/3 \end{array} \right],$$
so that jumps are only allowed between neighboring states. Since the
state $2$ has small probability, there are two metastable states for
this dynamics: state $1$ and state $3$. One can
check that the number of iterations $n^{\rm MH}(\varepsilon)$ needed
by the unbiased dynamics to go from state 1 to
state 3 is of order 
\begin{equation}\label{eq:n_MH}
n_{1 \to 3}^{\rm MH}(\varepsilon) \simeq \frac{6}{\varepsilon}
\end{equation}
in the limit $\varepsilon \to 0$ (see~\cite[Proposition 3.1]{amrx} for precise statements). 
One important feature of adaptive dynamics is that this exit time is drastically reduced
(compare~\eqref{eq:n_MH} and~\eqref{eq:n_alpha_1}-\eqref{eq:n_alpha_12} below).

Let us now consider the SHUS$^\alpha_\rho$ algorithm applied to this
toy model. Set $\tu_0=(\frac{1}{3},\frac{1}{3},\frac{1}{3})$ and
$X_0=1$. The main output of the paper~\cite{amrx} is that the number
of iterations needed to go from $1$ to $3$ using the adaptive
algorithm can be estimated in the limit $\varepsilon \to 0$ by
considering the minimum number of iterations required to reach the
metastable state $2$. This is what we estimate in the following, using
some formal derivations, which could be made rigorous using the same
techniques as in~\cite{amrx}. For the sake of conciseness, we do not
give here rigorous proofs of these results. We will however check the
consistency of these results with what is observed numerically in the
next section.

As long as $(X_n)_{n \ge 0}$ stays in
state $1$, $\tu_n(2)=\tu_0(2)=\frac 1 3$, $\tu_n(3)=\tu_0(3)=\frac 1
3$ and $\tu_n(1)=u_n$ where $(u_n)_{n \ge 0}$ is a sequence
satisfying, in view of~\eqref{eq:def:thetatilde}:
\begin{equation}\label{eq:un}
u_{n+1}=u_n+\gamma\frac{u_n+\frac{2}{3}}{g_\alpha\left(u_n+\frac{2}{3}\right)}
\left(\frac{3u_n}{3u_n+2}\right)^a,
\end{equation}
with initial condition $u_0=\frac13$. Let us denote by 
$$\zeta_n \eqdef \frac{(u_n,\frac13,\frac13)}{u_n+\frac23}$$
the associated normalized vector in $\Theta$.
Since $P^{t^a}_{\t}(1,2)=\frac 1 3 \left(\frac{\varepsilon
    \tn^a(1)}{\tn^a(2)}\wedge 1\right)=\frac 1 3 \left(\frac{\varepsilon
    \tu^a(1)}{\tu^a(2)}\wedge 1\right)$, it follows that $P^{t^a}_{\zeta_n}(1,2)=\frac 1 3
\left(\varepsilon (3u_n)^a\wedge 1\right)$.
Thus, the probability of staying in state $1$ from iteration $0$ up to
iteration $n+1$ is
$$\prod_{k=0}^n\left(
  1-P^{t^a}_{\zeta_k}(1,2)\right)=\exp\left(\sum_{k=0}^n\ln\bigg(1-\frac 1 3
  (\varepsilon (3u_k)^a\wedge 1)\bigg)\right).$$
For $\varepsilon$ small, expanding the logarithm,
this probability is of order $1/2$ when $\sum_{k=0}^nu_k^a$ is of
order $1/\varepsilon$.

\paragraph{The case $\alpha=1$.}
When $\alpha=1$, in view of~\eqref{eq:un}, 
$$u_{n+1}=u_n+\gamma\left(u_n+\frac{2}{3}\right)^{1-\mu} \left(\frac{3u_n}{3u_n+2}\right)^a.$$
When $n$ is large, $u_n$ becomes large and 
$$u_{n+1}\simeq u_n+\gamma u_n^{1-\mu}.$$
By
comparison with the ordinary differential equation
$\frac{dy}{dt}(t)=\gamma(y(t))^{1-\mu}$ with solution
$y(t)=\left(y(0)^{\mu}+\gamma\mu t\right)^{1/\mu}$, one obtains
$u_n\simeq (\gamma \mu n)^{1/\mu}$ and then $\sum_{k=0}^n u_k^a\simeq (\gamma \mu )^{a/\mu}
\frac{\mu}{\mu+a}n^{\frac{\mu+a}{\mu}}$. The probability to reach state $2$  for
the first time after at least 
$n+1$ iterations attains $1/2$ when this sum becomes of order $1/\varepsilon$ i.e. when $n$
is of order $\gamma^{-\frac{a}{\mu+a}}\mu^{-1}(\mu+a)^{\frac{\mu}{\mu+a}}\varepsilon^{-\frac{\mu}{\mu+a}}$.  We thus obtain that the
number of iterations $n^1_{1 \to 3}$ needed to go from $1$ to $3$ in the case
$\alpha=1$ satisfies, in the limit $\varepsilon \to 0$,
\begin{equation}\label{eq:n_alpha_1}
n^1_{1 \to 3}(\varepsilon) \simeq C^1_{1 \to 3} \, \varepsilon^{-\frac{\mu}{\mu+a}}.
\end{equation}
where $C^1_{1 \to 3} \eqdef \gamma^{-\frac{a}{\mu+a}}\mu^{-1}(\mu+a)^{\frac{\mu}{\mu+a}}$ is a constant.

\paragraph{The case $\alpha \in (\frac12,1)$.}
When $\alpha\in (\frac12,1)$, in view of~\eqref{eq:un}, 
\[
u_{n+1}=u_n+\gamma \frac{u_n+2/3}{\ln^{\frac{\alpha}{1-\alpha}}(1+u_n+2/3)} \left(\frac{3u_n}{3u_n+2}\right)^a.
\]
When $n$ is large, $u_n$ becomes large and 
$$u_{n+1}\simeq u_n\left(1+\frac{\gamma}{\ln^{\frac{\alpha}{1-\alpha}}(u_n)}\right),$$
so that
$$\ln\left(u_{n+1}\right)\simeq
\ln\left(u_n\right)+\frac{\gamma}{\ln^{\frac{\alpha}{1-\alpha}}(u_n)}.$$
By
comparison with the ordinary differential equation
$\frac{dy}{dt}(t)=\gamma(y(t))^{-\frac{\alpha}{1-\alpha}}$ with solution
$y(t)=\left(y(0)^{1/(1-\alpha)}+\frac{\gamma t}{1-\alpha}\right)^{1-\alpha}$, one obtains
$u_n\simeq \exp\left[ \left(\frac{\gamma
      n}{1-\alpha}\right)^{1-\alpha}\right]$. 
We now want to
estimate $\sum_{k=0}^n u^a_k$. Since for $c=a\left(\frac{\gamma
  }{1-\alpha}\right)^{1-\alpha}$,
$$\int_0^x\mathrm{e}^{cy^{1-\alpha}} \rmd y=\frac{x^\alpha
  \mathrm{e}^{cx^{1-\alpha}}}{c(1-\alpha)}-\frac{\alpha}{c(1-\alpha)}\int_0^xy^{\alpha-1}\mathrm{e}^{cy^{1-\alpha}}\rmd
y,$$
$\int_0^x\mathrm{e}^{cy^{1-\alpha}} \rmd y\sim \frac{x^\alpha}{c(1-\alpha)}\mathrm{e}^{cx^{1-\alpha}}$
as $x\to +\infty$. Hence $\sum_{k=0}^nu^a_k\sim \frac{n^{\alpha} \mathrm{e}^{c
    n^{1-\alpha}}}{c(1-\alpha)}$. The probability to reach state $2$  for
the first time after at least 
$n+1$ iterations attains $1/2$ when this sum becomes of order $1/\varepsilon$ i.e. when $n$
is of order $\frac{1-\alpha}{\gamma
}\left(-\frac{\ln\varepsilon}{a}\right)^{1/(1-\alpha)}$.  We thus obtain that the
number of iterations $n^\alpha_{1 \to 3}$ needed to go from $1$ to $3$ in the case
$\alpha\in (1/2,1)$ satisfies, in the limit $\varepsilon \to 0$,
\begin{equation}\label{eq:n_alpha_12}
n_{1 \to 3}^{\alpha}(\varepsilon) \simeq C^\alpha_{1 \to 3} \,
\left|\ln \varepsilon\right|^{1/(1-\alpha)}.
\end{equation}
where $C^\alpha_{1 \to 3} \eqdef \frac{1-\alpha}{\gamma
} a^{1/(\alpha-1)}$.

By comparing the transition time from $1$ to $3$ on the 
unbiased Metropolis-Hastings dynamics (see~\eqref{eq:n_MH}), and for the
SHUS$^\alpha_\rho$ algoritm (see~\eqref{eq:n_alpha_1}
and~\eqref{eq:n_alpha_12}), the interest of using the adaptively
biased dynamics is obvious: the exit times are much smaller for
SHUS. Moreover, one can see that it is interesting to consider $\alpha
< 1$ in order to reduce the exit time compared to $\alpha=1$.

\subsubsection{Numerical results on exit times}
\label{sec:exit_numerics}

Average exit times for the two-dimensional model described at the beginning of Section~\ref{sec:application}
are obtained by performing independent realizations of the
following procedure, for given values of~$a,\alpha,\beta$: initialize the system in the state $X_0=(-1,0)$, and run
the dynamics until the first time index $\mathcal{N}$ such that
the first component of $X_\mathcal{N}$
is larger than~1. We perform $K$ independent
realizations of this process. The corresponding empirical average
first exit time is denoted by $t^{\alpha}_\beta$. Since we work with a fixed maximal computational time 
(of about a week or two on our computing machines with our implementation of the code),
$K$ turns out to be of the order of a few hundreds for the largest exit times, while $K=10^5$ in the
easiest cases corresponding to the shortest exit times. In our
numerical results, we checked that $K$ is always sufficiently large so that the relative error on~$t_\beta^\alpha$ is less than a
few percents in the worst cases.

The first task is to identify the equivalent of the parameter~$\varepsilon$ in the toy model from Section~\ref{sec:strata}.
In the large $\beta$ regime, using Laplace's method, the ratio between the
probability of the stratum in the transition region (around the vertical axis $(0,y)$ for $y \in \mathbb{R}$) and the
metastable states is of order $\bar{C} \exp(-\beta \delta_0)$ for some positive
constants $\bar{C}$ and $\delta_0$. In view of~\eqref{eq:def_pi_toy_model}, this suggests the following formal equivalence
between $\varepsilon$ and $\beta$:
\begin{equation}\label{eq:eps_beta}
\varepsilon(\beta)=\bar{C}\mathrm{e}^{-\beta \delta_0}.
\end{equation}
We next replace $\varepsilon$ by the right-hand side of the above equality in the scalings found at the end of Section~\ref{sec:strata}.

When $\alpha=1$, it is expected from~\eqref{eq:eps_beta} and~\eqref{eq:n_alpha_1} that, in the regime $\beta \to \infty$, 
\begin{equation}\label{eq:t_alpha_1}
t^{1}_{\beta} \simeq \widetilde{C} \mathrm{e}^{\beta\delta_0\frac{\mu}{\mu+a}}
\end{equation}
for some constant $\widetilde{C}$. When $\alpha \in (1/2,1)$, it is expected from~\eqref{eq:eps_beta} and~\eqref{eq:n_alpha_12} that, in the regime $\beta \to
\infty$, 
\begin{equation}\label{eq:t_alpha_12}
\ln(t^\alpha_\beta)\simeq\frac{1}{1-\alpha}\ln\beta.
\end{equation}
We check in the sequel the scalings~\eqref{eq:t_alpha_1}-\eqref{eq:t_alpha_12} by varying the parameters of the dynamics in several ways: (i) fix $\alpha$ (as well as $\mu$ for $\alpha=1$) and vary~$a$; (ii) fix $\alpha=1$ and $a$, and vary~$\mu$; (iii) fix $\alpha=1$ and vary $a = 1-\mu$, which corresponds to the well-tempered metadynamics.

\begin{figure}
\includegraphics[width=0.5\textwidth]{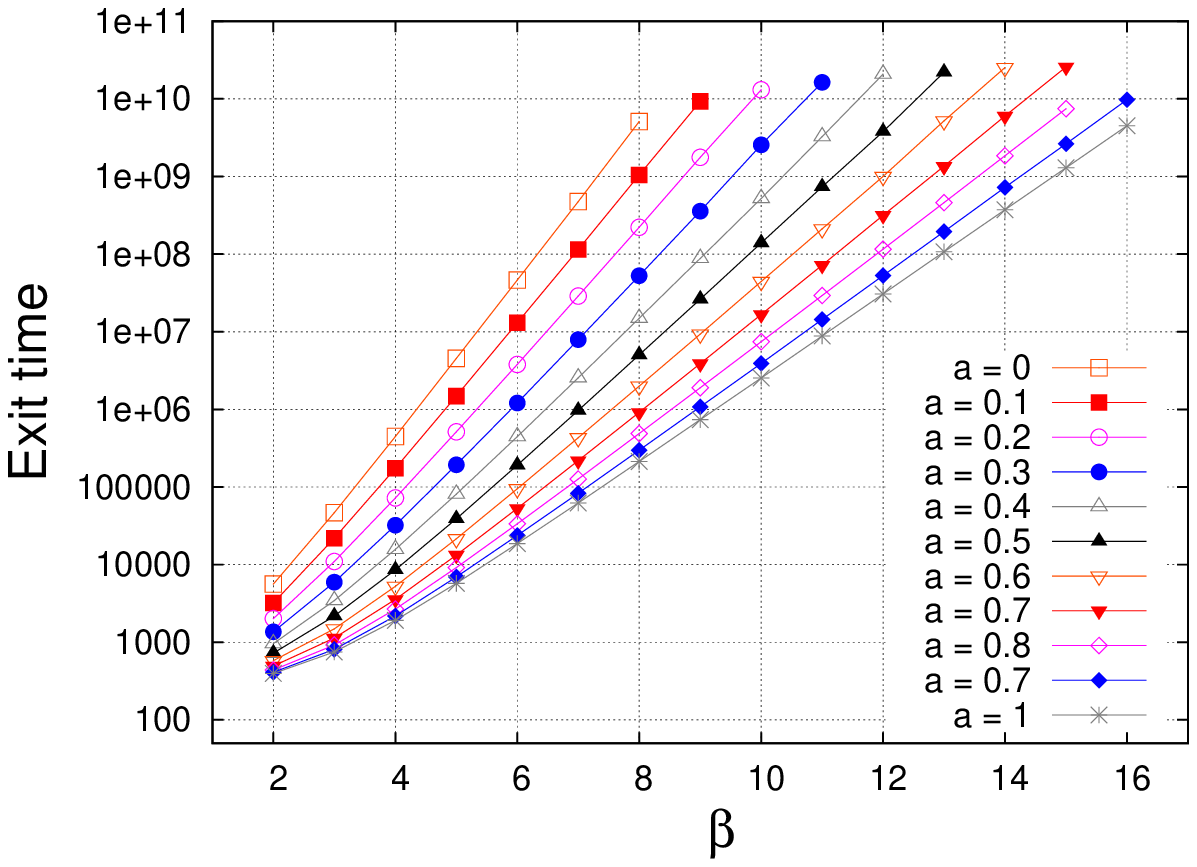}
\includegraphics[width=0.5\textwidth]{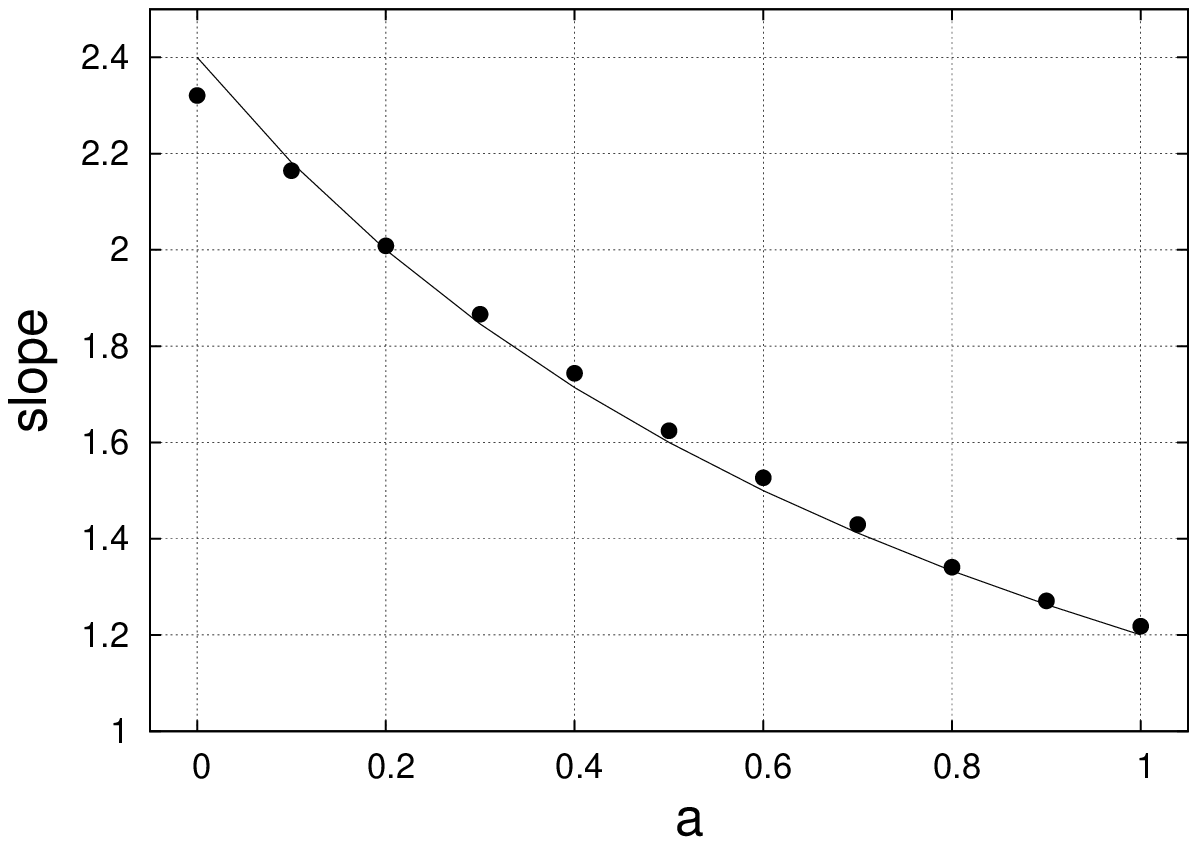}
\caption{Case $\alpha=1$. Left: Exit times for  various values of $a$. Right: Associated slopes $r(a)$, fitted by $2.4/(1+a)$.}
\label{fig:exit_1}
\end{figure}

\begin{figure}
  \includegraphics[width=0.5\textwidth]{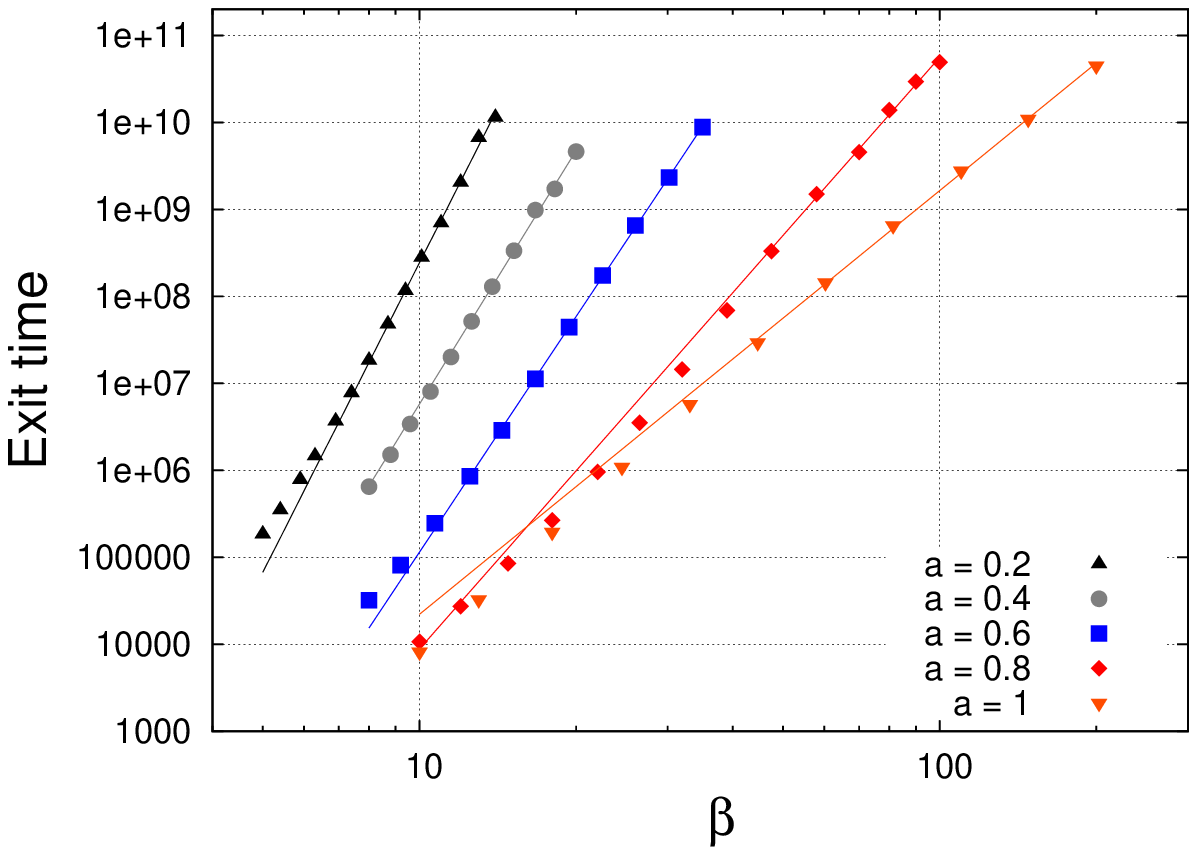}
   \includegraphics[width=0.5\textwidth]{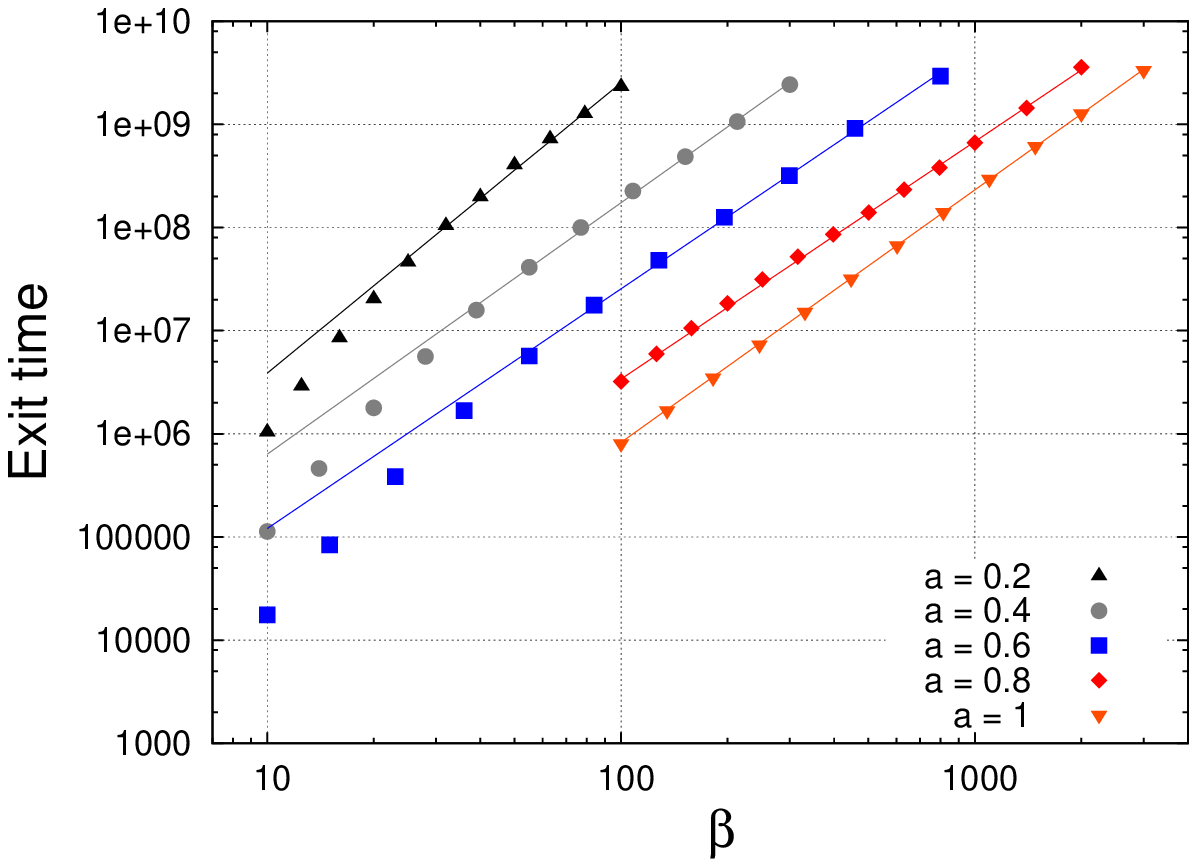}
\caption{Case $\alpha \in (1/2,1)$. Left: Exit times for $\alpha = 0.8$. Right: Exit times for $\alpha = 0.6$. In both cases, a linear fit is superimposed in solid lines to the data.}
\label{fig:exit_alpha}
\end{figure}

\paragraph{Dependence on~$a$.}
The dependence of the exit times on~$a$ is studied in Figure~\ref{fig:exit_1} for $\alpha=\mu=1$, and in Figure~\ref{fig:exit_alpha} for $\alpha\in (\frac12,1)$.

For $\alpha=\mu=1$, we perform, for each value of~$a$ a least-square fit of $\ln t^1_\beta$ in terms of $\beta$ to obtain the scaling $t^1_\beta \sim \mathrm{e}^{\beta r(a)}$. In view of~\eqref{eq:t_alpha_1}, we next compare $r(a)$ with $\delta_0/(1+a)$. Numerically, we estimate $\delta_0 \sim 2.4$, in accordance with the value found in~\cite{amrx}. The numerical results on Figure~\ref{fig:exit_1} are therefore in excellent agreement with the scaling expected from the toy model. 

For $\alpha \in (\frac12,1)$, we perform a fit of $\ln(t^\alpha_\beta)$ in terms of $\ln\beta$ to obtain the scaling
$\ln(t^\alpha_\beta) \sim s(\alpha,a) \ln \beta$. In view of~\eqref{eq:t_alpha_12}, 
we expect $s(\alpha,a)$ to be independent of~$a$ and close to $1/(1-\alpha)$. The slopes $s(\alpha,a)$ obtained in the numerical experiments displayed on 
Figure~\ref{fig:exit_alpha} for two values of $\alpha$ (namely $\alpha=0.6$ and
$\alpha=0.8$) are reported in Table~\ref{tab:exit_alpha}. They are 
close to the expected values for $\alpha=0.6$, as well as for
$\alpha=0.8$ when $a$ is close to $1$. The discrepancies observed for $\alpha=0.8$
and small values of $a$ may be due to the fact it is difficult to
reach the asymptotic regime $\beta \to \infty$ in this setting. It might be that the slopes of the curves would decrease for much larger values of $\beta$. 

\begin{table}
\begin{center}
  \begin{tabular}{|c||c c c c c | c|}
    \hline
& $a=0.2$ & $a=0.4$ & $a=0.6$ & $a=0.8$ & $a=1$ & expected \\
\hline
\hline
$\alpha=0.6$ & 2.72 & 2.38 & 2.27 & 2.33 & 2.44 & 2.5 \\
$\alpha=0.8$ & 11.8 & 9.69 & 9.00 & 6.80 & 4.87 & 5 \\
\hline
\end{tabular}
\end{center}
\caption{Numerically estimated slopes $s(\alpha,a)$ of the exit times $t^\alpha_\beta \sim
  \beta^{s(\alpha,a)}$, for various values of $a$ (see Figure~\ref{fig:exit_alpha}). The expected values
  are $\displaystyle s(\alpha) = 1/(1-\alpha)$, namely $s(0.6)=2.5$
  and $s(0.8)=5$.}\label{tab:exit_alpha}
\end{table}

\paragraph{Dependence on~$\mu$.}

\begin{figure}
\includegraphics[width=0.5\textwidth]{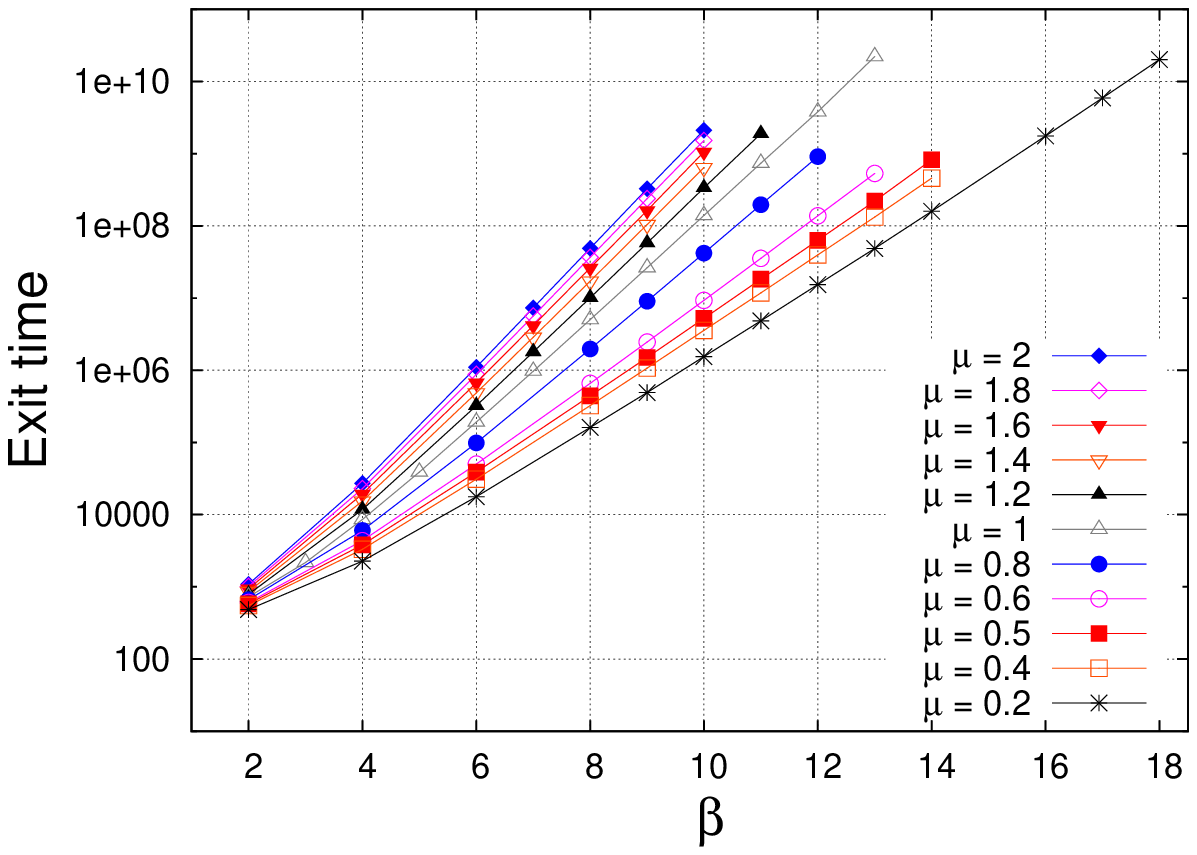}
\includegraphics[width=0.5\textwidth]{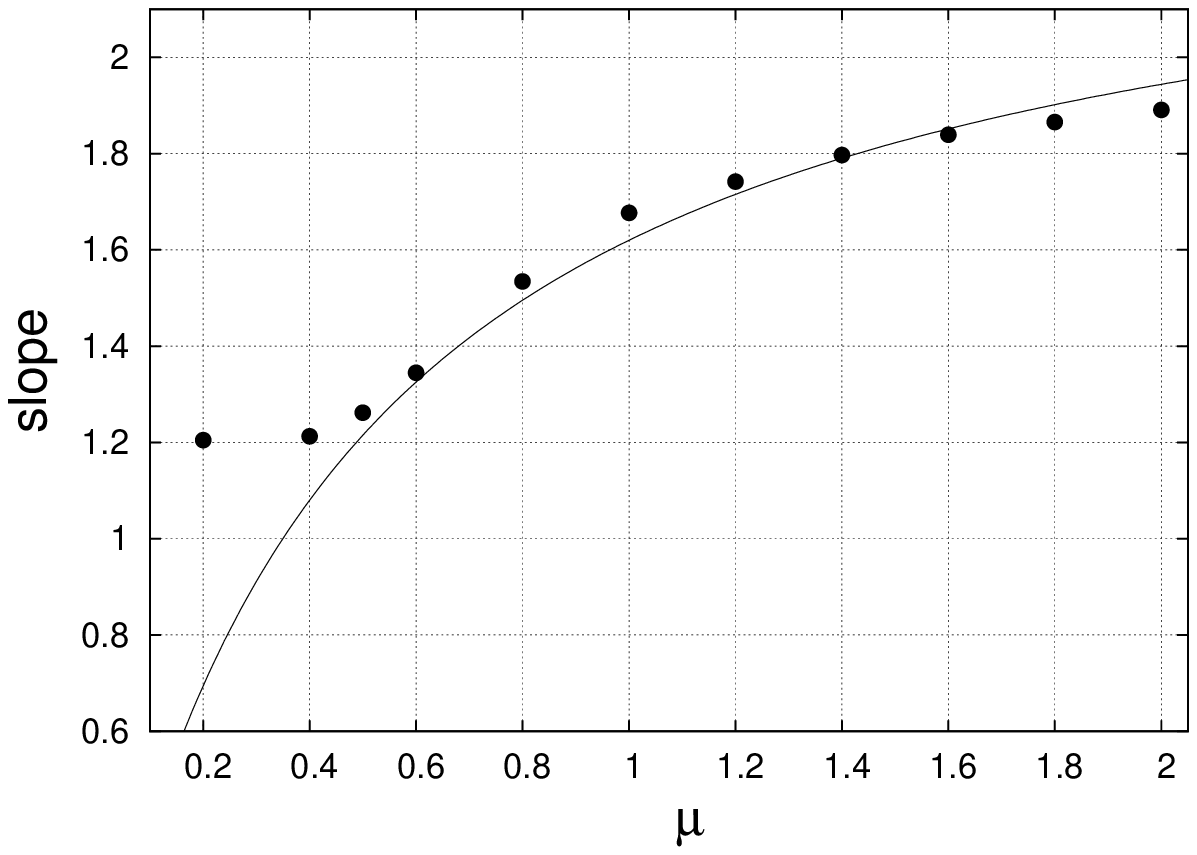}
\caption{Case $\alpha=1$ and $a = 0.5$. Left: Exit times for  various values of $\mu$. Right: Associated slopes $r(\mu)$, fitted by $2.3\mu/(\mu+0.5)$.}
\label{fig:exit_mu_varies}
\end{figure}

We next consider the case when $\alpha=1$ and $a=0.5$ are fixed, and compute scaling times for various values of~$\mu$; see Figure~\ref{fig:exit_mu_varies}. We again perform, for each value of~$a$ a least-square fit of $\ln t^1_\beta$ in terms of $\beta$ to obtain the scaling $t^1_\beta \sim \mathrm{e}^{\beta r(\mu)}$. In view of~\eqref{eq:t_alpha_1}, we next compare $r(\mu)$ with $\delta_0 \mu/(\mu+0.5)$. Numerically, we estimate $\delta_0 \sim 2.3$, in accordance with the value found in~\cite{amrx}. Note however that the predicted slopes $r(\mu)$ are somewhat off the prediction $\delta_0 \mu/(\mu+0.5)$ for small values of~$\mu$. This might be due to the fact that, as in the case $\alpha=0.8$ and $a$ small, the asymptotic regime has not yet been reached.

\paragraph{Well-tempered metadynamics.}

\begin{figure}
\includegraphics[width=0.5\textwidth]{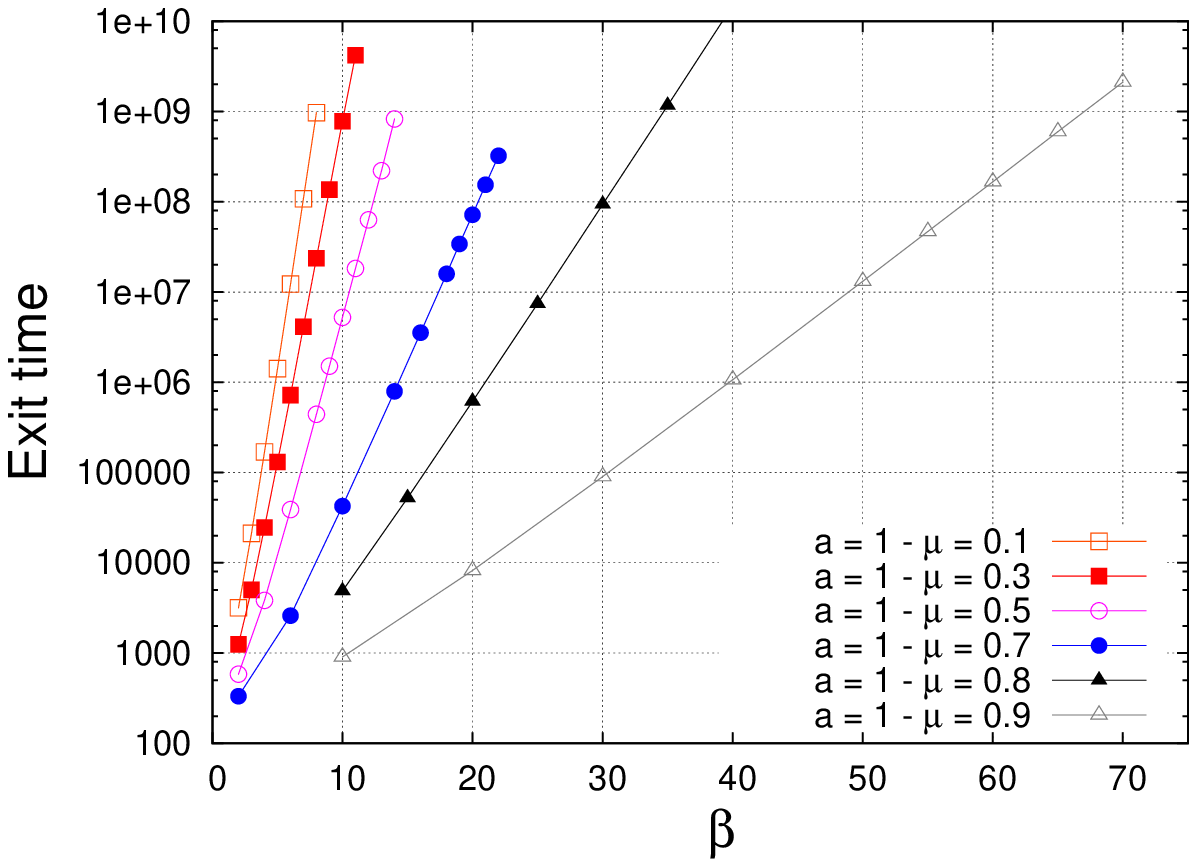}
\includegraphics[width=0.5\textwidth]{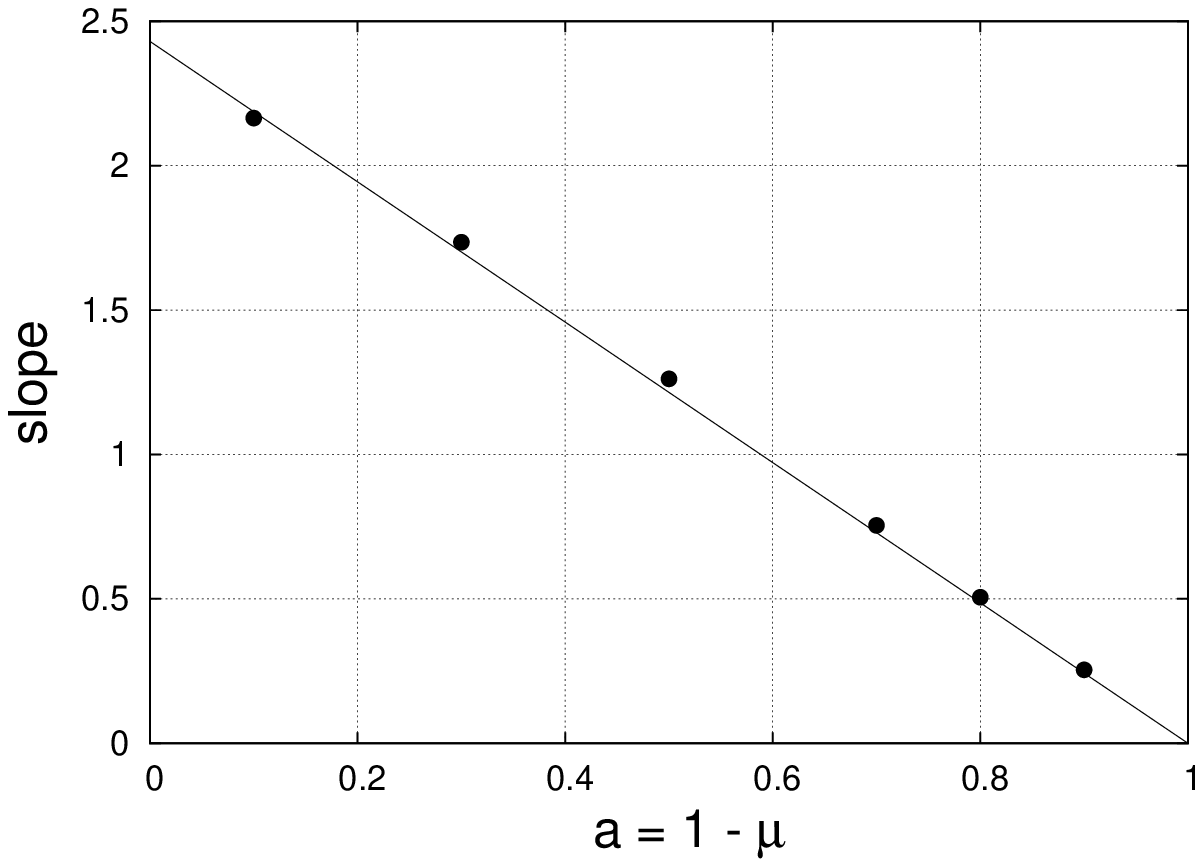}
\caption{Case $\alpha=1$ and $a = 1-\mu$ (original well-tempered dynamics). Left: Exit times for  various values of $a = 1-\mu$. Right: Associated slopes $r(a)$, fitted by $2.43(1-a)$.}
\label{fig:exit_WT_original}
\end{figure}

We finally turn to the well-tempered metadynamics case, which corresponds to $\alpha=1$ and varying values of $a = 1-\mu$. Exit times in this setting are reported in Figure~\ref{fig:exit_WT_original}, together with the slopes $r(a)$ fitted on the data as $t^1_\beta \sim \mathrm{e}^{\beta r(a)}$. The so-obtained values are in excellent agreement with the theoretical prediction $\delta_0(1-a)$ from~\eqref{eq:t_alpha_1}, for the choice $\delta_0 = 2.43$.

\subsection{Discussion on the effective sample size}
\label{sec:ESS}

As can be seen from the results, exit
times are drastically reduced as $\alpha$ decreases and $a$ increases. However,
the question arises whether the importance sampling strategy remains efficient
for large values of~$a$, for which the bias is larger hence the effective
sample size is smaller. 

Let us first recall the definition of the effective sample size. At
convergence (namely when the sequence $(\theta_n)_{n \ge 0}$ has reached its limiting value
$\theta_\star$), the weight of a sample $X$ is
$w(X)=\theta_\star(I(X))^a$. The effective sample size (ESS) of a
weighted ensemble of
$n$ i.i.d. samples $X_1, \ldots , X_n$ is defined as (see~\cite{KongLiuWong}) 
\begin{equation}\label{eq:ESS}
\mathrm{ESS}=\frac{\dps \left(\sum_{i=1}^n w(X_i)\right)^2}{\dps\sum_{i=1}^n w^2(X_i)}.
\end{equation}
The ESS is a real number in $[0,n]$. The more uniform the weights of the samples are, 
the closer to $n$ the ESS is. In order to normalize
this quantity, let us introduce the efficiency factor (EF) which is the
ESS divided by the number of samples:
\begin{equation}\label{eq:EF}
\mathrm{EF}=\frac{\dps \left(\sum_{i=1}^n w(X_i)\right)^2}{\dps n\sum_{i=1}^n w^2(X_i)}.
\end{equation}
The EF is close to one (respectively to zero) when the random variable $w(X)$
has a small (respectively a large) variance.

\begin{figure}
\begin{center}
\includegraphics[width=0.5\textwidth]{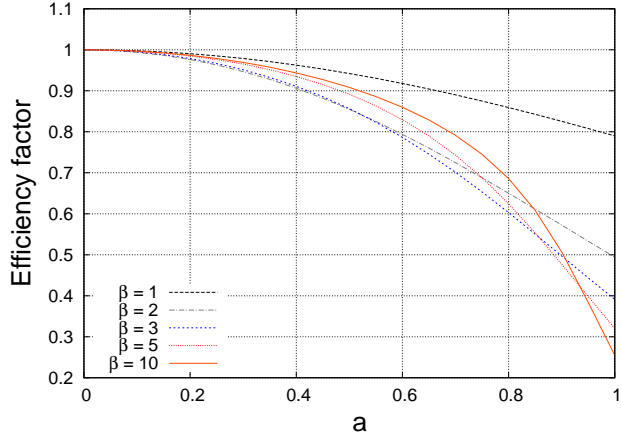}
\end{center}
\caption{Efficiency factors EF($a$) for various values of $\beta$.}
\label{fig:ESS}
\end{figure}

Following the strategy outlined
in~\cite[Section~4.1.2]{chopin-lelievre-stoltz-12}, it is possible to
give in our context the limit $\mathrm{EF}(a)$ of the efficiency factor as $n \to \infty$.
Indeed, at equilibrium, the samples are distributed according
to $\pi^\rho_{\theta_\star}$, where, we recall $\rho(t)=t^a$.
The probability of the $j$-th strata is thus $p_a(j) =
\theta_\star(j)^{1-a}/Z_a$ with $Z_a=\sum_{i=1}^d\theta_\star(i)^{1-a}$, the weight being
$w_a(j) = \theta_\star(j)^a$ in this region. Therefore,
\[
\mathrm{EF}(a) = \frac{\left(\sum_{j=1}^d p_a(j) w_a(j) \right)^2 }{\left(\sum_{j=1}^d p_a(j) \right)\left(\sum_{j=1}^d p_a(j) w^2_a(j) \right)} = \frac{1}{\left(\sum_{j=1}^d \theta_\star(j)^{1-a} \right)\left(\sum_{j=1}^d \theta_\star(j)^{1+a} \right)}.
\]
These functions are plotted in Figure~\ref{fig:ESS} for various values of
$\beta$. As expected, the efficiency factor decreases as $a$ increases.

\section{Proofs}\label{sec:proofs}

Throughout this section, $|\cdot|$ will be used both to denote the absolute
value and the Euclidean norm in $\R^p$. For two real numbers $a$ and
$b$, we denote $a \vee b = \max(a,b)$ and $a \wedge b=\min(a,b)$. Let us recall the
definition~\eqref{eq:Fn} of the filtration $(\mathcal{F}_n)_{n \ge 0}$:
\[
\mathcal{F}_n = \sigma\left(\tu_0, X_0, X_1, \cdots, X_n \right).
\]

The constants $C$ appearing in the proofs are finite constants which may change from line to line.

\subsection{Proof of Lemma~\ref{lem:SHUSalpharho:AS}}
\label{sec:proof:lem:SHUSalpharho:AS}
Let us rewrite~\eqref{eq:evol_algo_gen} as:
$$\tu_{n+1}(i)=\tu_n(i) + \gamma_{n+1} \, S_n \, \rho(\tn_n(i)) 1_{\Xset_i}(X_{n+1}).$$
By summing over $i=1,\ldots,d$, we have
\begin{equation}
S_{n+1}=S_n\left[1+
  \gamma_{n+1} \rat\left(\tn_n(I(X_{n+1}))\right)\right]\label{evosn1}.
\end{equation}
This implies that the sequence $(\tn_n)_{n \ge 0}$ evolves according to
\begin{equation}\label{eq:AlgoMultiplicatif}
\tn_{n+1}(i)= \tn_n(i) \  \left(1+\gamma_{n+1}\frac{\rat(\tn_n(i))}{\tn_n(i)}
  \un_{\Xset_i}(X_{n+1}) \right) \left[1+
  \gamma_{n+1}\rat(\tn_n(I(X_{n+1}))\right]^{-1}.
\end{equation}
By using $(1+a)/(1+b) = 1+ a-b + b(b-a)/(1+b)$, we have
\begin{align*}
  \tn_{n+1}(i)= \tn_n(i) & +\gamma_{n+1}\, \rat(\tn_n(i))
  \un_{\Xset_i}(X_{n+1}) -
  \gamma_{n+1} \tn_n(i) \, \rat(\tn_n(I(X_{n+1}))) \\
&  + \gamma_{n+1} ^2 \rat(\tn_n(I(X_{n+1}))) \frac{\tn_n(i)\,
    \rat(\tn_n(I(X_{n+1}))) - \rat(\tn_n(i)) \un_{\Xset_i}(X_{n+1})}{1 + \gamma_{n+1}
    \rat(\tn_n(I(X_{n+1}))) }.
\end{align*}
This can be reformulated as $\tn_{n+1}(i)= \tn_n(i) + \gamma_{n+1}
H_i(X_{n+1},\tn_n)+\gamma_{n+1}\Lambda_{n+1}(i)$, where the $i$-th component of
$\Lambda_{n+1}$ is given by
\begin{equation}
  \label{eq:DefinitionLambda}
  \Lambda_{n+1}(i)  \eqdef  \gamma_{n+1} \, \frac{ \left[ \rat(\tn_n(I(X_{n+1}))) \right]^2}{1+ \gamma_{n+1} \, \rat(\tn_n(I(X_{n+1})))} \ 
\left( \tn_n(i)-\un_{\Xset_i}(X_{n+1}) \right),
\end{equation}
Since $0 \leq \tn_n(i) \leq 1$ for all $i=1,\dots,d$, it holds
\[
\sum_{i=1}^d \left|\theta_n(i) - \un_{\Xset_i}(X_{n+1})\right|^2 = 1 + \sum_{i=1}^d \theta_n(i)^2 - 2 \sum_{i=1}^d \theta_n(i) \un_{\Xset_i}(X_{n+1}) \leq 1 + \sum_{i=1}^d \theta_n(i) = 2.
\]
It therefore follows that $|\Lambda_{n+1}| \leq \gamma_{n+1} \, \sqrt{2} \, \left( \sup_{(0,1)} \rat^2
\right)$.

\subsection{Proof of Proposition~\ref{prop:cvggal}}
The proof of Proposition~\ref{prop:cvggal} is performed by extending the
technique of proof used in~\cite{fort:jourdain:kuhn:lelievre:stoltz:2014} for
the convergence of the Wang-Landau algorithm to a {\em random} sequence of
stepsizes $(\gamma_n)_{n \ge 0}$.

Let us first recall the result \cite[Proposition
3.1]{fort:jourdain:kuhn:lelievre:stoltz:2014} on the Metropolis-Hastings
transition kernel $P^\rat_\t$, which will be useful below.
\begin{aprop}
  \label{prop:unifergo}
  Under A\ref{hyp:targetpi} and A\ref{hyp:kernel}, there exists $\delta \in
  (0,1)$ such that for all $\t \in \Theta$, for all $x \in \Xset$ and for all
  measurable set $A \subset \Xset$, it holds:
  \begin{align} 
    & P^\rat_\t(x,A) \geq \delta \, \int_A
    \pi^\rat_\t(x) \, \rmd \lambda(x) \eqsp, \label{eq:minorization_bound} \\
    & \forall n \ge 0, \, \sup_{\t \in \Theta} \sup_{x \in \Xset} \left\| (P^\rat_\t)^n(x, \cdot) -
      \pi^\rat_\t \, \rmd \lambda \right\|_\tv \leq 2 (1-\delta)^n,
    \label{eq:uniformeergodicite}
  \end{align}
where for a signed measure $\mu$, the total variation norm
is defined as
\[
\|\mu\|_\tv \eqdef \sup_{\{f \, : \, \sup_\Xset |f| \leq 1 \}} |\mu(f)| \eqsp.
\]
\end{aprop}
The proof is now organized as follows. We first state three lemmas which
quantify the dependence on $\theta$ of the invariant measure
$\pi^\rat_{\t}$, the transition kernel $P^\rat_{\t}$, and the solution to a
Poisson equation associated with $P^\rat_{\t}$. We then give the proof of Proposition~\ref{prop:cvggal}.
\begin{alem}
  \label{lem:diffdespi}
  For all $\t, \t' \in \Theta$, 
\[
\| \pi^\rat_{\t} \, \rmd \lambda - \pi^\rat_{\t'}\, \rmd \lambda \|_{\tv} \leq
2 (d-1) \ \sum_{i=1}^d \left|1 - \frac{\rat(\t'(i))}{\rat(\t(i))} \right|.
\]
\end{alem}
Notice that the result of course also holds with the symmetrized
right-hand side $$2 (d-1) \min \left( 
\sum_{i=1}^d \left|1 - \frac{\rat(\t'(i))}{\rat(\t(i))} \right|,
\sum_{i=1}^d \left|1 - \frac{\rat(\t(i))}{\rat(\t'(i))} \right|
\right).$$
\begin{adem}
  The proof is adapted from~\cite[Lemma
  4.6]{fort:jourdain:kuhn:lelievre:stoltz:2014}. By definition of $\pi^\rat_\t$,
\[
\pi^\rat_\t(x) = \sum_{i=1}^d \frac{\t_\star(i)/ \rat(\t(i))}{\sum_{j=1}^d[\t_\star(j)/
  \rat(\t(j))]} \ \frac{\pi(x)}{\t_\star(i)} \ \un_{\Xset_i}(x) \eqsp.
\]
Hence,
\begin{align*}
  \|\pi^\rat_\t \, \rmd \lambda - \pi^\rat_{\t'} \, \rmd \lambda \|_\tv & \leq \sum_{i=1}^d
  \left| \frac{\t_\star(i)/ \rat(\t(i))}{\sum_{k=1}^d[\t_\star(k)/
      \rat(\t(k))]} - \frac{\t_\star(i)/
      \rat(\t'(i))}{\sum_{k=1}^d[\t_\star(k)/ \rat(\t'(k))]} \right| \\
  & \leq \frac{\sum_{j=1}^d \sum_{i=1}^d
    \t_\star(i)\t_\star(j)\left|1/[\rat(\t(i)) \rat(\t'(j))] - 1/[\rat(\t'(i))
      \rat( \t(j))]\right|}{\sum_{k=1}^d[\t_\star(k)/ \rat(\t(k))] \ \ 
    \sum_{l=1}^d[\t_\star(l)/ \rat(\t'(l))]} \eqsp.
\end{align*}
We denote by $N(\t, \t')$ the numerator of the expression of the right-hand side of the previous inequality. Then,
\begin{align*}
  N(\t, \t') &= \sum_{j=1}^d \sum_{i\neq j} \t_\star(i)\t_\star(j) \ 
  \frac{\left| \rat(\t'(i)) \, \rat( \t(j)) - \rat(\t(i)) \,  \rat(\t'(j)) \right|}{\rat(\t(i)) \, \rat(\t'(i)) \, \rat(\t(j)) \, \rat( \t'(j)) } \\
  &\leq \sum_{j=1}^d \sum_{i\neq j} \t_\star(i)\t_\star(j) \ \left|\frac{
      \rat(\t(j)) - \rat(\t'(j))}{\rat(\t(i)) \, \rat(\t(j)) \, \rat(\t'(j)) }
  \right|+ \sum_{j=1}^d \sum_{i\neq j} \t_\star(i)\t_\star(j) \ \left| \frac{
      \rat(\t(i)) - \rat(\t'(i)) }{\rat(\t(i)) \, \rat(\t'(i)) \, \rat(\t(j)) }
  \right| \eqsp.
\end{align*}
For the denominator, we use the lower bound
\[
\forall i,j \in \{ 1,\dots,d \}, \qquad
\sum_{k=1}^d[\t_\star(k)/ \rat(\t(k))] \ \ \sum_{l=1}^d[\t_\star(l)/ \rat(\t'(l))] \geq
\frac{\t_\star(i) \, \t_\star(j)}{\rat(\t(i)) \,  \rat(\t'(j))}.
\]
Therefore,
\[
\|\pi^\rat_\t \, \rmd \lambda - \pi^\rat_{\t'} \, \rmd \lambda \|_\tv \leq 2\sum_{j=1}^d \sum_{i\neq j} \left|\frac{\rat( \t(j))
    -\rat(\t'(j))}{ \rat(\t(j)) } \right| \leq 2 (d-1) \sum_{j=1}^d \frac{\left|\rat(\t(j))
    -\rat(\t'(j))\right|}{\rat(\t(j))} \eqsp,
\]
which gives the claimed result.
\end{adem}

\begin{alem} \label{lem:diffdesnoyaux}
  For all $\t, \t' \in \Theta$ and $x \in \Xset$
\begin{equation}\label{eq:diffdesnoyaux}
\| P^\rat_{\t}(x,\cdot) - P^\rat_{\t'}(x,\cdot) \|_{\tv} \leq 4 \   \sup_{i\in \{1,\ldots,d\}}
\left|1 - \frac{\rat(\t'(i))}{\rat(\t(i))}  \right|.
\end{equation}
\end{alem}
\begin{adem}
  The proof is adapted from~\cite[Lemma
  4.7]{fort:jourdain:kuhn:lelievre:stoltz:2014}.  As $P^\rat_\t$ is a Metropolis
  kernel, for any bounded measurable function $f$,
\begin{align*}
  \left| P^\rat_\t f(x) - P^\rat_{\t'} f(x) \right| & = \left| \int_\Xset q(x,y) \left(
      \alpha_\t(x,y) -
      \alpha_{\t'}(x,y) \right) \left(f(y) - f(x) \right) \rmd \lambda(y)  \right| \\
  & \leq 2 \sup_\Xset |f| \ \sup_{\Xset^2} \left|\alpha_\theta - \alpha_{\t'}
  \right| 
\end{align*}
where $\alpha_\t(x,y) \eqdef 1 \wedge (\pi^\rat_\t(y) / \pi^\rat_\t(x))$ since $q$ is
symmetric. For notational simplicity, we do not indicate explicitly
the dependence of $\alpha_\t$ on $\rat$.

 Let us introduce the unnormalized measure
$$\forall x \in \Xset,\, \tilde{\pi}^\rat_\t(x) \eqdef \sum_{i=1}^d
\frac{\pi(x)}{ \rat(\theta(i))}\un_{\Xset_i}(x),$$
which is such that $\forall x \in
\Xset, \pi^\rat_\t(x)=(Z_\t^\rho)^{-1} \tilde{\pi}^\rat_\t(x)$ where
$Z_\t^\rho=\int_{\Xset} \tilde{\pi}^\rat_\t \rmd \lambda=\sum_{j=1}^d\frac{\theta_\star(j)}{
  \rat(\theta(j))}$. Notice that
\[
\forall x,y \in \Xset, \, \alpha_\t(x,y) = 1 \wedge \frac{\tilde{\pi}^\rat_\t(y)}{\tilde{\pi}^\rat_\t(x)}.
\] 
We now show that 
\begin{equation}\label{eq:diff_alpha}
\forall \t, \t' \in \Theta, \, \forall x,y \in \Xset, \, |\alpha_\t(x,y)-\alpha_{\t'}(x,y)| \le  2
\sup_\Xset \left| 1 - \frac{\tilde\pi^\rat_{\t}}{\tilde\pi^\rat_{\t'}}\right|
\end{equation}
which yields the result~\eqref{eq:diffdesnoyaux} since 
$$\sup_{\Xset} \left| 1 - \frac{\tilde{\pi}^\rat_{\t}}{\tilde{\pi}^\rat_{\t'}}
  \right| = \sup_{i
  \in \{1, \ldots, d\} } \left| 1 - \frac{\rat(\t'(i))}{\rat(\t(i))}
  \right|.
$$
The proof of~\eqref{eq:diff_alpha} is performed by distinguishing between four cases:

$\bullet$ $\tilde{\pi}^\rat_\t(y) \leq \tilde{\pi}^\rat_\t(x)$ and $\tilde{\pi}^\rat_{\t'}(y) \leq \tilde{\pi}^\rat_{\t'}(x)$.
Then,
\begin{align*}
  \left| \alpha_\t(x,y) - \alpha_{\t'}(x,y) \right| &= \left|
    \frac{\tilde{\pi}^\rat_\t(y)}{\tilde{\pi}^\rat_\t(x)}-\frac{\tilde{\pi}^\rat_\t(y)}{\tilde{\pi}^\rat_{\t'}(x)}+\frac{\tilde{\pi}^\rat_\t(y)}{\tilde{\pi}^\rat_{\t'}(x)}- \frac{\tilde{\pi}^\rat_{\t'}(y)}{\tilde{\pi}^\rat_{\t'}(x)} \right|
  \\&\leq\frac{\tilde{\pi}^\rat_\t(y)}{\tilde{\pi}^\rat_{\t}(x)}\frac{
    \left| \tilde{\pi}^\rat_{\t'}(x) -\tilde{\pi}^\rat_\t(x)  \right|}{\tilde{\pi}^\rat_{\t'}(x)}  + \frac{ \left| \tilde{\pi}^\rat_\t(y) - \tilde{\pi}^\rat_{\t'}(y) \right|}{\tilde{\pi}^\rat_{\t'}(x)} \\
  & \leq \frac{
    \left| \tilde{\pi}^\rat_\t(x) -
      \tilde{\pi}^\rat_{\t'}(x) \right|}{\tilde{\pi}^\rat_{\t'}(x)}+\frac{ \left| \tilde{\pi}^\rat_\t(y) - \tilde{\pi}^\rat_{\t'}(y) \right|}{\tilde{\pi}^\rat_{\t'}(y)}  \\
  &\leq 2 \sup_\Xset \left| 1 - \frac{\tilde{\pi}^\rat_{\t}}{\tilde{\pi}^\rat_{\t'}}\right| \eqsp.
\end{align*}

$\bullet$ $\tilde{\pi}^\rat_\t(x) \leq \tilde{\pi}^\rat_\t(y)$ and $\tilde{\pi}^\rat_{\t'}(x) \leq \tilde{\pi}^\rat_{\t'}(y)$.
Then, $ \left| \alpha_\t(x,y) - \alpha_{\t'}(x,y) \right| =0$.

$\bullet$ $\tilde{\pi}^\rat_\t(x) \leq \tilde{\pi}^\rat_\t(y)$ and $\tilde{\pi}^\rat_{\t'}(y) \leq \tilde{\pi}^\rat_{\t'}(x)$.
Then, 
\begin{align*}
   \left| \alpha_\t(x,y) - \alpha_{\t'}(x,y) \right| &= 1 -
\frac{\tilde{\pi}^\rat_{\t'}(y)}{\tilde{\pi}^\rat_{\t'}(x)} \leq 1 -
\frac{\tilde{\pi}^\rat_{\t'}(y)}{\tilde{\pi}^\rat_{\t'}(x)}\frac{\tilde{\pi}^\rat_{\t}(x)}{\tilde{\pi}^\rat_{\t}(y)} 
\leq \frac{|\tilde{\pi}^\rat_{\t'}(x)-\tilde{\pi}^\rat_{\t}(x)|}{\tilde{\pi}^\rat_{\t'}(x)}+\frac{\tilde{\pi}^\rat_{\t}(x)}{\tilde{\pi}^\rat_{\t'}(x)}\frac{|\tilde{\pi}^\rat_{\t}(y)-\tilde{\pi}^\rat_{\t'}(y)|}{\tilde{\pi}^\rat_{\t}(y)}\\
& \leq  \frac{|\tilde{\pi}^\rat_{\t'}(x)-\tilde{\pi}^\rat_{\t}(x)|}{\tilde{\pi}^\rat_{\t'}(x)}+\frac{\tilde{\pi}^\rat_{\t}(x)}{\tilde{\pi}^\rat_{\t}(y)} \frac{\tilde{\pi}^\rat_{\t'}(y)}{\tilde{\pi}^\rat_{\t'}(x)}\frac{|\tilde{\pi}^\rat_{\t}(y)-\tilde{\pi}^\rat_{\t'}(y)|}{\tilde{\pi}^\rat_{\t'}(y)} \\
&\leq 2 \sup_{\Xset} \left| 1 - \frac{\tilde{\pi}^\rat_{\t}}{\tilde{\pi}^\rat_{\t'}}
  \right| \eqsp.
\end{align*}

$\bullet$ $\tilde{\pi}^\rat_\t(y) \leq \tilde{\pi}^\rat_\t(x)$ and $\tilde{\pi}^\rat_{\t'}(x) \leq \tilde{\pi}^\rat_{\t'}(y)$.
If $\tilde{\pi}^\rat_{\t}(x) \leq \tilde{\pi}^\rat_{\t'}(y)$, it holds
\[
\left| \alpha_\t(x,y) - \alpha_{\t'}(x,y) \right| = 1 -
\frac{\tilde{\pi}^\rat_\t(y)}{\tilde{\pi}^\rat_\t(x)} \leq 1 - \frac{\tilde{\pi}^\rat_\t(y)}{\tilde{\pi}^\rat_{\t'}(y)} \leq
\sup_\Xset \left|1 - \frac{\tilde{\pi}^\rat_\t}{\tilde{\pi}^\rat_{\t'}} \right|.
\]
Otherwise, we have $\tilde{\pi}^\rat_{\t'}(x) \leq \tilde{\pi}^\rat_{\t'}(y) \leq \tilde{\pi}^\rat_\t(x)$ and we write
\begin{align*}
  \left| \alpha_\t(x,y) - \alpha_{\t'}(x,y) \right| =  1 - \frac{\tilde{\pi}^\rat_\t(y)}{\tilde{\pi}^\rat_\t(x)} & = \frac{\tilde{\pi}^\rat_{\t'}(x)}{\tilde{\pi}^\rat_\t(x)} \left(\frac{\tilde{\pi}^\rat_\t(x)}{\tilde{\pi}^\rat_{\t'}(x)} - \frac{\tilde{\pi}^\rat_{\t'}(y)}{\tilde{\pi}^\rat_{\t'}(x)} \right) +  \frac{\tilde{\pi}^\rat_{\t'}(y)}{\tilde{\pi}^\rat_\t(x)}    \left(1 - \frac{\tilde{\pi}^\rat_\t(y)}{\tilde{\pi}^\rat_{\t'}(y)}\right)  \\
  & \leq  \left(\frac{\tilde{\pi}^\rat_\t(x)}{\tilde{\pi}^\rat_{\t'}(x)} - 1 \right) + \left|1 - \frac{\tilde{\pi}^\rat_\t(y)}{\tilde{\pi}^\rat_{\t'}(y)}\right| \\
  & \leq 2 \sup_\Xset \left| 1 - \frac{\tilde{\pi}^\rat_\t}{\tilde{\pi}^\rat_{\t'}} \right| \eqsp.
\end{align*}
This shows that~\eqref{eq:diff_alpha} holds, and thus concludes the proof.
\end{adem}

\begin{alem}
  \label{lem:RegularitePoisson}
  Assume A\ref{hyp:targetpi}, A\ref{hyp:kernel} and R\ref{hyp:BoundFctRat}.
  Then, $\sup_{\theta \in \Theta} \sup_{x \in \Xset} |H(x,\theta)| \leq
  \sqrt{2} \left( \sup_{(0,1)} \rat \right)$ where~$H$ is defined
  by~\eqref{eq:DefinitionH}.  In addition, for any $\t \in \Theta$, there
  exists a unique function $\hatH_\t$ solving the Poisson equation
\begin{equation}\label{eq:Poisson}
\hatH_\t - P^\rat_\t \hatH_\t = H(\cdot, \t) - h(\t), \qquad \pi_\t^\rat\left(\hatH_\t\right) = 0,
\end{equation}
where, we recall, $h=\pi^\rat_\t(H(\cdot,\t))$ is defined by~\eqref{eq:Definitionh}.
Moreover, $\hatH_\t$ is uniformly bounded: $$\sup_{\t \in \Theta, x \in \Xset}
  \left|\hatH_\t(x)\right| < \infty,$$
and there exists a positive constant $C$ such that, for any $\t,\t' \in
  \Theta$,
  \[
  \sup_{\Xset} \left\{ \left| \hatH_\t - \hatH_{\t'} \right| + \left| P^\rat_\t
      \hatH_\t - P^\rat_{\t'} \hatH_{\t'} \right| \right \} \leq C \left( |\t -
    \t'| + \sum_{i=1}^d \left|1 - \frac{\rat(\t'(i))}{\rat(\t(i))}\right|
  \right) \eqsp.
  \]
\end{alem}
Notice that for notational simplicity, we do not indicate explicitly
the dependence of~$\hatH_\t$ on~$\rat$.
\begin{adem}
   Using the Euclidean norm: for all
  $(x,\t) \in \Xset \times \Theta $,
\begin{align}
|H(x, \t)|^2&=\sum_{i=1}^d (H_i(x,\t))^2 \le \sum_{i=1}^d \rat^2(\theta(i))\un_{\Xset_i}(x)+  \t(i)^2 \rat^2
(\theta(I(x)))\nonumber\\
&\le \sup _{(0,1)}  \rat^2 + \sum_{i=1}^d  \t(i) \sup _{(0,1)}
\rat^2= 2 \sup _{(0,1)}  \rat^2,
\label{eq:majoH}
\end{align}
so that $\sup_{\t \in \Theta} \sup_{x \in \Xset} |H(x, \t)| \leq
\sqrt{2}\sup_{(0,1)}  \rat$. Set
\[
\hatH_\t(x) \eqdef \sum_{n \geq 0} \left(\int_{\Xset} (P^\rat_\theta)^n(x, \rmd y) H(y,\theta)
  - h(\theta) \right).
\]
Proposition~\ref{prop:unifergo} shows that $\hatH_\t$ exists
for any $\t \in \Theta$. It is easily seen that this function satisfies~\eqref{eq:Poisson}. 
Moreover, there exists a  constant $C$ such that  (see
  \textit{e.g.}~\cite[Section~17.4.1]{meyn:tweedie:2009})
\begin{equation}
  \label{eq:upperbound:HatH}
\sup_{\t \in \Theta} \sup_{x \in \Xset} \left|\hatH_\t(x)\right| 
\leq \sup_{\t \in \Theta} \sup_{x \in \Xset} \sum_{n \geq 0} \left| (P^\rat_\t)^n
H(\cdot,\t)(x) - \pi^\rat_\t(H(\cdot, \t)) \right| \leq C \ \sup_{(0,1)} \rat \eqsp. 
\end{equation}
Notice that $\hatH_\t(x)$ of course depends on the choice of the
function $\rho$, even if we do not indicate it explicitly for the
ease of notation. In view of~\cite[Lemma~4.2]{fort:moulines:priouret:2011}
(using the constant function equal to 1 as a Lyapunov function, thanks
to Proposition~\ref{prop:unifergo}), there exists a  constant $C$
such that, for any $\t, \t' \in \Theta$,
\begin{multline*}
  \sup_\Xset \left| P^\rat_{\t}\hatH_{\t} - P^\rat_{\t'}\hatH_{\t'} \right| +
  \sup_{\Xset} \left| \hatH_\t - \hatH_{\t'} \right|
\\ \leq C \left(
    \sup_\Xset \left| H(\cdot,\t) - H(\cdot,{\t'}) \right| + \sup_{x \in \Xset}
    \| P^\rat_\t(x,\cdot) - P^\rat_{\t'}(x,\cdot) \|_\tv + \|\pi^\rat_\t \, \rmd \lambda - \pi^\rat_{\t'} \, \rmd\lambda\|_\tv
  \right) \eqsp.
\end{multline*}
By the definition \eqref{eq:DefinitionH} of $H$, $$\sup_\Xset \left| H_i(\cdot,\t) -
    H_i(\cdot,{\t'}) \right| \leq \sup_{j \in \{1, \ldots ,d\}} |\rat(\t(j)) - \rat(\t'(j))| + |\t(i) -
  \t'(i)| \, \sup_{(0,1)}\rat.$$ Therefore, by R\ref{hyp:BoundFctRat},  there exists a constant $C'$ such that for any $\t, \t' \in \Theta$,
\[
\sup_\Xset \left| H(\cdot,\t) - H(\cdot,{\t'}) \right| \leq C' \left\{  \sum_{i=1}^d \left| 1 - \frac{\rat(\t'(i))}{\rat(\t(i))} \right| + | \t -
  \t'| \right\}  \eqsp.
\]
The proof is then concluded by Lemmas~\ref{lem:diffdespi} and
\ref{lem:diffdesnoyaux}.
\end{adem}

We are now in position to prove Proposition~\ref{prop:cvggal}, by
considering successively the three items in Theorem~\ref{theo:cvggal}.

\textit{(i)} The proof of the first item consists in verifying the
sufficient conditions given in
\cite[Theorems~2.2 and~2.3]{andrieu:moulines:priouret:2005} for the
convergence of stochastic approximation algorithms.
Remember the definition~\eqref{eq:Definitionh} of the mean field function $h: \Theta \to \R^d$:
\[
h(\tn) = \int_{\Xset} H(x,\tn) \, \pi^\rat_\tn( \rmd x)
=\frac{\tstar-\theta}{\sum_{i=1}^d\frac{\tstar(i)}{\rat(\theta(i))}}\eqsp.
\]
By R\ref{hyp:LypFctRat}, the function $h$ is continuous on $\Theta$.
By~\cite[Proposition~4.5]{fort:jourdain:kuhn:lelievre:stoltz:2014},
the
function $U$ defined on $\Theta$ by
\[
U(\tn) \eqdef - \sum_{i=1}^d \tstar(i) \ \ln\left(
  \frac{\tn(i)}{\tstar(i)}\right)
\]
is non negative (thanks to Jensen's inequality),
continuously differentiable on $\Theta$ and the level sets $\left(\{\tn \in \Theta:
U(\tn) \leq M \}\right)_{M > 0}$ are a family of closed compact
neighborhood of $\tstar$ in the open set~$\Theta$.
We also have $\left<\nabla U(\tn),h(\tn)\right> \leq 0$ and $\left<\nabla
  U(\tn),h(\tn)\right> = 0$ if and only if $\tn= \tstar$.  Hence, the
assumption A1 of \cite{andrieu:moulines:priouret:2005} is satisfied with
$\mathcal{L} = \{\tstar \}$.

Moreover, under our assumptions, the conditions on the stepsize sequence $(\gamma_n)_{n
  \geq 1}$ in~\cite[Theorems~2.2 and~2.3]{andrieu:moulines:priouret:2005} hold
almost-surely.
To apply these two theorems which respectively show the stability and
the convergence of the algorithm, it is thus enough to prove that for
any compact subset ${\calK}$ of $\Theta$,
\begin{equation}\label{lastcond}
  \P-a.s. \qquad   \lim_k \sup_{ \ell \geq k} \left| \sum_{n=k}^{\ell} \gamma_{n+1} \Big(
  H(X_{n+1},\tn_n) - h(\tn_n) + \Lambda_{n+1} \Big) \un_{\t_n \in \calK} \right| =0.
\end{equation}
Indeed, a slight adaptation
of~\cite[Theorem~2.2]{andrieu:moulines:priouret:2005} shows that $\P$-a.s. the
sequence $(\tn_n)_{n \ge 0}$ remains in a compact subset of $\Theta$ under the
conditions A1 of \cite{andrieu:moulines:priouret:2005} together with $\lim_n \gamma_n=0$
$\P$-a.s.,  the recurrence property~\eqref{condstab} and~\eqref{lastcond}. Then,
\cite[Theorem~2.3]{andrieu:moulines:priouret:2005} ensures the a.s convergence
of $(\tn_n)_{n \ge 0}$ to $\tstar$ under the additional assumption $\sum_n \gamma_n =
\infty$, $\P$-a.s.

Let us now check \eqref{lastcond}.  By Lemma~\ref{lem:SHUSalpharho:AS},
\[
\P\left( \forall k, \sup_{\ell \geq k} \left|\sum_{n=k}^\ell \gamma_{n+1}
    \Lambda_{n+1} \right|\leq \ \left( \sup_{(0,1)} \rat \right)^2 \, \sqrt{2}
  \ \sum_{n \geq k}\gamma_{n+1}^2 \right) =1 \eqsp
\]
so that R\ref{hyp:BoundFctRat} and \eqref{hyp:Gamma} imply that $\sup_{\ell
  \geq k} \left|\sum_{n=k}^\ell \gamma_{n+1} \Lambda_{n+1} \right|$ converges
to $0$ a.s. as $k\to\infty$.

To deal with $H(X_{n+1},\tn_n) - h(\tn_n)$, for each $\theta\in\Theta$, we introduce the Poisson equation 
\[
\forall x\in\Xset, \quad g(x) - P^\rat_\tn g(x) = H(x,\tn) - h(\tn), \qquad \pi_\t^\rat(g) = 0,
\]
whose unknown is the function $g:\Xset\to\R$.  By
Lemma~\ref{lem:RegularitePoisson}, this equation admits a unique solution
$\hatH_\tn(x)$, which is moreover uniformly bounded in $(\theta,x)$.
We write \begin{align*} H(X_{n+1},\tn_n) - h(\tn_n) &=
  \hatH_{\tn_n}(X_{n+1})-P^\rat_{\tn_n}\hatH_{\tn_n}(X_{n+1})=\mathcal{E}_{n+1} +
  R_{n+1}^{(1)} +R_{n+1}^{(2)},
\end{align*} with
\begin{align*}
  \mathcal{E}_{n+1} & =  \hatH_{\tn_n}(X_{n+1}) - P^\rat_{\tn_n}\hatH_{\tn_n}(X_{n})  \eqsp, \\
  R_{n+1}^{(1)} & = P^\rat_{\tn_n}\hatH_{\tn_n}(X_{n}) - P^\rat_{\tn_{n+1}}\hatH_{\tn_{n+1}}(X_{n+1}) \eqsp, \\
  R_{n+1}^{(2)} & = P^\rat_{\tn_{n+1}}\hatH_{\tn_{n+1}}(X_{n+1}) -
  P^\rat_{\tn_{n}}\hatH_{\tn_{n}}(X_{n+1}) \eqsp.
\end{align*}
Recall that $\gamma_{n+1}$ is $\mathcal{F}_n$-measurable. Let us first check that the
martingale $(M_k)_{k \geq 1}$ defined by $M_k \eqdef
\sum_{n=1}^k\gamma_n\mathcal{E}_n$ converges a.s. as $k\to\infty$, which will
imply that a.s.
\begin{equation}\label{eq:rest}
\lim_k \sup_{\ell \geq k} \left|\sum_{n=k}^\ell \gamma_{n+1} \mathcal{E}_{n+1}
\right|= 0 \ .
\end{equation}
Let us first prove the result assuming that $\gamma_1$ is square integrable,
which implies that $M_k$ is also square integrable. Indeed, for all $k \ge 1$,
$|M_k|\leq 2k \gamma_1 \sup_{\tn \in \Theta} \sup_{x \in \Xset}
\left|\hatH_\tn(x)\right|$.  The latter inequality holds since $(\mathcal{E}_{n})_{n
  \geq 0}$ is bounded by $2\sup_{\tn \in \Theta} \sup_{x \in \Xset}
|\hatH_\tn(x)|$ and, by \eqref{hyp:Gamma}, $(\gamma_n)_{n \geq 1}$ is bounded
by $\gamma_1$. Moreover, $\gamma_{n+1}$ is ${\cal F}_n$-measurable and the
conditional distribution of $X_{n+1}$ given ${\cal F}_n$ is
$P^\rat_{\tn_n}(X_n,\cdot)$, so that
\[
  \E(\gamma_{n+1}\mathcal{E}_{n+1} \,| \, {\cal F}_n) =\gamma_{n+1}\left[
    \E\left( \left. \hatH_{\tn_n}(X_{n+1}) \right|{\cal F}_n\right) -
    P^\rat_{\tn_n}\hatH_{\tn_n}(X_{n}) \right]=0 \eqsp.
\]
In conclusion, $(M_k)_{k \geq 1}$ is a square integrable ${\cal
  F}_k$-martingale.  Since
\[
\sum_{n}\E\left[\left. (M_{n+1}-M_n)^2 \right|{\cal F}_n\right]=\sum_{n
}\gamma_{n+1}^2\E(\mathcal{E}_{n+1}^2|{\cal F}_n)
\]
is smaller than $C\sum_{n\geq 1}\gamma_n^2$ which is a.s. finite by
\eqref{hyp:Gamma}, $(M_k)_{k \geq 1}$ converges a.s. by \cite[Theorem
2.15]{hall:heyde:1980} and this implies~\eqref{eq:rest}. Now, if $\gamma_1$ is
not square integrable, one can apply the above argument upon replacing $\gamma_1$ by
$\gamma_1 \wedge \Gamma$ where $\Gamma \in \N$ is a constant.  This shows
that~\eqref{eq:rest} holds almost surely on the event $\{\gamma_1<\Gamma\}$,
and thus on the event $\cup_{\Gamma=1}^\infty \{\gamma_1<\Gamma\}$. Since the set
$\cup_{\Gamma=1}^\infty \{\gamma_1<\Gamma\}=\{\gamma_1 < \infty\}$ is of
probability one,~\eqref{eq:rest} holds almost surely.

We now consider the term $R_{n+1}^{(1)}$. By the monotonic property of
$(\gamma_n)_{n \geq 1}$ and since $\hatH_\t$ is uniformly
bounded in $(\theta,x)$, following the same lines
as in the proof of \cite[Proposition~4.10]{fort:jourdain:kuhn:lelievre:stoltz:2014}, it can be checked that there
exists a constant $C$ such that
\begin{equation}\label{eq:rest1}
\P\left( \forall k,   \sup_{\ell \geq k} \left|\sum_{n=k}^\ell \gamma_{n+1} R_{n+1}^{(1)}
  \right|\leq C \ \gamma_{k+1} \right) =1 \eqsp. 
\end{equation}
The argument is based on a summation by parts and the fact that the
series $\sum_n R_{n}^{(1)}$ is telescoping. From~\eqref{eq:rest1} and~\eqref{hyp:Gamma}, $\sup_{\ell \geq k} \left|\sum_{n=k}^\ell \gamma_{n+1} R_{n+1}^{(1)}
\right|$ tends to zero a.s. as $k \to \infty$.

We now consider the term $R_{n+1}^{(2)}$.  By
Lemma~\ref{lem:RegularitePoisson}, there exists a constant $C$ such that for
any $\tn, \tn' \in \Theta$
\begin{equation}\label{eq:majoTerme2}
\sup_\Xset \left| P^\rat_{\tn}\hatH_{\tn} - P^\rat_{\tn'}\hatH_{\tn'} \right| \leq C
\left( |\tn - \tn'| + \sum_{i=1}^d \left|1 -
    \frac{\rat(\tn'(i))}{\rat(\tn(i))} \right| \right)\eqsp.
\end{equation}
By R\ref{hyp:LypFctRat}, for any compact subset $\calK$ of $\Theta$, there exists a constant $C$ such that for any $ n\geq 0$,
\begin{equation}\label{eq:controle:fluctuation:rat}
\sum_{i=1}^d \left| 1 - \frac{\rat(\t_{n+1}(i))}{\rat(\t_n(i))} \right|
\un_{\t_n \in \calK} \leq C \, \left| \t_{n+1} -\t_n\right|.
\end{equation}
Moreover, by Lemma~\ref{lem:SHUSalpharho:AS}, R\ref{hyp:BoundFctRat} and the
boundedness of $H$ (see Lemma~\ref{lem:RegularitePoisson}) there exists a
constant $C$ such that with probability one, for any $n \geq 0$,
\begin{equation}\label{eq:controle:fluctuation:rat1}
\left| \t_{n+1} -\t_n\right| \leq C \, \gamma_{n+1}\eqsp.
\end{equation}
Therefore,
combining~\eqref{eq:majoTerme2}--\eqref{eq:controle:fluctuation:rat}--\eqref{eq:controle:fluctuation:rat1},
there exists a constant $C$ such that
\[
\P\left( \forall k, \, \sup_{\ell \geq k} \left|\sum_{n=k}^\ell \gamma_{n+1}
    R_{n+1}^{(2)} \un_{\t_n \in \calK} \right|  \leq C \, \sum_{n \geq k} \gamma_{n+1}^{2} \right)
=1 \eqsp.
\]
By \eqref{hyp:Gamma},  $\sup_{\ell \geq k} \left|\sum_{n=k}^\ell \gamma_{n+1}
  R_{n+1}^{(2)} \un_{\t_n \in \calK} \right|$ tends to zero a.s. as $k \to \infty$.
This concludes the proof of the a.s. convergence: $\lim_{n \to \infty} \t_n=\t_\star$.

\textit{(ii)} The proof follows the same lines as the proof of \cite[Theorem
3.4]{fort:jourdain:kuhn:lelievre:stoltz:2014} and details are omitted.  The
only result which has to be adapted is \cite[Corollary
4.8]{fort:jourdain:kuhn:lelievre:stoltz:2014}. Combining
Lemmas~\ref{lem:SHUSalpharho:AS}, ~\ref{lem:diffdespi}, \ref{lem:diffdesnoyaux}
and the
estimates~\eqref{eq:controle:fluctuation:rat}--\eqref{eq:controle:fluctuation:rat1},
we easily obtain the existence of a constant $C$ such that almost
surely, for any $n \geq 1$ on the set $\{\t_n \in \calK \}$,
  \begin{align*}
    & \| \pi^\rat_{\tn_{n+1}} \, \rmd \lambda - \pi^\rat_{\tn_{n}} \, \rmd \lambda
    \|_{\mathrm{TV}} + \sup_{x \in \Xset} \|P^\rat_{\tn_n}(x, \cdot) -
    P^\rat_{\tn_{n+1}}(x,\cdot) \|_{\mathrm{TV}} \leq C \gamma_{n+1} \eqsp.
   \end{align*}

\textit{(iii)} The proof is very similar to the proof of 
\cite[Theorem~3.5]{fort:jourdain:kuhn:lelievre:stoltz:2014} and is therefore omitted.

\subsection{Proof of Proposition~\ref{prop:recurrence}: recurrence
  of the algorithm}\label{sec:rec}

In all this section, we consider that the sequence is generated by the
WL$_\rat$ algorithm~\ref{algo:WL} (see Section~\ref{sec:relation_FE}).

The aim of this section is to give some sufficient conditions on $\left(\gamma_{n}\right)_{n \ge 1}$ such that $\P$-a.s., the sequence $(\t_n)_{n \ge
  0}$ visits a.s. infinitely often a compact subset of $\Theta$. For $n\ge 0$,
we set
\begin{equation*}
\underline{\tn}_n \eqdef \min_{1\le i\le d}\tn_n(i).
\end{equation*}
The objective is thus to verify
that a.s.  the sequence $(\underline{\tn}_n)_{n \ge 0}$ takes infinitely often
values in a compact subset of $(0,1)$. We will show this property along a sequence of well
chosen stopping times $(T_k)_{k\geq 0}$ defined inductively as follows.

We set $T_0=0$ and for $k\in\N$, $T_{k+1}=\infty$ if $T_k=\infty$. Otherwise
when $T_k<+\infty$, let for $m\in\N$, $\tn_{T_k+md}((1)_m)\leq
\tn_{T_k+md}((2)_m)\leq\hdots\leq \tn_{T_k+md}((d)_m)$ denote the increasing
reordering of $(\tn_{T_k+md}(i))_{1\leq i\leq d}$ and $$i_m \eqdef \max\{i\leq
d:\tn_{T_k+md}((i)_m)<R\underline{\tn}_{T_k+md}\}\mbox{ where $R$ is given by
  R\ref{hyp:fdec}}.$$
We then introduce an event corresponding to visiting
successively the strata of small weights with indices $(i)_m$ for
$i\in\{1,\hdots,i_m\}$, in decreasing order:
\begin{equation}
  \label{defam}
   A_m \eqdef  \bigg\{ X_{T_k+md+1}\in\Xset_{(i_m)_m},X_{T_k+md+2}\in\Xset_{(i_m-1)_m},
   \hdots,X_{T_k+md+i_m}\in\Xset_{(1)_m} \bigg\}.
\end{equation}
The next stopping $T_{k+1}$ is then defined by 
$$T_{k+1}(\omega)=T_k(\omega)+d\times\inf\{m\geq 1:\omega\in A_{m-1}\}\mbox{
  with convention }\inf\emptyset=+\infty.$$
Note that $T_k \geq kd$ by definition. 
Let us first show some additional properties on this sequence of stopping times.

\begin{alem}\label{lemgeo}
  Assume A\ref{hyp:targetpi}, A\ref{hyp:kernel} and R\ref{hyp:fdec}. Then,
  $\P(\forall k\in\N,\;T_k<+\infty)=1$ and
  $$\exists p\in(0,1),\;\forall k,m\in\N,\,\P\left(T_{k+1}-T_k > md|{\cal
      F}_{T_k}\right)\leq (1-p)^m.$$
  In addition, \begin{equation}
    \P\left(\exists C_\star<+\infty,\;\forall k\in\N,\;T_k\leq C_\star
      k\right)=1.\label{majotk}
\end{equation}
\end{alem}

\begin{adem}
  The first two statements are a consequence of
\begin{equation}\label{eq:minop}
  \exists p\in(0,1),\;\forall k,m\in\N\mbox{ with
  }T_k<\infty,\;\P(A_m|{\cal F}_{T_k+md})\geq p.
\end{equation}
  This inequality is proved as follows (see the proof of \cite[Lemma
  4.2]{fort:jourdain:kuhn:lelievre:stoltz:2014} for a similar reasoning).
The main ingredient in the proof is the
  following inequality:
  \begin{align}
    \forall x\in\Xset,\;\forall
    i\in\{1,\hdots,d\},\;P^\rat_{\theta}(x,\Xset_{i})&=
    \int_{\Xset_{i}}q(x,y)\left(1\wedge
      \frac{\rat(\theta(I(x))) \, \pi(y)}{\rat(\theta(i)) \,
        \pi(x)}\right) \rmd \lambda(y)\notag\\
&\geq \left( \inf_{\Xset^2} q \right)
    \int_{\Xset_{i}}\left(1\wedge
      \frac{\rat(\theta(I(x))) \, \pi(y)}{\rat(\theta(i)) \,
        \sup_{\Xset}\pi}\right) \rmd \lambda(y)\notag\\
&\geq \left( \inf_{\Xset^2} q \right)
    \int_{\Xset_{i}}\left(\frac{\pi(y)}{\sup_{\Xset} \pi}\wedge
      \frac{\rat(\theta(I(x))) \, \pi(y)}{\rat(\theta(i)) \,
        \sup_{\Xset}\pi}\right)\rmd \lambda(y)\notag\\
&= c \,
    \theta_\star(i)\left(\frac{\rat(\theta(I(x)))}{\rat(\theta(i))}\wedge
      1\right), \label{minotrans}
\end{align}
where $c=\frac{\inf_{\Xset^2} q}{\sup_{\Xset}\pi}>0$ by A\ref{hyp:targetpi} and
A\ref{hyp:kernel}.  Now, for $j\in\{1,\hdots,i_m-1\}$, it holds on the event
$\{X_{T_k+md+1}\in\Xset_{(i_m)_m},\hdots,X_{T_k+md+j}\in\Xset_{(i_m+1-j)_m}\}$,
\[
\begin{aligned}
  & \frac{\theta_{T_k+md+j}((i_m+1-j)_m)}{\theta_{T_k+md+j}((i_m-j)_m)} \\
  &
  =\frac{\theta_{T_k+md}((i_m+1-j)_m)}{\theta_{T_k+md}((i_m-j)_m)}\times\left(
    1 + \gamma_{T_k+md+j} \frac{\rat(\theta_{T_k+md+j-1}((i_m+1-j)_m))}{\theta_{T_k+md+j-1}((i_m+1-j)_m)}\right) \eqsp.
\end{aligned}
\]
Both factors on the right-hand side are larger than~1 (the first one by definition of the ordered indices $(i)_m$), so that, by \eqref{minotrans} and the monotonicity of $\rho$,
\[
P^\rat_{\theta_{T_k+md+j}}\left(X_{T_k+md+j},\Xset_{(i_m-j)_m}\right) \geq
c\,\theta_\star((i_m-j)_m) \geq c \, \underline{\t}_\star \eqsp,
\] 
where
$\underline{\t}_\star=\min_{1 \le i \le d} \t_\star(i)$. Note that this implies in particular that $c \, \underline{\t}_\star \leq 1$. Using successively the strong Markov property
of the chain $(X_n,\theta_n)_{n \ge 0}$, a backward induction on $n$, the definition of
$i_m$, together with~\eqref{minotrans}, we have
\begin{align*}
  &\P(A_m|{\mathcal F}_{T_k+md})\\
  &=\E\left(\un_{\{X_{T_k+md+1}\in\Xset_{(i_m)_m},\hdots,X_{T_k+md+i_m-1}\in\Xset_{(2)_m}\}}P^\rat_{\theta_{T_k+md+i_m-1}}(X_{T_k+md+i_m-1},\Xset_{(1)_m})|{\mathcal F}_{T_k+md}\right)\\
  &\geq c \, \underline{\t}_\star\P\left(\left. \un_{\{X_{T_k+md+1}\in\Xset_{(i_m)_m},\hdots,X_{T_k+md+i_m-2}\in\Xset_{(3)_m}\}}P^\rat_{\theta_{T_k+md+i_m-2}}(X_{T_k+md+i_m-2},\Xset_{(2)_m})\right|{\mathcal F}_{T_k+md}\right)\\
  &\geq \left(c \, \underline{\t}_\star\right)^{i_m-1}P^\rat_{\theta_{T_k +md}}(X_{T_k+md},\Xset_{(i_m)_m})\\
  &\geq \left(c \, \underline{\t}_\star\right)^{i_m-1}c \,
  \underline{\t}_\star
  \frac{\rat(\theta_{T_k+md}(I(X_{T_k+md})))}{\rat(\theta_{T_k+md}((i_m)_m))} \\
&\geq\left(c \, \underline{\t}_\star\right)^{d}
\frac{\rat(\underline{\theta}_{T_k+md})}{\rat(\theta_{T_k+md}((i_m)_m)} \geq\left(c \, \underline{\t}_\star\right)^{d}
\frac{\rat(\theta_{T_k+md}((i_m)_m)/R)}{\rat(\theta_{T_k+md}((i_m)_m)}\ge \left(c \,
  \underline{\t}_\star\right)^{d} \inf_{t \in (0,1/R)} \frac{\rat(t)}{\rat(Rt)},
\end{align*}
where, for the last but one inequality, we used the monotonicity of $\rho$ and the fact that $c \, \underline{\t}_\star \in (0,1]$.
The proof is therefore concluded by setting $$p = \left(c \,
  \underline{\t}_\star\right)^{d} \inf_{t \in (0,1/R)} \frac{\rat(t)}{\rat(Rt)}$$ where $\inf_{t \in (0,1/R)} \frac{\rat(t)}{\rat(Rt)} >0$ by R\ref{hyp:fdec}.

This concludes the proof of~\eqref{eq:minop} and thus of the first two
statements of Lemma~\ref{lemgeo}.
The third statement can be
deduced from the second one by a coupling argument, as in the proof of \cite[Proposition
3]{fort:jourdain:lelievre:stoltz:2015}. Indeed, it can be shown that there exists
two sequences $(\tilde{T}_k)_{k \ge 0}$ and $(\tau_k)_{k \ge 1}$ such
that: (i) $(\tilde{T}_k)_{k \ge 0}$ has the same law as $(T_k)_{k \ge
  0}$, (ii) $(\tau_k)_{k \ge 1}$ are independent geometric random
variables with parameter $p$ and (iii) $\forall k \in \N, \,
\tilde{T}_{k+1} - \tilde{T}_k \le d \tau_{k+1}$. As a consequence,
\begin{align*}
  \P\left( \limsup_{k\to\infty}\frac{T_k}{k}\leq \frac{d}{p}\right) =\P\left(
    \limsup_{k\to\infty}\frac{\tilde{T}_k}{k}\leq \frac{d}{p}\right) 
   \ge \P\left(\limsup_{k\to\infty}\frac{1}{k}\sum_{j=1}^k \tau_j\leq
    \frac{1}{p}\right)=1 \eqsp,
\end{align*}
the last equality being a consequence of
the strong law of large numbers. This concludes the proof of~\eqref{majotk}.
\end{adem}

\begin{arem} We proved Lemma~\ref{lemgeo} for a Metropolis-Hastings kernel,
  but it actually holds in a more general setting. Indeed,
   Assume $\min_{1\le i\le d}\tn_\star(i)>0$, \eqref{eq:minorization_bound} and R\ref{hyp:fdec}. Then the conclusion of Lemma~\ref{lemgeo} still holds.

The proof of this result is the following. By
\eqref{eq:minorization_bound} and the monotonicity of $\rat$, it holds:
$\forall \tn\in\Theta,\;\forall i\in\{1,\hdots,d\},\;\forall x\in\Xset,$
\begin{align*}
P^\rat_\tn(x,\Xset_i)\ge \delta \frac{\tn_\star(i)}{\rat(\tn(i))}\left(\sum_{j=1}^d\frac{\tn_\star(j)}{\rat(\tn(j))}\right)^{-1}\ge \delta \tn_\star(i)\frac{\rat(\min_{1\le j\le d}\tn(j))}{\rat(\tn(i))}.
\end{align*}
With the definition of $i_m$ and the monotonicity of $\rat$, one deduces that
\begin{align}
   P^\rat_{\t_{T_k+md}}(X_{T_k+md},\Xset_{(i_m)_m})&\ge\delta
   \tn_\star((i_m)_m)\frac{\rat(\underline{\t}_{T_k+md})}{\rat(\t_{T_k+md}((i_m)_m)}\notag
   \\
&\ge \delta
\tn_\star((i_m)_m)\frac{\rat(\t_{T_k+md}((i_m)_m)/R)}{\rat(\t_{T_k+md}((i_m)_m))}\notag\\
&\ge \delta\underline{\tn}_\star\inf_{t \in (0,1/R)} \frac{\rat(t)}{\rat(Rt)}.\label{eq:mintrantkmd}
\end{align}
  For $j\in\{1,\hdots,i_m-1\}$, it holds on the event
$\{X_{T_k+md+1}\in\Xset_{(i_m)_m},\hdots,X_{T_k+md+j}\in\Xset_{(i_m+1-j)_m}\}$  that $\underline{\t}_{T_k+md+j}=\t_{T_k+md+j}((1)_m)$, $$\frac{\underline{\t}_{T_k+md+j}}{\theta_{T_k+md+j}((i_m-j)_m)}=\frac{\t_{T_k+md+j}((1)_m)}{\theta_{T_k+md+j}((i_m-j)_m)}=\frac{\t_{T_k+md}((1)_m)}{\theta_{T_k+md}((i_m-j)_m)}\ge \frac{1}{R}$$
so that, following the derivation of \eqref{eq:mintrantkmd}, $P^\rat_{\t_{T_k+md+j}}(X_{T_k+md+j},\Xset_{(i_m-j)_m})\ge \delta\underline{\tn}_\star\inf_{t \in (0,1/R)} \frac{\rat(t)}{\rat(Rt)}$.
Therefore, the conclusions of Lemma~\ref{lemgeo} hold with $p=\left(\delta\underline{\tn}_\star\inf_{t \in (0,1/R)} \frac{\rat(t)}{\rat(Rt)}\right)^d$.
\end{arem}
We are now in position to state the main result of this section.

\begin{alem} \label{lem:LimSupPositive}
  Assume A\ref{hyp:targetpi}, A\ref{hyp:kernel}, R\ref{hyp:BoundFctRat} to
  R\ref{hyp:ffininf} and that the sequence $(\gamma_n)_{n\geq 1}$ is
  non-increasing, bounded from above by a deterministic sequence converging to
  $0$ as $n\to\infty$ and such that $\bar{r}_{d,\gamma}<\infty$,
  where, we recall (see~\eqref{eq:hypr})
  \begin{equation*}
\bar{r}_{d,\gamma} = \sup_{n\geq
    1}\frac{\gamma_{n}}{\gamma_{n+d-1}}.
\end{equation*} 
Then
  $$\P\left(\limsup_{k\to\infty}\underline{\tn}_{T_k-d}>0\right)=1.$$
\end{alem}
Notice that this lemma implies that, almost surely, the sequence
$(\tn_n)_{n \ge 1}$ returns infinitely often to a compact subset of $\Theta$
(namely~\eqref{condstab}) since $\lim_{k \to \infty} T_k=\infty$  and by Lemma~\ref{lemgeo}, $\forall k\ge 0$,
$T_k<\infty$ almost surely. Therefore,  Proposition~\ref{prop:recurrence} is an immediate
consequence of Lemma~\ref{lem:LimSupPositive}.

 \begin{adem}
   The argument follows the proof of the second statement in \cite[Proposition
   4.1]{fort:jourdain:kuhn:lelievre:stoltz:2014}. For $k\geq 1$, we set
   $Y_k\eqdef \underline{\tn}_{T_k-d}.$ As a preliminary result, let us first
   prove that there exists $\underline{k}\in\N \setminus \{0\}$ and $\bar{y}\in (0,1)$ such that
\begin{equation}
  \forall k\geq\underline{k}, \quad Y_k\leq\bar{y} \Longrightarrow
\E(\ln(Y_{k+1})|{\mathcal F}_{T_k})\geq \ln(Y_k).\label{surm}
\end{equation}
One has
$$\forall i\in\{1,\hdots,d\},\;\forall
n\in\N,\;\tn_{n+1}(i)=\tn_n(i)\frac{1+\gamma_{n+1}
  \un_{\Xset_i}(X_{n+1}) \rat(\tn_n(i))/ \tn_n(i) \,
  }{1+\gamma_{n+1} \rat(\tn_n(I(X_{n+1}))}.$$
One deduces
that, on the one hand, for any index $i\in\{1,\hdots,d\}$ such that $\theta_{T_k-d}(i)\geq
R \underline{\theta}_{T_k-d}$, one has
\begin{align}
  \theta_{T_{k+1}-d}(i)&\geq
  \theta_{T_{k}-d}(i)\prod_{j=T_k-d+1}^{T_{k+1}-d}\frac{1}{1+\gamma_{j} \, \rat(\theta_{j-1}(I(X_{j})))}\geq
  \theta_{T_{k}-d}(i)\left(\frac{1}{1+\gamma_{T_k-d+1}
      \sup_{(0,1)}\rat}\right)^{T_{k+1}-T_k}\nonumber\\
&\geq
  R \left(\frac{1}{1+\gamma_{T_k-d+1} \sup_{(0,1)}\rat}\right)^{T_{k+1}-T_k}\underline{\theta}_{T_k-d},\label{eq:case1}
\end{align}
where we used the monotonicity of the sequence $(\gamma_n)_{n\geq 1}$ for the
second inequality. On the other hand, by definition of $T_k$, any stratum with index
$i\in\{1,\hdots,d\}$ such that $\theta_{T_k-d}(i)< R \underline{\theta}_{T_k-d}$
is visited at least once between the times $T_{k}-d+1$ and $T_k$ so that, using that
$\tn_n(i)$ decreases for $n$ between $T_k-d$ and this visit, as well as the monotonicities of
$t \mapsto \rat(t) /t$ and $n \mapsto \gamma_n$,
\begin{align}
  \theta_{T_{k+1}-d}(i)&\geq \theta_{T_{k}-d}(i) \left(
    1+\gamma_{T_k}\frac{\rat(R \underline{\theta}_{T_k-d})}{R \underline{\theta}_{T_k-d}}
  \right) \ \left(\frac{1}{1+\gamma_{T_k-d+1} \,
      \sup_{(0,1)}\rat}\right)^{T_{k+1}-T_k}\nonumber \\&\geq
\left(
    1+\gamma_{T_k}\frac{\rat(R \underline{\theta}_{T_k-d})}{R \underline{\theta}_{T_k-d}}
  \right) \ \left(\frac{1}{1+\gamma_{T_k-d+1} \,
      \sup_{(0,1)}\rat}\right)^{T_{k+1}-T_k}   \ \underline{\theta}_{T_k-d}.\label{eq:case2}
\end{align}
Combining~\eqref{eq:case1} and~\eqref{eq:case2}, one deduces that
\begin{align*}
   Y_{k+1}\geq \left(R \wedge \left(1+\gamma_{T_k}\frac{\rat(R Y_k)}{R  Y_k}\right)\right)\left(\frac{1}{1+ \bar{r}_{d,\gamma}\gamma_{T_k} \sup_{(0,1)} \rat }\right)^{T_{k+1}-T_k}Y_k.
\end{align*}
Taking the logarithm and remarking that the second statement in Lemma
\ref{lemgeo} implies that $\E(T_{k+1}-T_k|{\cal F}_{T_k})\leq \frac{d}{p}$,
one obtains that
\begin{align*}
  \E(\ln(Y_{k+1})|{\cal F}_{T_k})-\ln(Y_k)\geq \left(\ln (R )\wedge\ln
    \left(1+\gamma_{T_k}\frac{\rat(R Y_k)}{R 
      Y_k}\right)\right)-\frac{d}{p}\ln\left(1+\bar{r}_{d,\gamma}\gamma_{T_k}
\sup_{(0,1)} \rat \right).
\end{align*}
Since $T_k\geq kd$ and the sequence $(\gamma_n)_{n\geq 1}$ is bounded from
above by a deterministic sequence converging to $0$ one can choose
$\underline{k}\in\N \setminus \{0\}$ such that $\forall k\geq \underline{k}$,
$\gamma_{T_k}\leq
\bar{\gamma} = \frac{R ^{\frac{p}{d}}-1}{\bar{r}_{d,\gamma} \sup_{(0,1)}\rat}$ so that
$\ln(R )-\frac{d}{p}\ln(1+\bar{r}_{d,\gamma}\gamma_{T_k} \sup_{(0,1)}\rat)\geq 0$. Last,
since $\lim_{t\to 0^+}\rat(t)/t=+\infty$, one may choose $\bar{y}\in
(0,1/R )$ such that 
\[
\inf_{t\in (0,R \bar{y})}\frac{\rat(t)}{t}\geq \frac{1}{\bar{\gamma}}\left[\exp\left(\frac{d\bar{\gamma} \bar{r}_{d,\gamma}
      \sup_{(0,1)}\rat}{p}\right)-1\right], 
\]
so that for all $(y,\gamma)\in
(0,\bar{y})\times (0,\bar{\gamma})$,
\begin{align*} \ln\left(1+\gamma\frac{\rat(R y)}{R y}\right)\geq
  \frac{\gamma}{\bar{\gamma}}\ln\left(1+\bar{\gamma}\frac{\rat(R y)}{R y}\right)\geq
  \frac{d\bar{r}_{d,\gamma} \, \sup_{(0,1)} \rat }{p}\gamma\geq
  \frac{d}{p}\ln(1+\bar{r}_{d,\gamma}\gamma \sup_{(0,1)}\rat).
\end{align*}
where, for the first inequality, we used the fact that for $\alpha \in (0,1)$, for all $x \ge 0$,
$\ln(1+\alpha x) \ge \alpha \ln(1+x)$ (by concavity of the logarithm).
This concludes the proof of \eqref{surm}.

Now, to prove that
a.s. $\limsup_{k\to\infty}\underline{\tn}_{T_k-d}>0$, let us introduce
the stopping times $(\sigma_m)_{m \ge 0}$ and $(\tau_m)_{m \ge 1}$
such that $\sigma_0=0$, and for $m\geq 1$ (with the convention $\inf\emptyset=\infty$),
\[
\tau_m=\inf\{k>\sigma_{m-1} \, : \, Y_k\leq \bar{y}\}, \qquad
\sigma_m=\inf\{k>\tau_m \, : \, Y_k>\bar{y}\} \eqsp,
\]
where $\bar{y}$ has been introduced in~\eqref{surm}.
On the event $\{Y_k > \overline{y} \text{ infinitely often} \}$, one
  has $\limsup_{k\to\infty}Y_k \ge \overline{y}>0$. Notice that the
  complementary of the previous event writes $\{Y_k > \overline{y}
  \text{ infinitely often} \}^c=\{\exists m \ge 1, \, \tau_m < \infty
  = \sigma_m \}$.  To prove the result on this event, let us consider,
  for any fixed $m \ge 1$ and $l \ge 1$, the
  process $(Z_k)_{k \ge \underline{k} \vee l}$
  defined by
$$\forall k \ge \underline{k} \vee l, \, Z_k=-\ln (Y_{k\wedge
  \sigma_m}) \, \un_{\tau_m \le l}$$
where $\underline{k}$ has been introduced in~\eqref{surm}.
The process $(Z_k)_{k \ge \underline{k} \vee l}$ is a non-negative ${\mathcal F}_{T_k}$-supermartingale by \eqref{surm}
and thus converges a.s. to a finite limit as $k \to \infty$. Hence,
for any fixed $m \ge 1$, on
$\{\tau_m < \infty\}=\cup_{l \ge 1} \{ \tau_m \le l\}$, the process $(-\ln (Y_{k\wedge
  \sigma_m}) )_{k \ge 1}$ converges a.s. to a finite limit $V_m$. As a consequence, on $\{\exists m\geq
  1:\tau_m<\infty=\sigma_m\}$, $(Y_k)_{k \ge 1}$ converges a.s. to
  $\sum_{m\geq 1}
  \un_{\{\tau_m<\infty=\sigma_m\}}\mathrm{e}^{-V_m}$ which is positive on the
  event $\{\exists m\geq
  1:\tau_m<\infty=\sigma_m\}$. In conclusion,
  almost surely, $\limsup_{k\to\infty}
Y_k>0$. This concludes the proof.
\end{adem}

  \subsection{Proof of Proposition \ref{prop:pasetstab}}\label{sec:gamma_SHUS}
 We have checked in the previous section that the WL$_\rat$ algorithm
  (which encompasses the SHUS$_\rat^\alpha$ algorithm, see Section~\ref{sec:relation_FE}) is recurrent under mild conditions on the stepsize sequence
  $\left(\gamma_{n}\right)_{n \ge 1}$. In this section, we verify that for any $\alpha
\in (\frac12,1]$, these conditions are satisfied for the
  stepsize sequence generated by the SHUS$_\rat^\alpha$ algorithm, as
  well as the usual
  conditions~\eqref{hyp:Gamma} of summability on the sequence
  $\left(\gamma_{n}\right)_{n \ge 1}$ (when $\alpha=1$ and $\mu \in
  (0,1)$, this requires $\rho(t)=t^a$ for some $a\in [0,1)$). This is the content
  of Proposition~\ref{prop:stepsize} below.
  Proposition~\ref{prop:pasetstab} is then deduced from Propositions~\ref{prop:recurrence}
and~\ref{prop:stepsize}  by conditioning w.r.t.
  $\mathcal{F}_0$.

\begin{aprop} \label{prop:stepsize}
  Assume A\ref{hyp:targetpi}, A\ref{hyp:kernel}, R\ref{hyp:BoundFctRat},
  R\ref{hyp:fdec} and R\ref{hyp:finf}. The random stepsize sequence
  $(\gamma_{n+1}=\frac{\gamma}{g_\alpha(S_n)})_{n \ge 0}$ generated by the SHUS$_\rat^\alpha$ algorithm
  started from a deterministic initial condition
  $(\tu_0,X_0)\in(\R_+^*)^d\times\Xset$ is decreasing, bounded from above by
  some deterministic sequence converging to $0$ as $n\to\infty$ and such that
$\P\left(\inf_{n\geq 1} n^\alpha\gamma_n>0\right)=1$. Moreover,
$\bar{r}_{d,\gamma}=\sup_{n\geq
  1}\frac{\gamma_n}{\gamma_{n+d-1}}<\infty$ with the explicit upper bounds:
\begin{itemize}
\item if $\alpha\in(\frac{1}{2},1)$, $$\bar{r}_{d,\gamma} \le \left(1+
    \frac{(d-1) \gamma
      \sup_{t\in(0,1)}\rat(t)}{\ln(1+S_0)^{\frac{1}{1-\alpha}}}\right)^{\frac{\alpha}{1-\alpha}}
  \, ,$$
\item if $\alpha=1$,
$$\bar{r}_{d,\gamma} \le 1+
\frac{\mu\left(1\vee\left(1+\gamma\frac{\sup_{t\in(0,1)}
        \rat(t)}{S_0^\mu}\right)^{\mu-1}\right) (d-1) \gamma
  \sup_{t\in(0,1)}\rat(t)}{S^\mu_0} \, .$$
\end{itemize}
Finally,
\begin{itemize}
\item if $\alpha\in(\frac{1}{2},1)$, $\P\left(\sup_{n\geq 1}
n^\alpha\gamma_n<+\infty\right)=1$,
\item if $\alpha=1$ and $\mu\ge 1$, 
there exists a random variable $C$ such that $$\P\left(C>\frac{\mu}{1+\mu}\mbox{ and }\sup_{n\geq 1}
n^{C}\gamma_n<+\infty\right)=1,$$
\item if $\alpha=1$ and $\rat(t)=t^a$ for some $a\in[0,1)$, then $\P\left(\sup_{n\geq 1}
n\gamma_n<+\infty\right)=1$.
\end{itemize}
   \end{aprop}
The property $\inf_{n\geq 1} n^\alpha\gamma_n>0$ implies that $\sum_{n
  \ge 1} \gamma_n = \infty$, while the last three items provide 
sufficient assumptions to prove  $\sum_{n
  \ge 1} \gamma_n^2 < \infty$ (see  the assumptions~\eqref{hyp:Gamma}
required on the stepsize sequence to prove convergence). In particular, since $\frac{\mu}{1+\mu}\ge
\frac 12$ is equivalent to $\mu\ge 1$, the second item shows that $\sum_{n
  \ge 1} \gamma_n^2 < \infty$ when $\alpha=1$ and $\mu\ge 1$. When
$\alpha=1$ and $\mu<1$, we have not been able to prove that
$\P\left(\sum_{n\ge 1}\gamma_n^2<\infty\right)=1$ without supposing
that $\rat(t)=t^a$ for some $a\in [0,1)$.

Let us also mention that,  when $\alpha=1$ and $\rat(t)=t^a$ for some $a\in[0,1)$, the
proof we give below implies that $\P(\inf_{n \ge 0}\min_{1\le i\le d}\t_n(i)>0)=1$ (using
Equations~\eqref{majoSnmu} and~\eqref{eq:stab_ta} below),
{\em i.e.} that the SHUS$^1_{t^a}$ algorithm is stable. This gives
another way to prove the stability of the method in this specific setting, without following the
two-step argument that we used for a general $\rat$, namely first proving  the
recurrence of the algorithm (see Propositions~\ref{prop:recurrence} and~\ref{prop:pasetstab}), and then using~\cite[Theorem
2.2]{andrieu:moulines:priouret:2005} (see the proof of Proposition~\ref{prop:cvggal}).
\begin{adem}
We decompose the proof in several steps.

\noindent
{\bf Deterministic upper bound on $(\gamma_n)_{n \ge 1}$.} By \eqref{evosn1}, $(S_n)_{n\geq 0}$ is increasing so that
  $(\gamma_n)_{n\geq 1}$ is decreasing since $g_\alpha$ is increasing. Using~\eqref{evosn1} again, we have
  $$S_{n+1}=S_n+ \frac{\gamma}{g_\alpha(S_n)}
  \frac{\rat(\tn_n(I(X_{n+1})))}{\tn_n(I(X_{n+1}))}\tu_n(I(X_{n+1})),$$
  and
  since by \eqref{eq:def:thetatilde}, for all $i\in\{1,\hdots,d\}$,
  $(\tu_n(i))_{n\geq 0}$ is non-decreasing, it holds
\begin{equation}\label{eq:UpDownSn}
S_n + \frac{ \gamma}{g_\alpha(S_n)} \left( \inf_{t \in (0,1)} \frac{\rat(t)}{t} \right)
\min_{i=1,\cdots, d} \tu_0(i) \leq S_{n+1} \leq S_n \left(1 + \gamma \,
  \frac{\sup_{t\in(0,1)} \rat(t)}{g_\alpha(S_n)} \right).
\end{equation}
  The lower bound on $S_{n+1}$ implies that the sequence $(\gamma_n)_{n \ge 1}$ is bounded from above by a deterministic
    sequence converging to $0$ as $n\to\infty$, following the
    arguments in the proof of \cite[Lemma
    1]{fort:jourdain:lelievre:stoltz:2015}.
  
\noindent
{\bf Lower bound on $(\gamma_n)_{n \ge 0}$.} When $\alpha \in
  (\frac12,1)$, one obtains the lower bound on
  $(n^\alpha \gamma_n)_{n\ge 1}$  from the upper bound in \eqref{eq:UpDownSn} by an easy adaptation of the proof of~\cite[Lemma
  1]{fort:jourdain:lelievre:stoltz:2015}. When $\alpha=1$, since $g_1(s)=s^\mu$ with $\mu>0$, the upper bound in \eqref{eq:UpDownSn} and the inequality $(1+x)^\mu\le 1+\mu(1\vee (1+x_0)^{\mu-1}) x$ for $0\le x\le x_0$ imply that
\begin{equation}
   S_{n+1}^\mu\le S_n^\mu\left(1+\gamma\frac{\sup_{t\in(0,1)} \rat(t)}{S_n^\mu}\right)^\mu\le S_n^\mu+\tilde\mu\gamma\sup_{t\in(0,1)} \rat(t).\label{majodifsmun}
\end{equation}
with $\tilde{\mu}=\mu\left(1\vee\left(1+\gamma\frac{\sup_{t\in(0,1)} \rat(t)}{S_0^\mu}\right)^{\mu-1}\right)$.
By induction on $n$, one deduces that, for all $n\in{\mathbb N}$, 
\begin{equation}
   S_n^\mu\le S_0^\mu+n\tilde\mu\gamma\sup_{t\in(0,1)} \rat(t) \, ,\label{majoSnmu}
\end{equation} and therefore 
$$\gamma_{n+1}=\frac{\gamma}{S_n^\mu}\ge \frac{\gamma}{S_0^\mu+n\tilde\mu\gamma\sup_{t\in(0,1)} \rat(t)}=\frac{\gamma_1}{1+n\tilde\mu\gamma_1\sup_{t\in(0,1)} \rat(t)}.$$

\noindent
{\bf Upper bounds on $\bar{r}_{d,\gamma}$.} Let us now derive the upper
bounds on $\bar{r}_{d,\gamma}$.
  When $\alpha=1$, by \eqref{majodifsmun}, $S^\mu_{n+d-1}\leq
  S^\mu_n+\tilde \mu \gamma (d-1) \sup_{t\in(0,1)}\rat(t)$, which yields
\[ 
\frac{\gamma_{n+1}}{\gamma_{n+d}}=\frac{S^\mu_{n+d-1}}{S^\mu_n}
\leq 1+\frac{\tilde \mu \gamma (d-1) \sup_{t\in(0,1)}\rat(t)}{S^\mu_n}  
\leq 1+\frac{\tilde \mu \gamma (d-1) \sup_{t\in(0,1)}\rat(t)}{S^\mu_0}.  
\] 
When $\alpha\in(\frac{1}{2},1)$, by
\eqref{eq:UpDownSn}, the inequality $\ln(1+x) \leq x$ on $\R^+$ and the monotonicity
of $(S_n)_{n \geq 0}$, we have for any $0 \leq q \leq n$,
\begin{align*}
  \ln(1+S_{n+1}) &=\ln(1+S_n)+\ln\left(1+\frac{\gamma S_n
      \, \sup_{t\in(0,1)}\rat(t)}{(1+S_n)\ln(1+S_n)^{\frac{\alpha}{1-\alpha}}}\right) 
  \leq\ln(1+S_n)+\frac{\gamma \sup_{t\in(0,1)}\rat(t)
  }{\ln(1+S_0)^{\frac{\alpha}{1-\alpha}}}\\&\leq\ln(1+S_n)+\ln(1+S_q)\frac{\gamma \sup_{t\in(0,1)}\rat(t)
  }{\ln(1+S_0)^{\frac{1}{1-\alpha}}}.
\end{align*}
Therefore,
\[
\frac{\gamma_{n+1}}{\gamma_{n+d}}=\frac{\ln(1+S_{n+d-1})^{\frac{\alpha}{1-\alpha}}}{\ln(1+S_n)^{\frac{\alpha}{1-\alpha}}}\leq
\left(1+\frac{\gamma(d-1) \sup_{t\in(0,1)} \rat(t)}{\ln(1+S_0)^{\frac{1}{1-\alpha}}}\right)^{\frac{\alpha}{1-\alpha}}.
\]

\noindent
{\bf Upper bounds on $(\gamma_n)_{n \ge 1}$: the two cases $\alpha \in
  (\frac12,1)$ or $\alpha=1$ and $\mu \ge 1$.}
To deal with the last assertion, we are going to derive lower bounds on  $$\underline{\tu}_n=\min_{1\le i\le d} \tu_n(i).$$
By~\eqref{eq:def:thetatilde}, for all $n\ge 0$ and
  for all $i \in \{1, \ldots, d\}$,
\begin{equation}\label{eq:rectu}
  \tu_{n+1}(i) = \tu_n(i)\left( 1 +  \gamma_{n+1}
  \, \frac{\rat(\t_n(i))}{\t_n(i)}\,  \un_{\Xset_i}(X_{n+1}) \right).
\end{equation}
As a consequence, for all $k\ge 0$ and
  for all $i \in \{1, \ldots, d\}$,
$$\tu_{T_{k+1}}(i) \ge \tu_{T_{k+1}-d}(i)
\prod_{n=T_{k+1}-d+1}^{T_{k+1}}\left( 1 +  \gamma_{n+1} \inf_{t \in (0,1)} \frac{\rat(t)}{t}  \un_{\Xset_i}(X_{n}) \right).
$$
Now, if $i$ is the index of a stratum with large weight, namely
$\tu_{T_{k+1}-d}(i) \ge R \underline{\tu}_{T_{k+1}-d}$, one
simply uses the lower bound:
\begin{equation}\label{eq:l1}
\tu_{T_{k+1}}(i) \ge R \underline{\tu}_{T_{k+1}-d}.
\end{equation}
If  $i$ is the index of a stratum with small weight, namely
$\tu_{T_{k+1}-d}(i) < R \underline{\tu}_{T_{k+1}-d}$, by
definition of the sequence $(T_k)_{k\in\N}$, this stratum is visited
at least once between $T_{k+1}-d$ and $T_{k+1}$, and thus, using the monotonicity of the sequence $(\gamma_n)_{n\ge 1}$, we get
\begin{equation}\label{eq:l2}
\tu_{T_{k+1}}(i) \ge \left( 1 +  \gamma_{T_{k+1}} \inf_{t \in (0,1)} \frac{\rat(t)}{t}  \right) \underline{\tu}_{T_{k+1}-d}.
\end{equation}
By combining~\eqref{eq:l1} and~\eqref{eq:l2}, one thus obtains
$$\underline{\tu}_{T_{k+1}}\geq \left(R \wedge
    \left(1+ \gamma_{T_{k+1}} \inf_{t \in (0,1)} \frac{\rat(t)}{t} 
     \right)\right)\underline{\tu}_{T_{k+1}-d}.$$
Since $T_k \le T_{k+1}-d$, by the monotonicity of the sequence $(\tu_n)_{n
  \ge 0}$, one concludes that 
\begin{equation*}
   \forall k\in\N,\;\underline{\tu}_{T_{k+1}}\geq \left(R \wedge
    \left(1+\gamma_{T_{k+1}}\inf_{t \in (0,1)} \frac{\rat(t)}{t} \
      \right)\right)\underline{\tu}_{T_{k}}.
  \end{equation*}
Since the sequence $(\gamma_n)_{n \ge 1}$ is bounded from
above by a deterministic sequence converging to $0$, one deduces that
there exists $K$ such that
\begin{equation}\label{evolmintuentredeuxtk}
   \forall k \ge K,\;\underline{\tu}_{T_{k+1}}\geq 
    \left(1+\gamma_{T_{k+1}}\inf_{t \in (0,1)} \frac{\rat(t)}{t} \
      \right)\underline{\tu}_{T_{k}}.
  \end{equation}
  This inequality together with \eqref{majotk}, the lower bound on $(n^\alpha\gamma_n)_{n\ge 1}$ and 
\begin{equation}
   S_{n+1}\ge S_n + \frac{ \gamma}{g_\alpha(S_n)} \left( \inf_{t \in (0,1)} \frac{\rat(t)}{t} \right)
\underline{\tu}_n\label{minodiffsn}
\end{equation} permits to use the arguments of the proof of \cite[Lemma
  1]{fort:jourdain:lelievre:stoltz:2015} to get $\P\left(\sup_{n\geq 1}
n^\alpha\gamma_n<+\infty\right)=1$ when
  $\alpha\in(\frac12,1)$. 

When $\alpha=1$, \eqref{minodiffsn} together with the inequality $(x+y)^{1+\mu}\ge x^{1+\mu}+(1+\mu)x^\mu y$ valid for $x,y>0$ leads to
$$S_{n+1}^{1+\mu}\ge S_n^{1+\mu}+(1+\mu)\gamma\left( \inf_{t \in (0,1)} \frac{\rat(t)}{t} \right)
\underline{\tu}_n.$$
Together with \eqref{majotk} and \eqref{evolmintuentredeuxtk}, this inequality permits to adapt the arguments of the proof of \cite[Proposition 
  1]{fort:jourdain:lelievre:stoltz:2015} to obtain the existence of a
  random variable $\varepsilon > 0$ such that a.s. $\inf_{n \ge 1}
  n^{-1-\varepsilon} S_n^{1+\mu} >0$. Therefore, there exists a random
  variable $C$ such that $\P\left(C>\frac{\mu}{1+\mu}\mbox{ and }\sup_{n\geq 1}
n^{C}\gamma_n<+\infty\right)=1$. 

\noindent
{\bf Upper bounds on $(\gamma_n)_{n \ge 1}$: the case $\alpha=1$ and $\rho(t)=t^a$.}
Let us now suppose that $\alpha=1$ and $\rat(t)=t^a$ for some $a\in
[0,1)$, so that $\gamma_{n+1}=\frac{\gamma}{S_n^\mu}$. The objective
is to show that a.s. $\sup_{n \ge 1} n \gamma_n < \infty$. Since
$S_n=\sum_{i=1}^d\tu_{n}(i)\ge \underline{\tu}_{n}$, it is
sufficient to prove that
\begin{equation}\label{eq:stab_ta}
\P\left(\inf_{n \ge 1}n^{-\frac{1}{\mu}}\underline{\tu}_{n}>0\right)=1.
\end{equation}

 Writing $\tu_n(i)\gamma_{n+1}\frac{\rat(\t_n(i))}{\t_n(i)}=
\gamma \tu_n(i)^a S_n^{1-(a+\mu)}$ in~\eqref{eq:rectu} and using the
monotonicity of the sequences $(\tu_n(i))_{n\ge 0}$ and $(S_n)_{n\ge
  0}$, one deduces that~\eqref{eq:l2} and~\eqref{evolmintuentredeuxtk}
may respectively
be replaced by:
\begin{itemize}
\item  For all $k\in\N$, if  $i$ is the index of a stratum with small weight, namely
$\tu_{T_{k+1}-d}(i) < R \underline{\tu}_{T_{k+1}-d}$,
\begin{equation}\label{eq:l2_ta}
\tu_{T_{k+1}}(i) \ge
\tu_{T_{k+1}-d}(i)+\gamma\tu^a_{T_{k+1}-d}(i)\left(S_{T_{k+1}-d}^{1-(\mu+a)}\wedge
  S_{T_{k+1}-1}^{1-(\mu+a)}\right).
\end{equation}
\item By combining~\eqref{eq:l1} and~\eqref{eq:l2_ta}, one thus obtains
\begin{equation}\label{evolmintuentredeuxtkratpuis}
\underline{\tu}_{T_{k+1}}\geq \left(R\underline{\tu}_{T_{k}}\right) \wedge \left(\underline{\tu}_{T_{k}}+\gamma\underline{\tu}^a_{T_{k}}(S_{T_{k}}^{1-(\mu+a)}\wedge
  S_{T_{k+1}-1}^{1-(\mu+a)})\right)
\end{equation} 
\end{itemize}

Let us suppose that $\mu+a\ge 1$. Since by~\eqref{evolmintuentredeuxtkratpuis}, $k\mapsto\underline{\tu}_{T_{k}}$ grows at least geometrically with ratio $R$ as
long as the decreasing sequence $\left(1+\gamma \left(\underline{\tu}_{T_{k}}\right)^{a-1}S_{T_{k+1}-1}^{1-(\mu+a)}\right)_{k
\ge 1}$ is larger than $R$, there exists a random variable $K\in \N
\setminus \{0\}$ such that a.s. 
$$\forall k\ge K,\;\underline{\tu}_{T_{k+1}}\geq \underline{\tu}_{T_{k}}+\gamma\underline{\tu}^a_{T_{k}}S_{T_{k+1}-1}^{1-(\mu+a)}.$$

Since by \eqref{majoSnmu} and \eqref{majotk}, there exists a positive
random variable $C$ such that a.s. $\forall k,\;\gamma
S_{T_{k+1}-1}^{1-(\mu+a)}\ge Ck^{\frac{1-(\mu+a)}{\mu}}$ we deduce that a.s.,
\begin{align*}
   \forall k\ge K,\;\underline{\tu}^{1-a}_{T_{k+1}}\geq \underline{\tu}^{1-a}_{T_{k}}\left(1+C\left(\underline{\tu}_{T_{k}}\right)^{a-1}k^{\frac{1-(\mu+a)}{\mu}}\right)^{1-a}.
\end{align*}
With the monotonicity of $(\underline{\tu}_n)_{n \ge 0}$ and
the inequality $(1+x)^{1-a}\ge 1+(1-a)(1+x_0)^{-a}x$ valid for $0\le
x\le x_0$,  we deduce that a.s.
$$\forall k\ge K,\;\underline{\tu}^{1-a}_{T_{k+1}}\geq
\underline{\tu}^{1-a}_{T_{k}}+(1-a)\left(1+C\left(\underline{\tu}_{0}\right)^{a-1}\right)^{-a} C k^{\frac{1-(\mu+a)}{\mu}}.$$
Therefore $\P(\inf_{k\ge 1}k^{\frac{a-1}{\mu}}\underline{\tu}^{1-a}_{T_{k}}>0)=1$. Since the inequality $\forall
n,\;T_{\lfloor n/C_\star\rfloor}\le n$ deduced from \eqref{majotk}
and the monotonicity of $(\underline{\tu}_n)_{n \ge 0}$ imply
that $\underline{\tu}_n\ge\underline{\tu}_{T_{\lfloor
    n/C_\star\rfloor}}$, we deduce that $\P\left(\inf_{n \ge 1}n^{-\frac{1}{\mu}}\underline{\tu}_{n}>0\right)=1$ which gives~\eqref{eq:stab_ta}.

Let us finally suppose that $a+\mu<1$. Since $S_{T_{k+1}}^{1-(\mu+a)}\ge S_{T_{k}}^{1-(\mu+a)}\ge(\underline{\tu}_{T_{k}})^{1-(\mu+a)}$, \eqref{evolmintuentredeuxtkratpuis} implies 
$$\underline{\tu}_{T_{k+1}}\geq \left(R\underline{\tu}_{T_{k}}\right) \wedge \left(\underline{\tu}_{T_{k}}+\gamma\left(\underline{\tu}_{T_{k}}\right)^{1-\mu}\right).$$
Reasoning like in the case $\mu+a\ge 1$, we obtain the existence of a random variable $K<\infty$ such that a.s. 
$$\forall k\ge K,\;\underline{\tu}_{T_{k+1}}\geq \underline{\tu}_{T_{k}}+\gamma\left(\underline{\tu}_{T_{k}}\right)^{1-\mu},$$
and deduce that $\P\left(\inf_{n\ge 1} n^{-\frac{1}{\mu}}\underline{\tu}_{n}>0\right)=1$.



\end{adem}


\subsection{Proof of Corollary~\ref{cor:asymptotic}}

\paragraph{Case $\alpha=1$.} Using~\eqref{evosn1} together with $\gamma_{n+1}=\frac{\gamma}{S_n^\mu}$, one obtains that for $k\ge 1$, 
$$S_{k}^\mu=S_{k-1}^\mu\left(1+\gamma\frac{\rho(\t_{k-1}(I(X_{k})))}{S_{k-1}^\mu}\right)^\mu.$$
Since for $0\le x\le x_0$, $$|(1+x)^\mu-(1+\mu
x)|=\left|\mu(\mu-1)\int_0^x(x-y)(1+y)^{\mu-2} \rmd y\right|\le \frac{\mu|\mu-1|}{2}(1\vee (1+x_0)^{\mu-2}) x^2,$$ one deduces that for $n\in{\mathbb N}$, 
$$\frac{S_{n}^\mu}{n}=\frac{S_0^\mu}{n}+\frac{\gamma\mu}{n}\sum_{k=1}^n\rho(\t_{k-1}(I(X_{k})))-\frac{1}{n}\sum_{k=1}^nR_{k},$$
where $0\le R_{k}\le \frac{\mu|\mu -1|}{2}\left(1\vee \left(1+\gamma\frac{\sup_{t\in (0,1)}\rho(t)}{S_0^\mu}\right)^{\mu-2}\right) \gamma\left(\sup_{t\in (0,1)}\rat(t)\right)^2\gamma_{k}$.
By R\ref{hyp:BoundFctRat}, Proposition \ref{prop:stepsize} and
Ces\`aro Lemma, it holds: a.s., $\lim_{n\to \infty}\frac{1}{n}\sum_{k=1}^n R_k=0$. Therefore, to conclude the proof it is enough to check that a.s.,
  \begin{equation}\label{eq:ave_theta}
\lim_{n \to \infty} \frac{1}{n}\sum_{k=1}^n
\rat(\theta_{k-1}(I(X_k))) = \left(\sum_{j=1}^d\frac{\theta_\star(j)}{\rat(\theta_\star(j))} \right)^{-1}.
\end{equation}
This follows from the decomposition 
\begin{align*}
 \left( \sum_{j=1}^d
  \frac{\theta_{\star}(j)}{\rat(\theta_{\star}(j))}\right) \frac{1}{n} &\sum_{k=1}^n 
  \rat(\theta_{k-1}(I(X_k)))=\frac{1}{n} \sum_{k=1}^n \left( \sum_{j=1}^d
  \frac{\theta_{k-1}(j)}{\rat(\theta_{k-1}(j))} \right)
  \rat(\theta_{k-1}(I(X_k)))  \\&+\frac{1}{n} \sum_{k=1}^n \left(\sum_{j=1}^d
  \frac{\theta_{\star}(j)}{\rat(\theta_{\star}(j))}- \sum_{j=1}^d
  \frac{\theta_{k-1}(j)}{\rat(\theta_{k-1}(j))}   \right)
  \rat(\theta_{k-1}(I(X_k))).\end{align*}
where the first term in the right-hand side almost surely converges to $1$ by item
(iii) in Theorem~\ref{theo:cvggal} (choosing $f\equiv 1$). The absolute value of the second term in the right-hand side is bounded from above by
$$\left(\sup_{t\in(0,1)} \rat(t)\right)\frac{1}{n} \sum_{k=1}^n \left| \sum_{j=1}^d
  \left\{\frac{\theta_{\star}(j)}{\rat(\theta_{\star}(j))}-\frac{\theta_{k-1}(j)}{\rat(\theta_{k-1}(j))} 
     \right\} \right|$$
which
almost surely goes to zero thanks to item (i) in Theorem~\ref{theo:cvggal}, the continuity of $\rat$ deduced from R\ref{hyp:LypFctRat} and
Ces\`aro Lemma. 

\paragraph{Case $\alpha \in(\frac12,1)$.} For all $n \ge 0$,
$\gamma_{n+1}=\frac{\gamma}{g_\alpha(S_n)}=\frac{\gamma}{\ln(1+S_n)^{\frac{\alpha}{1-\alpha}}}$.
Following the proof of~\cite[Proposition
4]{fort:jourdain:lelievre:stoltz:2015}, using~\eqref{evosn1}, R\ref{hyp:BoundFctRat} and the fact that $\P(\sup_nn^\alpha\gamma_n<\infty)=1$ deduced from Proposition \ref{prop:stepsize}, one can prove that for all $n
\ge 1$,
$$\frac{1}{n} \left( \ln(1+S_n)\right)^{\frac{1}{1-\alpha}} = 
\frac{1}{n} \left( \ln(1+S_0)\right)^{\frac{1}{1-\alpha}} + \frac{\gamma}{(1-\alpha)n} \sum_{k=1}^n
\rat(\theta_{k-1}(I(X_k)))  +
\frac{1}{n} \sum_{k=1}^n R_k,$$
where for all $k \ge 1$, $|R_k| \le C k^{-\alpha}$, for a finite
random variable $C$. The first and last terms almost surely converge
to zero. Using~\eqref{eq:ave_theta}, one thus obtains, almost surely,
$$\lim_{n \to \infty} \frac{1}{n} \left( \ln(1+S_n)\right)^{\frac{1}{1-\alpha}}
= \frac{\gamma}{1-\alpha} \left( \sum_{i=1}^d
  \frac{\theta_\star(i)}{\rat(\theta_\star(i))} \right)^{-1}$$
which concludes
the proof.

\subsection*{Acknowledgements}

This work is supported by the European Research Council under the
European Union's Seventh Framework Programme (FP/2007-2013) / ERC
Grant Agreement number 614492 and by the French National Research
Agency under the grant ANR-14-CE23-0012 (COSMOS). We also benefited
from the scientific environment of the Laboratoire International
Associ\'e between the Centre National de la Recherche Scientifique and
the University of Illinois at Urbana-Champaign. We would like to thank
Brad Dickson, Alessandro Laio and Michele Parrinello for useful discussions.



\begin{thebibliography}{10}

\bibitem{andrieu:moulines:priouret:2005}
C.~Andrieu, E.~Moulines, and P.~Priouret.
\newblock Stability of stochastic approximation under verifiable conditions.
\newblock {\em SIAM J. Control Optim.}, 44:283--312, 2005.

\bibitem{atchade-liu-10}
Y.F. Atchad\'e and J.S. Liu.
\newblock The {Wang-Landau} algorithm for {Monte Carlo} computation in general
  state spaces.
\newblock {\em Stat. Sinica}, 20(1):209--233, 2010.


\bibitem{barducci-bussi-parrinello-08}
A.~Barducci, G.~Bussi, and M.~Parrinello.
\newblock Well-tempered metadynamics: A smoothly converging and tunable
  free-energy method.
\newblock {\em Phys. Rev. Lett.}, 100:020603, 2008.

\bibitem{benveniste:metivier:priouret:1987}
A.~Benveniste, M.~Metivier, and P.~Priouret.
\newblock {\em Adaptive {A}lgorithms and {S}tochastic {A}pproximations}.
\newblock Springer-Verlag, 1987.

\bibitem{borkar:2008}
V.S. Borkar.
\newblock {\em Stochastic Approximation: A Dynamical Systems Viewpoint}.
\newblock Cambridge University Press, 2008.

\bibitem{brooks:gelman:jones:meng:2011}
S.~Brooks, A.~Gelman, G.L. Jones, and X-L. Meng.
\newblock {\em {Handbook of Markov Chain Monte Carlo}}.
\newblock Chapman \& Hall, 2011.

\bibitem{bussi-laio-parrinello-06}
G.~Bussi, A.~Laio, and M.~Parrinello.
\newblock Equilibrium free energies from nonequilibrium metadynamics.
\newblock {\em Phys. Rev. Lett.}, 96:090601, 2006.

\bibitem{chen:2002}
H.~Chen.
\newblock {\em Stochastic {A}pproximation and {I}ts {A}pplications}.
\newblock Kluwer Academic Publishers, 2002.

\bibitem{chipot-pohorille-07}
C.~Chipot and A.~Pohorille, editors.
\newblock {\em Free Energy Calculations}, volume~86 of {\em Springer Series in
  Chemical Physics}.
\newblock Springer, 2007.

\bibitem{chopin-lelievre-stoltz-12}
N.~Chopin, T.~Leli\`{e}vre, and G.~Stoltz.
\newblock Free energy methods for {B}ayesian inference: efficient exploration
  of univariate {G}aussian mixture posteriors.
\newblock {\em Stat. Comput.}, 22(4):897--916, 2012.

\bibitem{crespo2010metadynamics}
Y.~Crespo, F.~Marinelli, F.~Pietrucci, and A.~Laio.
\newblock Metadynamics convergence law in a multidimensional system.
\newblock {\em Phys. Rev. E}, 81(5):055701, 2010.

\bibitem{DHSV15}
J.F. Dama, G.M. Hocky, R. Sun and G.A. Voth.
\newblock Exploring valleys without climbing every peak: More efficient and forgiving metabasin metadynamics via robust on-the-fly bias domain restriction.
\newblock {\em J. Chem. Theory Comput.}, 11(12):5638-5650, 2015


\bibitem{dama-parrinello-voth-14}
J.F. Dama, M.~Parrinello, and G.A. Voth.
\newblock Well-tempered metadynamics converges asymptotically.
\newblock {\em Phys. Rev. Lett.}, 112:240602(1--6), 2014.

\bibitem{DP01}
E.~Darve and A.~Pohorille.
\newblock Calculating free energies using average force.
\newblock {\em J. Chem. Phys.}, 115(20):9169--9183, 2001.

\bibitem{Dickson15}
B. Dickson.
\newblock $\mu$-tempered metadynamics: Artifact independent convergence times for wide hills.
\newblock {\em J. Chem. Phys.}, 143(23):234109, 2015.


\bibitem{amrx}
G.~Fort, B.~Jourdain, E.~Kuhn, T.~Leli\`evre, and G.~Stoltz.
\newblock Efficiency of the {W}ang-{L}andau algorithm: A simple test case.
\newblock {\em Appl. Math. Res. Express}, 2014(2):275--311, 2014.

\bibitem{fort-2015}
G.~Fort.
\newblock Central limit theorems for stochastic approximation with controlled {M}arkov chain dynamics.
\newblock {\em ESAIM: PS}, 19:60--80, 2015.

\bibitem{fort:jourdain:kuhn:lelievre:stoltz:2014}
G.~Fort, B.~Jourdain, E.~Kuhn, T.~Leli\`evre, and G.~Stoltz.
\newblock Convergence of the {Wang-Landau} algorithm.
\newblock {\em Math. Comput.}, 84(295):2297--2327, 2015.

\bibitem{fort:jourdain:lelievre:stoltz:2015}
G.~Fort, B.~Jourdain, T.~Leli\`evre, and G.~Stoltz.
\newblock {Self-Healing Umbrella Sampling: Convergence and efficiency}.
\newblock {\em Stat Comput.}, 27(1), 147-168, 2017. 

\bibitem{fort:moulines:priouret:2011}
G.~Fort, E.~Moulines, and P.~Priouret.
\newblock Convergence of adaptive and interacting {M}arkov chain {M}onte
  {C}arlo algorithms.
\newblock {\em Ann. Statist.}, 39(6):3262--3289, 2012.


\bibitem{fort:moulines:schreck:vihola:2015}
G.~Fort, E.~Moulines, A.~Schreck and M.~Vihola.
\newblock Convergence of Markovian Stochastic Approx-
imation with discontinuous dynamics. 
\newblock {\em SIAM J. Control Optim.}, 54(2):866--893, 2016.


\bibitem{hall:heyde:1980}
P.~Hall and P.P. Heyde.
\newblock {\em Martingale Limit Theory and its application}.
\newblock Academic Press, 1980.

\bibitem{Hastings70}
W.K. Hastings.
\newblock Monte {C}arlo sampling methods using {M}arkov chains and their
  applications.
\newblock {\em Biometrika}, 57:97--109, 1970.

\bibitem{HC04}
J.~H\'enin and C.~Chipot.
\newblock Overcoming free energy barriers using unconstrained molecular
  dynamics simulations.
\newblock {\em J. Chem. Phys.}, 121(7):2904--2914, 2004.

\bibitem{JR11}
P.E. Jacob and R.J. Ryder.
\newblock The {Wang-Landau} algorithm reaches the flat histogram criterion in
  finite time.
\newblock {\em Ann. Appl. Probab.}, 24(1):34--53, 2014.

\bibitem{JLR10}
B.~Jourdain, T.~Leli{\`e}vre, and R.~Roux.
\newblock Existence, uniqueness and convergence of a particle approximation for
  the adaptive biasing force process.
\newblock {\em ESAIM: M2AN}, 44(5):831--865, 2010.

\bibitem{KongLiuWong}
A.~Kong, J.~S. Liu, and W.H. Wong.
\newblock Sequential imputation and {B}ayesian missing data problems.
\newblock {\em J. Am. Statist. Assoc.}, 89:278--288, 1994.

\bibitem{kushner:2010}
H.~Kushner.
\newblock Stochastic approximation: a survey.
\newblock {\em Wiley Interdisciplinary Reviews: Computational Statistics},
  2(1):87--96, 2010.

\bibitem{kushner:yin:2003}
H.J. Kushner and G.G. Yin.
\newblock {\em Stochastic Approximation and Recursive Algorithms and Applications}.
\newblock Springer, 2003.

\bibitem{laio-parrinello-02}
A.~Laio and M.~Parrinello.
\newblock Escaping free-energy minima.
\newblock {\em Proc. Natl. Acad. Sci. U.S.A}, 99:12562--12566, 2002.

\bibitem{lelievre-minoukadeh-09}
T.~Leli\`evre and K.~Minoukadeh.
\newblock Long-time convergence of an adaptive biasing force method: {T}he
  bi-channel case.
\newblock {\em Arch. Ration. Mech. Anal.}, 202(1):1--34, 2011.

\bibitem{LRS08}
T.~Leli\`evre, M.~Rousset, and G.~Stoltz.
\newblock Long-time convergence of an adaptive biasing force method.
\newblock {\em Nonlinearity}, 21:1155--1181, 2008.

\bibitem{lelievre-rousset-stoltz-book-10}
T.~Leli\`evre, M.~Rousset, and G.~Stoltz.
\newblock {\em Free Energy Computations: A Mathematical Perspective}.
\newblock Imperial College Press, 2010.

\bibitem{SHUS}
S.~Marsili, A.~Barducci, R.~Chelli, P.~Procacci, and V.~Schettino.
\newblock Self-healing {Umbrella Sampling}: {A} non-equilibrium approach for
  quantitative free energy calculations.
\newblock {\em J. Phys. Chem. B}, 110(29):14011--14013, 2006.

\bibitem{MVTP15}
J. McCarty, O. Valsson, P. Tiwary and M. Parrinello.
\newblock Variationally optimized free-energy flooding for rate calculation,
\newblock {\em Phys. Rev. Lett.}, 115(7):070601, 2015

\bibitem{mcgovern2013boundary}
M.~McGovern and J.~de~Pablo.
\newblock A boundary correction algorithm for metadynamics in multiple
dimensions.
\newblock {\em J. Chem. Phys.}, 139(8):084102, 2013.

\bibitem{MRRTT53}
N.~Metropolis, A.W. Rosenbluth, M.N. Rosenbluth, A.H. Teller, and E.~Teller.
\newblock Equations of state calculations by fast computing machines.
\newblock {\em J. Chem. Phys.}, 21(6):1087--1091, 1953.

\bibitem{metzner-schuette-vanden-eijnden-06}
P.~Metzner, C.~Sch\"{u}tte, and E.~Vanden-Eijnden.
\newblock Illustration of transition path theory on a collection of simple
  examples.
\newblock {\em J. Chem. Phys.}, 125(1), 2006.

\bibitem{meyn:tweedie:2009}
S.~Meyn and R.L. Tweedie.
\newblock {\em Markov Chains and {S}tochastic {S}tability}.
\newblock Cambridge, 2009.

\bibitem{minoukadeh-chipot-lelievre-10}
K.~Minoukadeh, C.~Chipot, and T.~Leli\`evre.
\newblock Potential of mean force calculations: a multiple-walker adaptive
  biasing force approach.
\newblock {\em J. Chem. Th. Comput.}, 6(4):1008--1017, 2010.

\bibitem{PSLS03}
S.~Park, M.K. Sener, D.~Lu, and K.~Schulten.
\newblock Reaction paths based on mean first-passage times.
\newblock {\em J. Chem. Phys.}, 119(3):1313--1319, 2003.

\bibitem{PJ92}
B.T.~Polyak, and A.B. Juditsky.
\newblock Acceleration of stochastic approximation by averaging.
\newblock {\em SIAM J. Control Optim.}, 30(4):838--855,  1992.

\bibitem{R88}
D.~Ruppert.
\newblock Efficient estimations from a slowly convergent Robbins-Monro
process.
\newblock Technical Report 781, Cornell University Operations Research and Industrial Engineering, 1988.

\bibitem{wang-landau-01}
F. Wang and D.P. Landau.
\newblock Determining the density of states for classical statistical models: A
  random walk algorithm to produce a flat histogram.
\newblock {\em Phys. Rev. E}, 64:056101, 2001.

\bibitem{wang-landau-01-PRL}
F. Wang and D.P. Landau.
\newblock Efficient, multiple-range random walk algorithm to calculate the
  density of states.
\newblock {\em Phys. Rev. Lett.}, 86(10):2050--2053, 2001.

\end{thebibliography}
\end{document}